\documentclass[a4paper,11pt,reqno]{amsart}

\usepackage{amsmath,amsthm,amssymb,latexsym}
\newtheorem{thm}{Theorem}[section]

\newtheorem{prop}{Proposition}[section]
\newtheorem{lem}{Lemma}[section]
\theoremstyle{definition}

\newtheorem{rem}{Remark}[section]

\numberwithin{equation}{section}
\newcommand{\pt}{\partial}
\newcommand{\rre}{{\mathbb R}}
\newcommand{\cmx}{{\mathbb C}}

\newcommand{\zzz}{{\mathbb Z}}

\newcommand{\bmat}{\left(\begin{array}}
\newcommand{\emat}{\end{array} \right)}
\title[On Quantum Zakharov System]{
The Fourth Order Nonlinear Schr\"{o}dinger\\
limit for Quantum Zakharov System}

\author[Fang, Lin, and Segata]{Yung-Fu Fang${}^{\dagger}$, Chi-Kun Lin, and
Jun-ichi Segata${}^{\ddagger}$}

\date{}

\begin{document}

\maketitle

\vskip-5mm \centerline{${}^{\dagger}$Department of Mathematics
National Cheng Kung University} \centerline{No. 1, Dasyue Rd., Tainan
City 70101, Taiwan}

\vskip3mm \centerline{Department of Mathematical Sciences}
\centerline{Xi'an Jiaotong-Liverpool University}
\centerline{SIP. Suzhou, Jiansu 215123 P.R. China}

\vskip3mm \centerline{${}^{\ddagger}$ Mathematical Institute, Tohoku
University} \centerline{ 6-3, Aoba, Aramaki, Aoba-ku, Sendai 980-8578,
Japan}

\vskip5mm

\begin{abstract}
This paper is concerned with the quantum Zakharov system. We prove
that when the ionic speed of sound goes to infinity, the solution to the
fourth order Schr\"{o}dinger part of the quantum Zakharov system
converges to the solution to quantum modified nonlinear
Schr\"{o}dinger eqaution.


\end{abstract}

%----------------------%     Section 1.   Introduction     %------------------------------
\section{\bf Introduction} \label{Sec1}

We consider the quantum Zakharov system:
\begin{equation}
%\begin{eqnarray}
\left\{
\begin{array}{l}
\medskip
i\pt_{t}E_{\varepsilon,\lambda}
+\Delta E_{\varepsilon,\lambda}
-\varepsilon^{2}\Delta^{2}
E_{\varepsilon,\lambda}
=n_{\varepsilon,\lambda}E_{\varepsilon,\lambda}
\qquad t\in\rre,\ x\in\rre^{d},\\
\medskip
\displaystyle{
\lambda^{-2}\pt_{t}^{2}n_{\varepsilon,\lambda}
-\Delta n_{\varepsilon,\lambda}
+\varepsilon^{2}\Delta^{2}n_{\varepsilon,\lambda}
=\Delta|E_{\varepsilon,\lambda}|^{2}
\qquad t\in\rre,\ x\in\rre^{d},}\\
\displaystyle{
E_{\varepsilon,\lambda}(0,x)=E_{0}(x),
\ n_{\varepsilon,\lambda}(0,x)=n_{0}(x),
\ \pt_{t}n_{\varepsilon,\lambda}(0,x)
=n_{1}(x)\qquad x\in\rre^{d},}
\end{array}
\right. \label{ZK}
%\end{eqnarray}
\end{equation}
where $\varepsilon\in(0,1],\lambda\in[1,\infty)$,
$E_{\varepsilon,\lambda}:\rre\times\rre^{d}\to\cmx$ and
$n_{\varepsilon,\lambda}:\rre\times\rre^{d}\to\rre$ are unknown
functions, and $E_{0}:\rre^{d}\to\cmx$, $n_{0}:\rre^{d}\to\rre$ and
$n_{1}:\rre^{d}\to\rre$ are given functions.

This model was introduced by Garcia-Haas-Oliveira-Goedert \cite{GHPG}
and Haas-Shukla \cite{HS} to describe the nonlinear interaction between
high-frequency quantum Langmuir waves and low-frequency quantum
ion-acoustic waves. The physical background of the system (\ref{ZK}) can
be found in \cite{Haas}.

The classical Zakharov system
\begin{equation}
%\begin{eqnarray}
\left\{
\begin{array}{l}
\medskip
i\pt_{t}E_{0,\lambda}
+\Delta E_{0,\lambda}
=n_{0,\lambda}E_{0,\lambda}
\qquad t\in\rre,\ x\in\rre^{d},\\
\medskip
\displaystyle{
\lambda^{-2}\pt_{t}^{2}n_{0,\lambda}
-\Delta n_{0,\lambda}
=\Delta|E_{0,\lambda}|^{2}
\qquad t\in\rre,\ x\in\rre^{d}}
\end{array}
\right. \label{cZK}
%\end{eqnarray}
\end{equation}
was proposed by Zakharov \cite{Z} as a model for describing the
interaction between the Langmuir waves and ion-acoustic waves in a
plasma. In the system (\ref{cZK}), $E_{0,\lambda}$ denotes the slowly
varying envelope of the highly oscillatory electric field and
$n_{0,\lambda}$ denotes the deviation of the ion density from the
equilibrium, and $\lambda$ is the ionic speed of sound.

The classical Zakharov system (\ref{cZK}) has been extensively studied
from the point of view of local and global well-posedness
\cite{BC,BH,BHHT,GTV,KPV,K2,OT1}, blow-up of solutions
\cite{GM1,GM2,KM}, scattering \cite{GN,GNW,HPS,OT3,S}, and subsonic
limit of (\ref{cZK}) \cite{AA,MN1,MN2,OT2,SW}.

Let us review the local and global well-posedness %, namely, existence
% and uniqueness of solution and continuous dependence on initial data
% of solution
for the classical Zakharov system (\ref{cZK}). We will mention the Cauchy
problem of (\ref{cZK}) in $\rre^{d}$, $d=1,2,3$ only. For the high
dimensional case and the multi-dimensional torus case, see \cite{GTV}
and \cite{K2}, respectively.

Bourgain-Colliander \cite{BC} proved the local well-posedness of
(\ref{cZK}) in $H^{1}(\rre^{d})\times L^{2}(\rre^{d})
\times\dot{H}^{-1}(\rre^{d})$ for $d=2,3$ by using the Fourier restriction
norm associated to the free Schr\"{o}dinger and wave evolution groups.
By refining the Fourier restriction norm used in \cite{BC},
Ginibre-Tsutsumi-Velo \cite{GTV} have shown the local well-posedness
of (\ref{cZK}) in $H^{k}(\rre^{d})\times H^{\ell}(\rre^{d}) \times
H^{\ell-1}(\rre^{d})$ provided that $-1/2<k-\ell\le1$ and $0\le\ell+1/2\le
2k$ for $d=1$, $0\le k-\ell\le1$ and  $1\leq\ell+1\le 2k$ for $d=2,3$.

For one dimensional case, Colliander-Holmer-Tzirakis \cite{CHT} proved
the global well-posedness of (\ref{cZK}) in $L^{2}(\rre)\times
H^{-1/2}(\rre)\times H^{-3/2}(\rre)$. For the two dimensional case,
Bejenaru-Herr-Holmer-Tataru \cite{BHHT} showed the local
well-posedness of (\ref{cZK}) in $L^{2}(\rre^{2})\times H^{-1/2}(\rre^{2})
\times H^{-3/2}(\rre^{2})$ which is optimal from the point of view of the
scaling of Sobolev spaces. Furthermore, for the three dimensional case,
Bejenaru-Herr \cite{BH} have recently shown the local well-posedness of
(\ref{cZK}) in $H^{k}(\rre^{3})\times H^{\ell}(\rre^{3}) \times
H^{\ell-1}(\rre^{3})$ with $k>0$ and $\ell>-1/2$.

Like the classical model (\ref{cZK}), the quantum Zakharov system (\ref{ZK})
possesses the conservation of mass
\begin{eqnarray*}%%%%%%%%%%%%%%%%%%%%%%%%%%%%%%%%%
\|E_{\varepsilon,\lambda}(t)\|_{L_{x}^{2}}^{2}=\text{constant}
\end{eqnarray*}%%%%%%%%%%%%%%%%%%%%%%%%%%%%%%%%%
and the conservation of the Hamiltonian
\begin{equation}\begin{array}{rll}%%%%%%%%%%%%%%%%%%%%%%%%%%%%%%%%%
\lefteqn{\|\nabla E_{\varepsilon,\lambda}(t)\|_{L_{x}^{2}}^{2}
+\varepsilon^{2}
\|\Delta E_{\varepsilon,\lambda}(t)\|_{L_{x}^{2}}^{2}
+\frac12\lambda^{-2}\|\pt_{t}\nabla^{-1}
n_{\varepsilon,\lambda}(t)\|_{L_{x}^{2}}^{2}}\\
& &
+\frac12\|n_{\varepsilon,\lambda}(t)\|_{L_{x}^{2}}^{2}
+\frac{\varepsilon^{2}}{2}\|
\nabla n_{\varepsilon,\lambda}(t)\|_{L_{x}^{2}}^{2}
+\int_{\rre^{d}}n_{\varepsilon,\lambda}
|E_{\varepsilon,\lambda}|^{2}dx
=\text{constant}.
\end{array}\label{Hami} \end{equation} %%%%%%%%%%%%%%%%%%%%%%%%%%%%%%%%%
Compared to the classical Zakharov system (\ref{cZK}),
there are few results for the quantum Zakharov system (\ref{ZK}). We
summarize the known results for (\ref{ZK}).

In one space dimension, the system (\ref{ZK}) is studied from the point
of view of the existence of exact solution and the  local and global
well-posedness. El-Wakil and Abdou \cite{EA} constructed the exact
traveling solutions of (\ref{ZK}) by using the improved $tanh$ function
method. Jiang-Lin-Shao \cite{JLP} proved the local well-posedness of
(\ref{ZK}) in $H^{k}(\rre)\times H^{\ell}(\rre) \times H^{\ell-2}(\rre)$
provided that $-3/2<k-\ell<3/2$, $-3/2<2k-\ell$, $-3/2<k+\ell$ and
$k>-3/4$. Chen-Fang-Wang \cite{CFW} have recently shown the global
well-posedness of (\ref{cZK}) in $L^{2}(\rre)\times H^{\ell}(\rre)\times
H^{\ell-2}(\rre)$ with $-3/2\le\ell\le3/2$. The proofs given in \cite{JLP} and
\cite{CFW} are based on the Fourier restriction method associated to the
fourth order Schr\"{o}dinger and the fourth order wave evolution
groups.

For space dimensions $d=1,2,3$, Guo-Zhang-Guo \cite{GZG} have proved
the global well-posedness of (\ref{ZK}) with initial data in
$H^{k}(\rre^{d})\times H^{k-1}(\rre^{d}) \times
(H^{k-3}\cap\dot{H}^{-1})(\rre^{d})$ and in $H^{k}(\rre^{d})\times
H^{k-1}(\rre^{d}) \times H^{k-3}(\rre^{d})$ for $d=1,2,3$ and $k\ge2$,
respectively. Especially, they proved the global well-posedness of the
quantum Zakharov system (\ref{ZK}) in the energy space 
$H^{2}(\rre^{d})\times H^{1}(\rre^{d}) \times
(H^{-1}\cap\dot{H}^{-1})(\rre^{d})$. However it is interesting that the
classical Zakharov system (\ref{cZK}) has the blow up solution in the
energy space which is $H^{1}(\rre^{d})\times L^{2}(\rre^{d}) \times
(H^{-1}\cap\dot{H}^{-1})(\rre^{d})$, see Glangetas-Merle
\cite{GM1,GM2}. As pointed out in \cite{GZG}, this difference is caused
by the strong dispersion stems from the quantum effect, namely the
quantum effect stabilizes the solution.
%Guo-Zhang-Guo \cite{GZG} also studied the classical limit
%behavior of the solution to the quantum
%Zakharov system (\ref{ZK}).

In this paper, we consider the convergence of the solution to the
quantum Zakharov system (\ref{ZK}) as $\lambda\to\infty$. Before
stating the main result for (\ref{ZK}), we review the convergence of the
solution to the classical Zakharov system (\ref{cZK}) as
$\lambda\to\infty$.

Let us formally take $\lambda\to\infty$ for the second equation in
(\ref{cZK}). If $n_{0,\infty} +|E_{0,\infty}|^{2}$ vanishes at space infinity,
we obtain the relation $n_{0,\infty}=-|E_{0,\infty}|^{2}$. Substituting
this relation into the first equation in (\ref{cZK}), we see that
$E_{0,\infty}$ satisfies the focusing cubic nonlinear Schr\"{o}dinger 
equation
\begin{equation} %%%%%%%%%%%%%%%%%%%%%%%%%%%%%%%%%
%\begin{eqnarray}
i\pt_{t}E_{0,\infty}+\Delta
E_{0,\infty}=-|E_{0,\infty}|^{2}E_{0,\infty}.
\label{NLS}
%\end{eqnarray}
\end{equation} %%%%%%%%%%%%%%%%%%%%%%%%%%%%%%%%%%
The local existence and uniqueness of (\ref{cZK}) is shown by
Schochet-Weinstein \cite{SW} when the initial datum lies in
$H^{k}(\rre^{d})\times H^{k-1}(\rre^{d}) \times
(H^{k-2}\cap\dot{H}^{-1})(\rre^{d})$ with $d=1,2,3$ and $k\ge2$.
%Schochet-Weinstein \cite{SW} showed a local existence and uniqueness
%issue in $H^{k}(\rre^{d})\times H^{k-1}(\rre^{d}) \times
%(H^{k-2}\cap\dot{H}^{-1})(\rre^{d})$ with $d=1,2,3$ and $k\ge2$,
Furthermore, they proved that the existence time interval is
independent of $\lambda\in[1,\infty)$ and the solution
$(E_{0,\lambda},n_{0,\lambda})$ of (\ref{cZK}) converges to
$(E_{0,\infty},-|E_{0,\infty}|^{2})$ as $\lambda\to\infty$, where
$E_{0,\infty}$ is a solution to the focusing cubic nonlinear Schr\"{o}dinger
equation (\ref{NLS}). Added-Added \cite{AA} studied the rate of
convergence of solution. The optimal rates of convergence for
$E_{0,\lambda}$ is given by Ozawa-Tsutsumi \cite{OT2}, including some
discussion of initial layer phenomenon. Notice that in \cite{SW,AA,OT2},
they imposed the assumption $n_{1}\in \dot{H}^{-1}(\rre^{d})$.
Kenig-Ponce-Vega \cite{KPV} obtained the optimal rates of convergence
for $E_{0,\lambda}$ in the non-compatible case $n_{0}+|E_{0}|^{2}
\neq0$ without the assumption $n_{1}\in \dot{H}^{-1}(\rre^{d})$.
%\textrr{What about the compatible case?}

Let us turn to the quantum system (\ref{ZK}). In analogy with the
classical system (\ref{cZK}), taking $\lambda\to\infty$ for the second
equation in (\ref{ZK}) and together with the assumption
$n_{\varepsilon,\infty} +(1-\varepsilon^{2}
\Delta)^{-1}|E_{\varepsilon,\infty}|^{2}$ vanishes at space infinity, we
obtain the relation
\begin{eqnarray*}%%%%%%%%%%%%%%%%%%%%%%%%%%%%%%%%%
n_{\varepsilon,\infty}=-(1-\varepsilon^{2}
\Delta)^{-1}|E_{\varepsilon,\infty}|^{2}.
\end{eqnarray*}%%%%%%%%%%%%%%%%%%%%%%%%%%%%%%%%%
Substituting this relation into the first equation in (\ref{ZK}), we see that
$E_{\varepsilon,\infty}$ satisfies the quantum modified nonlinear
Schr\"{o}dinger type equation
\begin{equation} %%%%%%%%%%%%%%%%%%%%%%%%%%%%%%%%%
%\begin{eqnarray}
i\pt_{t}E_{\varepsilon,\infty}-(-\Delta
+\varepsilon^{2}\Delta^{2})
E_{\varepsilon,\infty}=-\{(1-\varepsilon^{2}
\Delta)^{-1}|E_{\varepsilon,\infty}|^{2}\}
E_{\varepsilon,\infty}.
\label{4NLS}
%\end{eqnarray}
\end{equation} %%%%%%%%%%%%%%%%%%%%%%%%%%%%%%%%%
The main purpose of this paper is to prove that the solution
$E_{\varepsilon,\lambda}$ to the quantum Zakharov system converges
to the solution $E_{\varepsilon,\infty}$ to the fourth order nonlinear
Schr\"{o}dinger type equation (\ref{4NLS}).

We introduce several function spaces and notations. For non-negative
integers $m,n$, the Sobolev space $H^{m}$ and the weighted Sobolev
space $H^{m,n}$ are defined by
\begin{eqnarray*}%%%%%%%%%%%%%%%%%%%%%%%%%%%%%%%%%
H^{m}(\rre^{d})&=&
\Big\{u\in{{\mathcal S}}'(\rre^{d});
\|u\|_{H^{m}}=
\sum_{k=0}^{m}\|\nabla^{k}u\|_{L^{2}}<\infty\Big\},\\
H^{m,n}(\rre^{d})&=&
\Big\{u\in{{\mathcal S}}'(\rre^{d});
\|u\|_{H^{m,n}}=
\sum_{k=0}^{m}\sum_{\ell=0}^{n}
\||x|^{\ell}\nabla^{k}u\|_{L^{2}}
<\infty\Big\},
\end{eqnarray*}%%%%%%%%%%%%%%%%%%%%%%%%%%%%%%%%%
where the k-$th$ derivative $\nabla^{k}$ is defined by
\begin{eqnarray*}%%%%%%%%%%%%%%%%%%%%%%%%%%%%%%%%%
\nabla^{k}=\left\{
\begin{array}{l}
\Delta^{(k-1)/2}\nabla\qquad\ \ \ \text{if}\ k\ \text{is\ odd},\\
\Delta^{k/2}\qquad \qquad \quad \ \text{if}\ k\ \text{is\ even.}
\end{array}
\right.
\end{eqnarray*}%%%%%%%%%%%%%%%%%%%%%%%%%%%%%%%%%
We also define the following notations
$$%%%%%%%%%%%%%%%%%%%%%%%%%%%%%%%%%
\Delta_\varepsilon = \Delta-\varepsilon^2\Delta^2,\quad
|\nabla| = \sqrt{-\Delta}, \quad\text{and}\quad I_\varepsilon = (1-\varepsilon^2\Delta)^{-1}.
$$ %%%%%%%%%%%%%%%%%%%%%%%%%%%%%%%%%
Following Kishimoto-Maeda \cite{KM}, we define the operator
$\nabla^{-1}$ via the Helmholtz decomposition.
% $L^{2}(\rre^{d})^{d}=G\oplus L_{\sigma}^{2}$, where $G=\{\nabla
% g:\,g\in\dot{H}^{1}(\rre^{d})\}$ and $L_{\sigma}^{2}=\{h\in
% L^{2}(\rre^{d})^{d}:\,\text{div}h=0\}$. Let ${{\mathbb
% P}}_{G}:L^{2}(\rre^{d})^{d}\to G$ be the orthogonal projection onto $G$.
% For any function $f$ given by $f=-\nabla\cdot q$, $q\in
% L^{2}(\rre^{d})^{d}$, we define $\nabla^{-1}f$ by
% $\nabla^{-1}f={{\mathbb P}}_{G} q$. The orthogonal projection $P_G$ is
% also called the Leray projector. It plays the most important role in
% mathematical study of fluid dynamics \cite{Lions}. Furthermore, we
% define the function space $\dot{H}^{-1}(\rre^{d})$ by
% \begin{eqnarray*}%%%%%%%%%%%%%%%%%%%%%%%%%%%%%%%%%
% \dot{H}^{-1}(\rre^{d})
% =\{-\nabla\cdot{{\mathbb P}}_{G}q:\,q
% \in L^{2}(\rre^{d})^{d}\}.
% \end{eqnarray*}%%%%%%%%%%%%%%%%%%%%%%%%%%%%%%%%%
The homogeneous Sobolev space $\dot{H}^{-\sigma}$ is defined by
 \begin{eqnarray*}%%%%%%%%%%%%%%%%%%%%%%%%%%%%%%%%%
\dot{H}^{-\sigma}(\rre^{d})=\{f\in{{\mathcal S}}'(\rre^{d});
\|f\|_{\dot{H}^{-\sigma}}=\||\xi|^{-\sigma}\hat{f}(\xi)\|_{L_{\xi}^{2}}<\infty\}.
\end{eqnarray*}%%%%%%%%%%%%%%%%%%%%%%%%%%%%%%%%%
%The space $\dot{H}^{-1}$ is equivalent to $\dot{H}^{-\sigma}$ for $\sigma=1$.
Let $m$ and $M$ be integers throughout the paper. For the sake of
convenience, we define
\begin{eqnarray*}%%%%%%%%%%%%%%%%%%%%%%%%%%%%%%%%
C([0,\infty);\,X\times Y\times Z)=C([0,\infty);\,X)\times C([0,\infty);\,Y)
\times C([0,\infty);\,Z),
\end{eqnarray*} %%%%%%%%%%%%%%%%%%%%%%%%%%%%%%%%%
where $X$, $Y$, and $Z$ are Sobolev spaces. From now on, we drop the
parameter $\varepsilon$ of the solutions $E_{\varepsilon,\lambda}$ and
$n_{\varepsilon,\lambda}$. We also denote the space, for $d=1, 2, 3$,
 \begin{equation} %%%%%%%%%%%%%%%%%%%%%%%%%%%%%%
X_{M, d} := \,H^{M}(\rre^d)\times H^{M-1}(\rre^d) \times
(H^{M-3}(\rre^d)\cap \dot{H}^{-1}(\rre^d)). \label{XMd}
\end{equation}  %%%%%%%%%%%%%%%%%%%%%%%%%%%%%%%%%

%%%%%%%%%%%%%% (-\Delta+\varepsilon^{2}\Delta^{2}) --> (-\Delta_\varepsilon)
%%%%%%%%%%%%%% -\Delta+\varepsilon^{2}\Delta^{2} --> -\Delta_\varepsilon
%%%%%%%%%%%%%% {-\Delta+\varepsilon^{2}\Delta^{2}} --> {-\Delta_\varepsilon}
%%%%%%%%%%%%%% \Delta-\varepsilon^{2}\Delta^{2} --> \Delta_\varepsilon
%%%%%%%%%%%%%% (1-\varepsilon^{2}\Delta) --> I_\varepsilon

Our main results are as follows.
\begin{thm}[Case $d=1$]\label{main1}%%%%%%%%%%%%%%%%%%%%%%%% T1
Let $d=1$ and $M\ge2$. Then for any $(E_{0},n_{0},n_{1})\in X_{3M, 1}$
and $\lambda\in[1,\infty)$, there exist a unique solution to (\ref{ZK})
satisfying
%\begin{eqnarray*}%%%%%%%%%%%%%%%%%%%%%%%%%%%%%%
%(E_{\varepsilon,\lambda},n_{\varepsilon,\lambda},
%\pt_{t}n_{\varepsilon,\lambda})
%&\in&C([0,\infty);H^{3M}(\rre))
%\times
%C([0,\infty);H^{3M-1}(\rre))\\
%& &\qquad\qquad\times C
%([0,\infty);H^{3M-3}(\rre)\cap
%\dot{H}^{-1}(\rre))
%\end{eqnarray*}%%%%%%%%%%%%%%%%%%%%%%%%%%%%%%%%%
 \begin{equation*} %%%%%%%%%%%%%%%%%%%%%%%%%%%%%%
(E_{\lambda},n_{\lambda},
\pt_{t}n_{\lambda}) \in C([0,\infty);\, X_{3M, 1})
\end{equation*} %%%%%%%%%%%%%%%%%%%%%%%%%%%%%%%%%
and a unique solution to (\ref{4NLS}) satisfying
\begin{equation*}%%%%%%%%%%%%%%%%%%%%%%%%%%%%%%%%%
E_{\infty}\in C([0,\infty);H^{3M}(\rre)).
\end{equation*}%%%%%%%%%%%%%%%%%%%%%%%%%%%%%%%%%
\vskip1mm \noindent (i) Assume
$n_{0}+I_{\varepsilon}|E_{0}|^{2}\neq0$. Let $m\ge3$ and $M\ge3$
satisfy $3M\ge m+6$. If we further assume that $E_{0}\in
H^{0,M}(\rre)$, $x^{j}(n_{0}+I_{\varepsilon}|E_{0}|^{2})\in H^{m+j-3}$
with $j=0,1,2$,
%and $n_{0}+{I_{\varepsilon}}|E_{0}|^{2}\in\dot{H}^{-\sigma} (\rre)$ with $1/2<\sigma<3/2$,
then for any $T\in(0,\infty)$,
% there exists a positive constant $C$ such that
$(E_{\lambda}, n_{\lambda})$ satisfies
\begin{equation}
%\begin{eqnarray}%%%%%%%%%%%%%%%%%%%%%%%%%%%%%%%%%
\sup_{0\le t\le T}
\|E_{\lambda}(t)-
E_{\infty}(t)\|_{H^{m}}
\le C\lambda^{-1}
\label{T1-1}
%\end{eqnarray}%%%%%%%%%%%%%%%%%%%%%%%%%%%%%%%%%
\end{equation}
and
\begin{equation}
%\begin{eqnarray}%%%%%%%%%%%%%%%%%%%%%%%%%%%%%%%%%
\sup_{0\le t\le T}
\|n_{\lambda}(t)
+{I_{\varepsilon}}|E_{\lambda}|^{2}(t)
-Q_{\lambda}^{(0)}(t)\|_{H^{m}}
\le C\lambda^{-1}
\label{T1-2}
%\end{eqnarray}%%%%%%%%%%%%%%%%%%%%%%%%%%%%%%%%%
\end{equation}
for any $\lambda\in[1,\infty)$, where $Q_{\lambda}^{(0)} =\cos(\lambda
t\sqrt{-\Delta_\varepsilon}) \{n_{0}+{I_{\varepsilon}}|E_{0}|^{2}\}$ and
the constant $C$ depends on $\varepsilon$ and $T$, but independent of
$\lambda\in[1,\infty)$.

\vskip1mm

\noindent (ii) Assume $n_{0}+{I_{\varepsilon}}|E_{0}|^{2}\equiv0$. Let
$m\ge6$ and $M\ge4$ satisfy $3M\ge m+5$. If we further assume that
$E_{0}\in H^{0,M}(\rre)$ and $x^{j}(\nabla^{-1}n_{1}+2Im\{ E_{0}\nabla
{I_{\varepsilon}^{-1}}\overline{E}_{0} +\varepsilon^{2} \nabla
E_{0}\Delta_{\varepsilon}\overline{E_{0}}\}) \in H^{m+j-2}$ with
$j=0,1,2$,
% and $\nabla^{-1}n_{1} +2Im\{E_{0}\nabla {I_{\varepsilon}^{-1}}\overline{E}_{0}
% +\varepsilon^{2} \nabla
% E_{0}\Delta_{\varepsilon}\overline{E_{0}}\}\in\dot{H}^{-\sigma}$ with
% $1/2<\sigma<3/2$,
then for any $T\in(0,\infty)$, % there exists a positive constant $C$ such that
$(E_{\lambda}, n_{\lambda})$ satisfies
\begin{equation}
%\begin{eqnarray}%%%%%%%%%%%%%%%%%%%%%%%%%%%%%%%%%
\sup_{0\le t\le T}
\|E_{\lambda}(t)-
E_{\infty}(t)\|_{H^{m}}
\le C\lambda^{-2}
\label{T1-3}
%\end{eqnarray}%%%%%%%%%%%%%%%%%%%%%%%%%%%%%%%%%
\end{equation}
for any $\lambda\in[1,\infty)$.
\end{thm}%%%%%%%%%%%%%%%%%%%%%%%%%%%%%%%%%%%%%% T1

\begin{thm}[Case $d=2$]\label{main2}%%%%%%%%%%%%%%%%%%%%% T2
Let $d=2$ and $M\ge2$. Then for any $(E_{0},n_{0},n_{1})\in{X_{3M, 2}}$
and $\lambda\in[1,\infty)$, there exists a unique solution to (\ref{ZK})
satisfying
\begin{equation*} %%%%%%%%%%%%%%%%%%%%%%%%%%%%%%
(E_{\lambda},n_{\lambda},
\pt_{t}n_{\lambda}) \in C([0,\infty);
\,{X_{3M, 2}})
\end{equation*} %%%%%%%%%%%%%%%%%%%%%%%%%%%%%%%%%
and a unique solution to (\ref{4NLS}) satisfying
\begin{eqnarray*}
E_{\infty}\in
C([0,\infty);H^{3M}(\rre^{2})).
\end{eqnarray*}
\vskip1mm \noindent (i) Assume
$n_{0}+{I_{\varepsilon}}|E_{0}|^{2}\neq0$. Let $m\ge3$ and $M\ge3$
satisfy $3M\ge m+6$. If we further assume that $E_{0}\in H^{0,M}$,
$|x|^{j} (n_{0}+{I_{\varepsilon}}|E_{0}|^{2})\in H^{m+j-2}$ with
$j=0,1,2$, and $n_{0}+{I_{\varepsilon}}|E_{0}|^{2}\in \dot{H}^{-\sigma}$,
with $0<\sigma<1$, then for any $T\in(0,\infty)$,
% there exists a positive constant $C$ such that
$(E_{\lambda}, n_{\lambda})$ satisfies (\ref{T1-1}) and (\ref{T1-2}) for
any $\lambda\in[1,\infty)$, where the constant $C$ depends on
$\varepsilon$ and $T$. \vskip1mm \noindent (ii) Assume
$n_{0}+{I_{\varepsilon}}|E_{0}|^{2}\equiv0$. Let $m\ge5$ and $M\ge4$
satisfy $3M\ge m+6$. If we further assume that $E_{0}\in H^{0,M}$,
$|x|^{j}(\nabla^{-1}n_{1}+2Im\{E_{0}\nabla
{I_{\varepsilon}^{-1}}\overline{E}_{0} +\varepsilon^{2} \nabla
E_{0}\Delta_{\varepsilon}\overline{E_{0}}\}) \in H^{m+j-2}$ with
$j=0,1,2$ and $\nabla^{-1}n_{1}+2Im\{E_{0}\nabla
{I_{\varepsilon}^{-1}}\overline{E}_{0} +\varepsilon^{2} \nabla
E_{0}\Delta_{\varepsilon}\overline{E_{0}}\}\in\dot{H}^{-\sigma}$ with
$0<\sigma<1$, then for any $T\in(0,\infty)$,
% there exists a positive constant $C$ such that
$(E_{\lambda}, n_{\lambda})$ satisfies
\begin{equation}
%\begin{eqnarray}%%%%%%%%%%%%%%%%%%%%%%%%%%%%%%%%%
\sup_{0\le t\le T}
\|E_{\lambda}(t)-
E_{\infty}(t)\|_{H^{m}}
\le C\lambda^{-2}\log\lambda
\label{T2-3}
%\end{eqnarray}%%%%%%%%%%%%%%%%%%%%%%%%%%%%%%%%%
\end{equation}
for any $\lambda\in[1,\infty)$, where the constant $C$ depends on
$\varepsilon$ and $T$.
\end{thm}%%%%%%%%%%%%%%%%%%%%%%%%%%%%%%%%%%%%%%%%% T2

\begin{thm}[Case $d=3$]\label{main3}%%%%%%%%%%%%%%%%%%%%% T3
Let $d=3$ and $M\ge2$. Then for any $(E_{0},n_{0},n_{1})\in {X_{3M,
3}}$ and $\lambda\in[1,\infty)$, there exists a unique solution to
(\ref{ZK}) satisfying
\begin{equation*} %%%%%%%%%%%%%%%%%%%%%%%%%%%%%%
(E_{\lambda},n_{\lambda},
\pt_{t}n_{\lambda}) \in C([0,\infty);
\,{X_{3M, 3}}))
\end{equation*} %%%%%%%%%%%%%%%%%%%%%%%%%%%%%%%%%
and a unique solution to (\ref{4NLS}) satisfying
\begin{eqnarray*}
E_{\infty}\in
C([0,\infty);H^{3M}(\rre^{3})).
\end{eqnarray*}
\vskip1mm \noindent (i) Assume
$n_{0}+{I_{\varepsilon}}|E_{0}|^{2}\neq0$. Let $m\ge3$ and $M\ge3$
satisfy $3M\ge m+6$. If we further assume that $E_{0}\in H^{0,M}$,
$|x|^{j} (n_{0}+{I_{\varepsilon}}|E_{0}|^{2})\in H^{m+j-2}$ with
$j=0,1,2$, then %there exists a $T^{\ast}>0$ depending only on
%$\varepsilon$ such that
for any $T\in(0,\infty)$,
% there exists a positive constant $C$ such that
$(E_{\lambda}, n_{\lambda})$ satisfies (\ref{T1-1}) and (\ref{T1-2}) for
any $\lambda\in[1,\infty)$, where the constant $C$ depends on
$\varepsilon$ and $T$.

\vskip1mm

\noindent (ii) Assume $n_{0}+{I_{\varepsilon}}|E_{0}|^{2}\equiv0$. Let
$m\ge5$ and $M\ge4$ satisfy $3M\ge m+6$. If we further assume that
$E_{0}\in H^{0,M}$, $|x|^{j}(\nabla^{-1}n_{1}+2Im\{E_{0}\nabla
{I_{\varepsilon}^{-1}}\overline{E}_{0} +\varepsilon^{2} \nabla
E_{0}\Delta_{\varepsilon}\overline{E_{0}}\}) \in H^{m+j-2}$ with
$j=0,1,2$, then
%there exists a $T^{\ast}>0$ depending only on
%$\varepsilon$ such that
for any $T\in(0,\infty)$,
% there exists a positive constant $C$ such that
$(E_{\lambda}, n_{\lambda})$ satisfies
\begin{equation}
%\begin{eqnarray}%%%%%%%%%%%%%%%%%%%%%%%%%%%%%%%%%
\sup_{0\le t\le T}
\|E_{\lambda}(t)-
E_{\infty}(t)\|_{H^{m}}
\le C\lambda^{-2}
\label{T3-3}
%\end{eqnarray}%%%%%%%%%%%%%%%%%%%%%%%%%%%%%%%%%
\end{equation}
for any $\lambda\in[1,\infty)$, where the constant $C$ depends on
$\varepsilon$ and $T$.
\end{thm}%%%%%%%%%%%%%%%%%%%%%%%%%%%%%%%%%%%%%%%%% T3

\vskip2mm

\begin{rem}%%%%%%%%%%%%%%%%%%%%%%%%%%%%%%%%%%%%% Rmk1.1
The function $Q_{\varepsilon, \lambda}^{(0)}$ is solution to the fourth
order wave equation
\begin{eqnarray*}
\left\{
\begin{array}{l}
\medskip
\displaystyle{
\lambda^{-2}\pt_{t}^{2}n
-\Delta_\varepsilon n
=0
\qquad t\in\rre,\ x\in\rre^{d},}\\
\displaystyle{
n(0,x)=n_{0}+{I_{\varepsilon}}|E_{0}|^{2},
\ \pt_{t}n(0,x)
=0\qquad x\in\rre^{d}.}
\end{array}
\right.
\end{eqnarray*}
Theorems \ref{main1}, \ref{main2} and \ref{main3} tell us that the term
$Q_{\varepsilon, \lambda}^{(0)}$ represents the initial layer for
(\ref{ZK}). Note that $Q_{0,\lambda}^{(0)}$ coincides with the initial
layer for the Zakharov system (\ref{cZK}), see \cite{AA,SW,OT2}.
\end{rem}%%%%%%%%%%%%%%%%%%%%%%%%%%%%%%%%%%%%%%%%

%\begin{rem} %%%%%%%%%%%%%%%%%%%%%%%%%%%%%%%%%%%% Rmk 1.2
%For $d= 1, 2$, we can take $T$
% in Theorems \ref{main1} and \ref{main2}
%arbitrary large, while in dimensions three, the range of $T$ is
%$(0,T^{\ast})$, where $T^{\ast}$ is bounded. This difference is due to the
%Sobolev embedding, namely the embedding $H^{1}\hookrightarrow
%L^{\infty}$ holds for $d=1$, but not for $d=2, 3$. For $d=2$, however,
%employing the Br\'{e}zis-Gallou\"{e}t inequality (Lemma \ref{BG}), we
%can overcome this difficulty. See the proof of Proposition \ref{well} in
%Section 2.
%
% In fact, we can remove the boundedness of $T^{\ast}$ in Theorem \ref{main3} if we use
%the Strichartz estimates for \eqref{ZK}, instead of using the Gronwall's
%inequality. However we will not present it here.
%\end{rem}%%%%%%%%%%%%%%%%%%%%%%%%%%%%%%%%%%%%

%\begin{rem} %%%%%%%%%%%%%%%%%%%%%%%%%%%%%%%%%%%% Rmk1.3
% In Theorems \ref{main2},
%we assume that $n_{0}+{I_{\varepsilon}}|E_{0}|^{2}\in
%\dot{H}^{-\sigma}$ with  $\sigma>0$ for $d=2$. Note that this
%assumption is not needed for dimension one and three, see the proof of
%\textrr{ Lemmas \ref{1f} and \ref{3f} } in Section 4.
%\end{rem} %%%%%%%%%%%%%%%%%%%%%%%%%%%%%%%%%%%%

Let us give an outline of the proofs for Theorems \ref{main1},
\ref{main2} and \ref{main3}. The proofs follow from the arguments due
to Ozawa-Tsutsumi \cite{OT2} and Ukai \cite{U} with several
modifications.

From (\ref{ZK}) and (\ref{4NLS}), the difference $E_{\lambda}- E_{\infty}$
satisfies
\begin{equation}
%\begin{eqnarray} %%%%%%%%%%%%%%%%%%%%%%%%%%%%%%%%%%%%
\left\{
\begin{array}{l}%%%%%%%%%%%%%%%%%%%%%%%%%%%%%%%%%%%%
\medskip
i\pt_{t}(E_{\lambda}-
E_{\infty})
+\Delta_\varepsilon
(E_{\lambda}-
E_{\infty})\\
\medskip
\qquad\qquad=
-\{{I_{\varepsilon}}|E_{\lambda}|^{2}\}E_{\lambda}
+\{{I_{\varepsilon}}|E_{\infty}|^{2}\}E_{\infty}+ Q_{\lambda}E_{\lambda},\\
\medskip
(E_{\lambda}-
E_{\infty})(0,x)=0,
\end{array}\right.\label{I}  %%%%%%%%%%%%%%%%%%%%%%%%%%%%
%\end{eqnarray}%%%%%%%%%%%%%%%%%%%%%%%%%%%%%%%%%%%%
\end{equation}
where $Q_{\lambda} =n_{\lambda}+{I_{\varepsilon}}|E_{\lambda}|^{2}$.
To evaluate $E_{\lambda} -E_{\infty}$, we rewrite  (\ref{I}) into the
integral equation
%\begin{eqnarray}%%%%%%%%%%%%%%%%%%%%%%%%%%%%%%%%%%%%
%& &E_{\lambda}(t)-E_{\infty}(t)=\nonumber\\
%& & \int_{0}^{t}
%iU_{\varepsilon}(t-s)\Big(\big[
%\{{I_{\varepsilon}}|E_{\lambda}|^{2}\}E_{\lambda}-
%\{{I_{\varepsilon}}|E_{\infty}|^{2}\}E_{\varepsilon,
%\infty}\big]
%- (Q_{\lambda}
%E_{\lambda})\Big)(s)ds,
%\label{EE}
%\end{eqnarray}%%%%%%%%%%%%%%%%%%%%%%%%%%%%%%%%%%%%
\begin{equation}%%%%%%%%%%%%%%%%%%%%%%%%%%%%%%%%%%%%
E_{\lambda}(t)-E_{\infty}(t)=
\int_{0}^{t} iU_{\varepsilon}(t-s)\Big(\big[
\{{I_{\varepsilon}}|E_{\lambda}|^{2}\}E_{\lambda}-
\{{I_{\varepsilon}}|E_{\infty}|^{2}\}E_{\infty}\big] -
(Q_{\lambda} E_{\lambda})\Big)ds, \label{E-E}
\end{equation}%%%%%%%%%%%%%%%%%%%%%%%%%%%%%%%%%%%%
where $U_{\varepsilon}(t)= \exp(it\Delta_\varepsilon)$ is an
$L^{2}$-unitary group generated by the differential operator
$i\Delta_\varepsilon$. It is easy to evaluate the first term in the integral
of (\ref{E-E}) in $H_{x}^{m}$. To estimate the second term in the integral
of (\ref{E-E}) in $H_{x}^{m}$, we derive the relevant equation for
$Q_{\lambda}$:
\begin{equation}\left\{\begin{array}{l}%%%%%%%%%%%%%%%%%%%%%%%
\medskip
\displaystyle{
\lambda^{-2}\pt_{t}^{2}Q_{\lambda}
-\Delta_\varepsilon Q_{\lambda}=
\lambda^{-2}\pt_{t}^{2}{I_{\varepsilon}}|E_{\lambda}|^{2}
\qquad t\in\rre,\ x\in\rre^{d},}\\
\medskip\displaystyle{
Q_{\lambda}(0,x)=n_{0}(x)+{I_{\varepsilon}}|E_{0}|^{2}(x)
\qquad x\in\rre^{d},}\\
\displaystyle{\pt_{t}Q_{\lambda}(0,x)
=n_{1}(x)+2Im
\{E_{0}\Delta_\varepsilon\overline{E}_{0}\}
\qquad x\in\rre^{d},}
\end{array}\right. \label{Q eq}
\end{equation}%%%%%%%%%%%%%%%%%%%%%%%%%%%
which implies that $Q_{\lambda}(t)$ can be written as
%\begin{eqnarray} %%%%%%%%%%%%%%%%%%%%%%%%%%%%%%%%%%%%
\begin{equation} %%%%%%%%%%%%%%%%%%%%%%%%%%%%%%%%%%%%
\begin{array}{ll}
& \displaystyle
\cos(\lambda t\sqrt{-\Delta_\varepsilon})\{n_{0}+{I_{\varepsilon}}|E_{0}|^{2}\}
+\frac{\sin(\lambda t\sqrt{-\Delta_\varepsilon})}{\lambda
\sqrt{-\Delta_\varepsilon}}\{n_{1}+2Im
(E_{0}\Delta_\varepsilon\overline{E}_{0})\}\\
&\displaystyle +
\int_{0}^{t}
\frac{\sin(\lambda (t-s)\sqrt{-\Delta_\varepsilon})}{\lambda
\sqrt{-\Delta_\varepsilon}}
\pt_{t}^{2}{I_{\varepsilon}}|E_{\lambda}|^{2}
(s)ds.
\end{array}\label{INT}
%\end{eqnarray} %%%%%%%%%%%%%%%%%%%%%%%%%%%%%%%%%%%%
\end{equation} %%%%%%%%%%%%%%%%%%%%%%%%%%%%%%%%%%%%
% It is not difficult to estimate the second and third terms on the right % hand side of (\ref{INT}).
The difficulty to evaluate the first term of (\ref{INT}) comes from the lack
of the explicit representation for the unitary group $\cos(\lambda
t\sqrt{-\Delta_\varepsilon})$. For the classical Zakharov system
(\ref{cZK}), the explicit formula is utilized for the unitary group
$\cos(\lambda t\sqrt{-\Delta})$ and the Schr\"{o}dinger part is localized
near the origin, see \cite{OT2}.

To overcome this difficulty, we employ the method of stationary phase.
The key point is that the wave part of  (\ref{ZK}) has no stationary point
near the origin. Therefore the wave part of  (\ref{ZK}) decays faster than
$(\lambda t)^{-1}$ as $\lambda\to\infty$, which guarantees that the
interaction between $Q_{\lambda}$ and $E_{\lambda}$ are weak in the
sense that the first term in (\ref{INT}) decays like $\lambda^{-1}$ in
$H^{m}$ as $\lambda\to\infty$.

The plan of this paper is as follows. In Section 2, we prove the solvability
and uniform estimates for the solutions to (\ref{ZK}) and (\ref{4NLS}). In
Section 3, we prove Theorems \ref{main1}, \ref{main2} and \ref{main3}
and finally in Section 4, we derive the interaction estimate between
$Q_{\lambda}$ and $E_{\lambda}$.

Throughout this paper, we use the notation $A\sim B$ to represent $C_{1}A\le B\le
C_{2} A$ for some constants $C_{1}$ and $C_{2}$.
We also use the notation $A\lesssim B$ to denote $A\le CB$ for
some constant $C$.

%----------------------%     Section 2.  Preliminaries      %------------------------
\section{\bf Preliminaries}

In this section we consider the solvability and uniform estimates for the
solutions to (\ref{ZK}) and (\ref{4NLS}). 
For the quantum Zakharov system (\ref{ZK}), we have

\begin{prop}\label{well} %%%%%%%%%%%%%%%%%%%%%%%%%%%%%%%% P2.1
Let $d=1,2,3$ and $M\ge2$. Then for any $(E_{0},n_{0},n_{1})\in
{X_{M,\,d}}$ and $\lambda\in[1,\infty)$, there exists a unique solution
to (\ref{ZK}) satisfying
\begin{eqnarray*} %%%%%%%%%%%%%%%%%%%%%%%%%%%%%%%%%%%%
(E_{\lambda},n_{\lambda},\pt_{t}n_{\lambda})
&\in&C([0,\infty);\, {X_{M,\,d}}).
\end{eqnarray*} %%%%%%%%%%%%%%%%%%%%%%%%%%%%%%%%%%%%
Furthermore, for any $T\in(0,\infty)$ and some $C>0$,
$(E_{\lambda}, n_{\lambda})$ satisfies
\begin{equation} %%%%%%%%%%%%%%%%%%%%%%%%%%%%%%%%%%%%
\sup_{0\le t\le T} (\|E_{\lambda}(t)\|_{H_{x}^{M}}
+\|n_{\lambda}(t)\|_{H_{x}^{M-1}})
\le C\label{b1}
\end{equation} %%%%%%%%%%%%%%%%%%%%%%%%%%%%%%%%%%%%
for any $\lambda\in[1,\infty)$.
%, and if $d=3$, then there exists a
%$T^{\ast}>0$ depending only on $\varepsilon$ such that $(E_{\lambda},
%n_{\lambda})$ satisfies the inequality (\ref{b1}) for any
%$T\in(0,T^{\ast})$.
\end{prop} %%%%%%%%%%%%%%%%%%%%%%%%%%%%%%%%%%%% P2.1

We need the following lemmas to show Proposition \ref{well}.

%\begin{lem}[Young's inequality]\label{Y} %%%%%%%%%%%%%%%%%%%% L2.1
%Let $x,y,z\ge0$ and $p,q,r\ge0$ satisfy $q+r<2$. Then for any
%$\delta,\eta>0$,
% there exists a positive constant $C$ such that
%\begin{eqnarray*} %%%%%%%%%%%%%%%%%%%%%%%%%%%%%
%x^{p}y^{q}z^{r}\le Cx^{\frac{2p}{2-q-r}}+\delta y^{2}+\eta z^{2}
%\end{eqnarray*} %%%%%%%%%%%%%%%%%%%%%%%%%%%%%%
%holds, where $C$ depends on $\delta$ and $\eta$.
%\end{lem} %%%%%%%%%%%%%%%%%%%%%%%%%%%%%%%%%%%%

\begin{lem}[Gagliardo-Nirenberg's inequality] \label{GN} %%%%%%%%%%%%%% L2.2
Let $2\le p\le\infty$ and $j,k\in\zzz_{+}$ satisfy
$ \displaystyle %%%%%%%%%%%%%%%%%%%%%%%%%%%%%%%%
\theta:=\frac{d}{k}\left(\frac12-\frac{1}{p}+\frac{j}{d}\right)\in(0,1).
$ %%%%%%%%%%%%%%%%%%%%%%%%%%%%%%%%%%%%%%%%
Then the inequality
\begin{eqnarray*} %%%%%%%%%%%%%%%%%%%%%%%%%%%%%%%%%%%%
\|\nabla^{j}f\|_{L_{x}^{p}}
\le C\|f\|_{L_{x}^{2}}^{1-\theta}\|\nabla^{k}f\|_{L_{x}^{2}}^{\theta}
\end{eqnarray*} %%%%%%%%%%%%%%%%%%%%%%%%%%%%%%%%%%%%
holds for any $f\in H^{k}(\rre^{d})$.
\end{lem} %%%%%%%%%%%%%%%%%%%%%%%%%%%%%%%%%%%%%%%%%

The next lemma is due to Br\'{e}zis-Gallou\"{e}t \cite{BGM}
which is needed to prove Proposition \ref{well} for $d=2$.

\begin{lem}[Br\'{e}zis-Gallou\"{e}t inequality] \label{BG} %%%%%%%%%%%%%% L2.3
We have the following inequalities.

\vskip1mm \noindent

(i) Let $f\in H^{2}(\rre^{2})$. Then we have $f\in L^{\infty}(\rre^{2})$ and
\begin{eqnarray*} %%%%%%%%%%%%%%%%%%%%%%%%%%%%%%%%%%%%
\|f\|_{L_{x}^{\infty}}\le C(\|f\|_{H_{x}^{1}}
\sqrt{\log(e+\|\nabla^{2}f\|_{L_{x}^{2}})}+1).
\end{eqnarray*} %%%%%%%%%%%%%%%%%%%%%%%%%%%%%%%%%%%%

\vskip1mm \noindent

(ii) Let $f\in H^{1}(\rre^{2})$. Then we have $f\in L^{4}(\rre^{2})$ and
\begin{eqnarray*} %%%%%%%%%%%%%%%%%%%%%%%%%%%%%%%%%%%%
\|f\|_{L_{x}^{4}}\le C(\|f\|_{L_{x}^{2}}^{1/2}
\|\nabla f\|_{L_{x}^{2}}^{1/2}
\sqrt{\log(e+\|\nabla f\|_{L_{x}^{2}})}+1).
\end{eqnarray*} %%%%%%%%%%%%%%%%%%%%%%%%%%%%%%%%%%%%
\end{lem}%%%%%%%%%%%%%%%%%%%%%%%%%%%%%%%%%%%%%%%%%

\vskip3mm \noindent
{\bf Proof of Lemma \ref{BG}.} %%%%%%%%%%%%%%%%%%%%%%%%%%
See \cite[Lemma 2]{BGM}. $\qed$

\vskip2mm

To prove Proposition \ref{well} for $d=3$, we employ the
Strichartz estimates for the Schr\"odinger equations
%\begin{equation} %%%%%%%%%%%%%%%%%%%%%%%%%%%%%%%%
%(i\partial_t -\Delta)u = h, \qquad u(0,\cdot) = u_0, \qquad x\in\mathbb{R}^3,
%\label{Sch} \end{equation}  %%%%%%%%%%%%%%%%% %%%%%%%%%%%%%%%
and the fouth-order Schr\"odinger equations
\begin{equation}
%\begin{eqnarray}
\left\{
\begin{array}{l}
i\partial_tE-\Delta E+\Delta^2 E = h(t,x)
\qquad t\in\rre,\ x\in\rre^{3},\\
E(0,x)=E_{0}(x)
\qquad\qquad\qquad\quad\ t\in\rre,\ x\in\rre^{3}.
\end{array}
\right. \label{4-Sch}
%\end{eqnarray}
\end{equation} %%%%%%%%%%%%%%%%% %%%%%%%%%%%%%%%
A pair $(q, r)$ is called (3 dimensional) Schr\"odinger admissible,
% for short $S$-admissible,
if $2\leq q, \, r\leq \infty$, %$(q,\, r,\, d)\ne (2,\, \infty,\,2)$,
and
\begin{equation} %%%%%%%%%%%%%%%%%%%%%%%%
\frac{2}{q}+\frac{3}{r}=\frac{3}{2}.
 \label{S-qrd} \end{equation}  %%%%%%%%%%%%%%%%%
A pair $(q, r)$ is called (3 dimensional) biharmionic admissible,
% for short $B$-admissible,
if $2\leq q, \, r\leq \infty$, %$(q,\, r,\, d)\ne (2,\, \infty,\,4)$,
and
\begin{equation} %%%%%%%%%%%%%%%%%%%%%%%%
\frac{4}{q}+\frac{3}{r}=\frac{3}{2}.
 \label{B-qrd} \end{equation}  %%%%%%%%%%%%%%%%%
We recall some known results.
%\begin{lem}\label{Str-est} %%%%%%%%%%%%%%%%%%%%%%%%%%%%%%%% L5.1
%(Keel and Tao \cite{KT}) Suppose that %$d\geq1$ and
%$(q, r)$ and
%$(\widetilde q, \widetilde r)$ are Schr\"odinger admissible. If $u$ is a
%solution of \eqref{Sch} for some data $u_0$, $h$ and time $0<T<\infty$,
%%then
%\begin{equation} %%%%%%%%%%%%%%%%%%%%%%%%%%%%%%%%%%%%
%\|u\|_{L_t^q([0,T]; L_x^r(\rre^{3}))}+ \|u\|_{C([0,T]; L_x^2(\rre^{3}))} \lesssim
%\|u_0\|_{L_x^2(\rre^{3})}+ \|h\|_{L_t^{\widetilde q'}([0,T]; L_x^{\widetilde r'}(\rre^{3}))},
%\label{Str-1} \end{equation} %%%%%%%%%%%%%%%%%%%%%%%%%%%%%%%%%%%%
%\end{lem} %%%%%%%%%%%%%%%%%%%%%%%%%%%%%%%%%%%% L5.1
\begin{lem}\label{4-Str-est} %%%%%%%%%%%%%%%%%%%%%%%%%%%%%%%% L5.2
(Pausader \cite{P}) Let $E\in C([0,T], H^{-4}(\rre^{3}))$ be a solution of
\eqref{4-Sch}. For any biharmonic pairs  $(q, r)$ and $(\widetilde q,
\widetilde r)$, it satisfies
\begin{equation} %%%%%%%%%%%%%%%%%%%%%%%%%%%%%%%%%%%%
\|E\|_{L_t^q([0,T]; L_x^r(\rre^{3}))} \leq
C\Big(\|E_0\|_{L_x^2(\rre^{3})}+
\|h\|_{L_t^{\widetilde q'}([0,T]; L_x^{\widetilde r'}(\rre^{3}))}\Big),
\label{Str-2} \end{equation} %%%%%%%%%%%%%%%%%%%%%%%%%%%%%%%%%%%%
where $C$ depends only on $\widetilde q'$, and $\widetilde r'$.
Besides, for any Schr\"odinger admissible pairs $(q, r)$ and $(a, b)$, and
any $s\geq0$, we have 
\begin{equation} %%%%%%%%%%%%%%%%%%%%%%%%%%%%%%%%%%%%
\||\nabla|^{s} E\|_{L_t^q([0,T]; L_x^r(\rre^{3}))} \leq
C\Big(\||\nabla|^{s-\frac{2}{q}} E_0\|_{L_x^2(\rre^{3})}+
\||\nabla|^{s-\frac{2}{q}-\frac{2}{a}} h\|_{L_t^{a'}([0,T]; L_x^{b'}(\rre^{3}))}\Big).
\label{Str-3} \end{equation} %%%%%%%%%%%%%%%%%%%%%%%%%%%%%%%%%%%%
where $C$ depends only on $a'$, and $b'$.
\end{lem} %%%%%%%%%%%%%%%%%%%%%%%%%%%%%%%%%%%% L5.2

\vskip3mm \noindent
{\bf Proof of Proposition \ref{well}.} %%%%%%%%%%%%%%%%%%%%%%%%%%
We prove this proposition by the induction argument on $M$. For the
simplicity, we abbreviate $(E_{\lambda},n_{\lambda})$ to $(E,n)$. By the
density $C_{0}^{\infty}(\rre^{d})\hookrightarrow H^{M}(\rre^{d})$, we
may assume that $(E,n)$ is smooth. The global existence and uniqueness
of the solution to (\ref{ZK}) is proved by Guo-Zhang-Guo \cite[Theorem
1.1]{GZG}. Hence we derive the uniform bound for solution (\ref{b1}).

We first derive the $L^{2}$ bound of $E$. Taking the imaginary part of
the inner product in $L_{x}^{2}$ between the first equation of (\ref{ZK})
and $E$, we have
$%\begin{eqnarray*} %%%%%%%%%%%%%%%%%%%%%%%%%%%%%%%%%%%%
\frac{d}{dt}\|E(t)\|_{L_{x}^{2}}^{2} = 0.
$ %\end{eqnarray*} %%%%%%%%%%%%%%%%%%%%%%%%%%%%%%%%%%%%
Therefore
\begin{equation}
%\begin{eqnarray} %%%%%%%%%%%%%%%%%%%%%%%%%%%%%%%%%%%%
\|E(t)\|_{L_{x}^{2}}=\|E_{0}\|_{L_{x}^{2}} \label{f30}
%\end{eqnarray} %%%%%%%%%%%%%%%%%%%%%%%%%%%%%%%%%%%%
\end{equation}
for any $t\in(0,\infty)$.

Next we derive the $H^{2}\times H^{1}$ bound for $(E,n)$. Taking the
real part of the inner product in $L_{x}^{2}$ between the first equation
of (\ref{ZK}) and $\pt_{t}E$, we obtain
\begin{equation}
%\begin{eqnarray} %%%%%%%%%%%%%%%%%%%%%%%%%%%%%%%%%%%%
\frac{d}{dt}
(\|\nabla E(t)\|_{L_{x}^{2}}^{2}
+\varepsilon^{2}
\|\Delta E(t)\|_{L_{x}^{2}}^{2})
=-2Re
\int_{\rre^{d}}n\overline{E}\pt_{t}Edx.
\label{f31}
%\end{eqnarray} %%%%%%%%%%%%%%%%%%%%%%%%%%%%%%%%%%%%
\end{equation}
Applying the operator $\nabla^{-1}$ to the second equation of (\ref{ZK})
and taking the inner product in $L_{x}^{2}$ between the resultant
equation and $\nabla^{-1}\pt_{t}n$, we have
\begin{equation}
%\begin{eqnarray} %%%%%%%%%%%%%%%%%%%%%%%%%%%%%%%%%%%%
\frac{d}{dt}
\left(\frac12\lambda^{-2}
\|\pt_{t}\nabla^{-1}n(t)\|_{L_{x}^{2}}^{2}
+\frac12\|n(t)\|_{L_{x}^{2}}^{2}
+\frac{\varepsilon^{2}}{2}\|
\nabla n(t)\|_{L_{x}^{2}}^{2}
\right)
=-Re\int_{\rre^{d}}\pt_{t}n|E|^{2}
dx.\label{f32}
%\end{eqnarray} %%%%%%%%%%%%%%%%%%%%%%%%%%%%%%%%%%%%
\end{equation}
%where we used the following calculation. Since $\Delta
%n=-\text{div}(-\nabla n)$ and $-\nabla n\in G$, we have
%$\nabla^{-1}\Delta n=-\nabla n$. On the other hand, $\pt_{t}n$ is
%decomposed as $\pt_{t}n=-\text{div}{{\mathbb F}}$, ${{\mathbb
%F}}=\nabla g+{{\mathbb H}}$ with $\text{div}{{\mathbb H}}=0$. Hence
%$\nabla^{-1}\pt_{t}n=\nabla g$. Therefore, integrating by parts, we have
%\begin{eqnarray*} %%%%%%%%%%%%%%%%%%%%%%%%%%%%%%%%%%%%
%\int_{\rre^{d}}\nabla^{-1}\Delta n\nabla^{-1}\pt_{t}ndx
%=-\int_{\rre^{d}}\nabla n\nabla gdx=\int_{\rre^{d}}n\Delta gdx.
%\end{eqnarray*} %%%%%%%%%%%%%%%%%%%%%%%%%%%%%%%%%%%%
%Noting $\pt_{t}n=-\Delta g$, we obtain
%\begin{eqnarray*} %%%%%%%%%%%%%%%%%%%%%%%%%%%%%%%%%%%%
%\int_{\rre^{d}}\nabla^{-1}\Delta n\nabla^{-1}\pt_{t}ndx
%=-\frac12\frac{d}{dt}\| n(t)\|_{L_{x}^{2}}^{2}.
%\end{eqnarray*} %%%%%%%%%%%%%%%%%%%%%%%%%%%%%%%%%%%%
%In a similar way, we have
%\begin{eqnarray*} %%%%%%%%%%%%%%%%%%%%%%%%%%%%%%%%%%%%
%\int_{\rre^{d}}\nabla^{-1}\Delta^{2} n\nabla^{-1}\pt_{t}ndx
%=-\frac12\frac{d}{dt}\|\nabla n(t)\|_{L_{x}^{2}}^{2}.
%\end{eqnarray*} %%%%%%%%%%%%%%%%%%%%%%%%%%%%%%%%%%%%
Let
\begin{eqnarray*} %%%%%%%%%%%%%%%%%%%%%%%%%%%%%%%%%%%%
H_{2}(E,n)(t)
&=&\|\nabla E(t)\|_{L_{x}^{2}}^{2}
+\varepsilon^{2}
\|\Delta E(t)\|_{L_{x}^{2}}^{2}
+\frac12\lambda^{-2}\|\pt_{t}\nabla^{-1}n(t)\|_{L_{x}^{2}}^{2}\\
& &
+\frac12\|n(t)\|_{L_{x}^{2}}^{2}
+\frac{\varepsilon^{2}}{2}\|
\nabla n(t)\|_{L_{x}^{2}}^{2}
+\int_{\rre^{d}}n|E|^{2}dx
+C_2,
\end{eqnarray*} %%%%%%%%%%%%%%%%%%%%%%%%%%%%%%%%%%%%
where $C_2=C^\prime_{2}\|E_{0}\|_{L_{x}^{2}}^{\frac{16-2d}{4-d}}$. The
positive constant $C^\prime_{2}$ is chosen so that the inequality
\begin{equation}
%\begin{eqnarray} %%%%%%%%%%%%%%%%%%%%%%%%%%%%%%%%%%%%
\|E(t)\|_{H_{x}^{2}}^{2}
+\|n(t)\|_{H_{x}^{1}}^{2}
\le CH_{2}(E,n)(t)
\label{v2}
%\end{eqnarray} %%%%%%%%%%%%%%%%%%%%%%%%%%%%%%%%%%%%
\end{equation}
holds for some constant $C$ independent of $\lambda\in[1,\infty)$.
Indeed, Lemma \ref{GN} and the Young inequality yield
\begin{eqnarray} %%%%%%%%%%%%%%%%%%%%%%%%%%%%%%%%%%%%
\int_{\rre^{d}}n|E|^{2} dx
%&\le&\|n\|_{L_{x}^{2}}\|E\|_{L_{x}^{4}}^{2}\nonumber\\
&\le&C'
\|n\|_{L_{x}^{2}}\|E\|_{L_{x}^{2}}^{2-\frac{d}{4}}
\|\Delta E\|_{L_{x}^{2}}^{\frac{d}{4}}
\nonumber\\
%&\le&C'
%\|E\|_{L_{x}^{2}}^{\frac{16-2d}{4-d}}
%+\frac14\|n\|_{L_{x}^{2}}^{2}
%+\frac{\varepsilon^{2}}{2}
%\|\Delta E\|_{L_{x}^{2}}^{2}
%\nonumber\\
&\le&\frac{\varepsilon^{2}}{2}
\|\Delta E\|_{L_{x}^{2}}^{2}
+\frac14\|n\|_{L_{x}^{2}}^{2}
+C'\|E_{0}\|_{L_{x}^{2}}^{\frac{16-2d}{4-d}}.
\label{v1}
\end{eqnarray} %%%%%%%%%%%%%%%%%%%%%%%%%%%%%%%%%%%%
Choosing $C^\prime_{2}=C'+1$, we have
\begin{eqnarray*} %%%%%%%%%%%%%%%%%%%%%%%%%%%%%%%%%%%%
\|E(t)\|_{H_{x}^{2}}^{2}
+\|n(t)\|_{H_{x}^{1}}^{2}
\le \frac{4}{\varepsilon^{2}}H_{2}(E,n)(t).
\end{eqnarray*} %%%%%%%%%%%%%%%%%%%%%%%%%%%%%%%%%%%%
From (\ref{f30}), (\ref{f31}) and (\ref{f32}), we get
$%\begin{eqnarray*} %%%%%%%%%%%%%%%%%%%%%%%%%%%%%%%%%%%%
\frac{d}{dt}H_{2}(E,n)(t)=0.
$ %\end{eqnarray*} %%%%%%%%%%%%%%%%%%%%%%%%%%%%%%%%%%%%
Hence
\begin{eqnarray} %%%%%%%%%%%%%%%%%%%%%%%%%%%%%%%%%%%%
H_{2}(E,n)(t)=H_{2}(E,n)(0).
\label{v3}
\end{eqnarray} %%%%%%%%%%%%%%%%%%%%%%%%%%%%%%%%%%%%
Furthermore from (\ref{v1}), we find
\begin{eqnarray} %%%%%%%%%%%%%%%%%%%%%%%%%%%%%%%%%%%%
H_{2}(E,n)(0)
\le C(
\|E_{0}\|_{L_{x}^{2}}^{\frac{16-2d}{4-d}}
+
\|E_{0}\|_{H_{x}^{2}}^{2}
+\|n_{0}\|_{H_{x}^{1}}^{2}
+\|n_{1}\|_{\dot{H}_{x}^{-1}}^{2}),
\label{v4}
\end{eqnarray} %%%%%%%%%%%%%%%%%%%%%%%%%%%%%%%%%%%%
where the constant $C$ is independent of $\lambda\in[1,\infty)$.
Combining (\ref{v2}), (\ref{v3}) and (\ref{v4}), we have
\begin{eqnarray*} %%%%%%%%%%%%%%%%%%%%%%%%%%%%%%%%%%%%
\|E(t)\|_{H_{x}^{2}}^{2}
+\|n(t)\|_{H_{x}^{1}}^{2}
\le
C(\|E_{0}\|_{L_{x}^{2}}^{\frac{16-2d}{4-d}}
+\|E_{0}\|_{H_{x}^{2}}^{2}
+\|n_{0}\|_{H_{x}^{1}}^{2}
+\|n_{1}\|_{\dot{H}_{x}^{-1}}^{2})
\end{eqnarray*} %%%%%%%%%%%%%%%%%%%%%%%%%%%%%%%%%%%%
for any $t\in[0,\infty)$. This proves $H^2\times H^1$ bound for $(E,n)$.

\vskip1mm

Finally we derive that for any integers $M\ge3$, $H^{M-1}\times
H^{M-2}$ bound for $(E,n)$ implies $H^{M}\times H^{M-1}$ bound for
$(E,n)$. Now we assume that
\begin{equation} %%%%%%%%%%%%%%%%%%%%%%%%%%%%%%%%%%%%
\sup_{0\leq t \leq T}
\big(\|E(t)\|_{H^{M-1}_x}+\|n(t)\|_{H^{M-2}_x}\big) \leq C. \label{Induction}
\end{equation} %%%%%%%%%%%%%%%%%%%%%%%%%%%%%%%%%%%%
Applying $\nabla^{M-2}$ to the first equation of (\ref{ZK}) and taking the
real part of the inner product in $L_{x}^{2}$ between the resulting
equation and $\pt_{t}\nabla^{M-2}E$, we obtain
\begin{eqnarray} %%%%%%%%%%%%%%%%%%%%%%%%%%%%%%%%%%%%
\lefteqn{\frac{d}{dt}
\left(\|\nabla^{M-1}E(t)\|_{L_{x}^{2}}^{2}
+\varepsilon^{2}
\|\nabla^{M}E(t)\|_{L_{x}^{2}}^{2}
\right)}\nonumber\\
&=&-2Re
\int_{\rre^{d}}
\nabla^{M-2}(n\overline{E})
\pt_{t}\nabla^{M-2}Edx.
\label{f33}
\end{eqnarray} %%%%%%%%%%%%%%%%%%%%%%%%%%%%%%%%%%%%
Similarly, applying $\nabla^{M-3}$ to the second equation of (\ref{ZK})
and taking the inner product in $L_{x}^{2}$ between the resulting
equation and $\pt_{t}\nabla^{M-3}n$, we find
\begin{eqnarray} %%%%%%%%%%%%%%%%%%%%%%%%%%%%%%%%%%%%
\lefteqn{\hskip-9mm\frac{d}{dt}
\left(\frac12\lambda^{-2}
\|\pt_{t}\nabla^{M-3}n(t)\|_{L_{x}^{2}}^{2}
+\frac12\|\nabla^{M-2}n(t)\|_{L_{x}^{2}}^{2}
+\frac{\varepsilon^{2}}{2}\|
\nabla^{M-1}n(t)\|_{L_{x}^{2}}^{2}
\right)}
\nonumber\\
&=&-\int_{\rre^{d}}
\nabla^{M-2}|E|^{2}
\pt_{t}\nabla^{M-2}ndx.
\qquad\qquad\qquad\qquad{}
\label{f34}
\end{eqnarray} %%%%%%%%%%%%%%%%%%%%%%%%%%%%%%%%%%%%
Define
\begin{eqnarray} %%%%%%%%%%%%%%%%%%%%%%%%%%%%%%%%%%%%
H_{M}(E,n)(t)
&=&\|\nabla^{M-1}E(t)
\|_{L_{x}^{2}}^{2}
+\varepsilon^{2}\|\nabla^{M}E(t)\|_{L_{x}^{2}}^{2}
+\frac12\lambda^{-2}\|\pt_{t}\nabla^{M-3}n(t)\|_{L_{x}^{2}}^{2}\nonumber\\
& &+
\frac12\|\nabla^{M-2}n(t)\|_{L_{x}^{2}}^{2}
+\frac{\varepsilon^{2}}{2}\|\nabla^{M-1} n(t)\|_{L_{x}^{2}}^{2}\\
& &+\int_{\rre^{d}}\nabla^{M-2}n\nabla^{M-2}|E|^{2}dx+C_{M}.\nonumber
\end{eqnarray} %%%%%%%%%%%%%%%%%%%%%%%%%%%%%%%%%%%%
In fact, the positive constant $C_{M}$ is chosen in a way such that
$H_{M}(E,n)(t)$ satisfies
\begin{eqnarray*} %%%%%%%%%%%%%%%%%%%%%%%%%%%%%%%%%%%%
\|\nabla^{M}E(t)\|_{L_{x}^{2}}^{2}
+\|\nabla^{M-1}n(t)\|_{L_{x}^{2}}^{2}+1
\le CH_{M}(E,n)(t)
\end{eqnarray*} %%%%%%%%%%%%%%%%%%%%%%%%%%%%%%%%%%%%
for any $t\in[0,T)$, where $C$ is independent of $\lambda\in[1,\infty)$.
Indeed, by the Sobolev embedding $H^{1}\hookrightarrow L^{4}$, we
have
\begin{eqnarray*} %%%%%%%%%%%%%%%%%%%%%%%%%%%%%%%%%%%%
\int_{\rre^{d}}
\nabla^{M-2}n\nabla^{M-2}|E|^{2}dx
%\le C'\|\nabla^{M-2}n\|_{L_{x}^{2}} \|E\|_{W_{x}^{M-2,4}}^{2}
\le
C'\|\nabla^{M-2}n\|_{L_{x}^{2}}
\|E\|_{H_{x}^{M-1}}^{2},
\end{eqnarray*} %%%%%%%%%%%%%%%%%%%%%%%%%%%%%%%%%%%%
hence it suffices to choose $C_{M}=C' \sup_{0\le t\le
T}\|n(t)\|_{H_{x}^{M-2}}\sup_{0\le t\le T}\|E(t)\|^{2}_{H_{x}^{M-1}}
+1$. Note that $C_{M}$ is independent of $\lambda\in[1,\infty)$.

From (\ref{f33}) and (\ref{f34}), we obtain
\begin{eqnarray*} %%%%%%%%%%%%%%%%%%%%%%%%%%%%%%%%%%%%
\lefteqn{\frac{d}{dt}H_{M}(E,n)(t)}\\
&=&
\int_{\rre^{d}}\nabla^{M-2}
\pt_{t}|E|^{2}\nabla^{M-2}ndx
-2Re\int_{\rre^{d}}\nabla^{M-2}
(n\overline{E})\pt_{t}\nabla^{M-2}
Edx.
\end{eqnarray*} %%%%%%%%%%%%%%%%%%%%%%%%%%%%%%%%%%%%
From (\ref{ZK}) we see $\pt_{t}|E|^{2} =-2Im(\overline{E}\Delta E)
+2\varepsilon^{2} Im (\overline{E}\Delta^{2}E)$. Hence we have
\begin{equation} %%%%%%%%%%%%%%%%%%%%%%%%%%%%%%%%%%%%
\frac{d}{dt}H_{M}(E,n)(t)\equiv I_{1}+I_{2}+I_{3}+I_{4},
\end{equation} %%%%%%%%%%%%%%%%%%%%%%%%%%%%%%%%%%%%
where
\begin{equation}\begin{array}{lll} %%%%%%%%%%%%%%%%%%%%%%%%%%%%%%%%%%%%
&\displaystyle I_{1}=-2Im\int_{\rre^{d}}\nabla^{M-2}
(\overline{E}\Delta E)\nabla^{M-2}n\,dx,\\
&\displaystyle I_{2}=2\varepsilon^{2}Im\int_{\rre^{d}}\nabla^{M-2}
(\overline{E}\Delta^{2} E)\nabla^{M-2}n\,dx,\\
&\displaystyle I_{3}=+2Im\int_{\rre^{d}}\nabla^{M-2}
(n\overline{E})\nabla^{M}E\,dx,\\
&\displaystyle I_{4}=-2\varepsilon^{2} Im\int_{\rre^{d}}\nabla^{M-2}
(n\overline{E})\nabla^{M+2}E\,dx.
\end{array}\label{I1234}\end{equation} %%%%%%%%%%%%%%%%%%%%%%%%%%%%%%%%%%%%
For $I_{1}$ and $I_{3}$, we use the Sobolev inequality and the induction
hypothesis \eqref{Induction} to obtain
\begin{eqnarray} %%%%%%%%%%%%%%%%%%%%%%%%%%%%%%%%%%%%
|I_{1}|
%&\le&C\|E(t)\|_{W_{x}^{M-2,\infty}}
%\|E(t)\|_{H_{x}^{M}}\|\nabla^{M-2}n(t)\|_{L_{x}^{2}}\\
&\le&
C\|E(t)\|_{H_{x}^{M}}^{2}
\|\nabla^{M-2}n(t)\|_{L_{x}^{2}}\nonumber\\
%&\le&C(\|\nabla^{M}E(t)\|_{L_{x}^{2}}
%+\|E(t)\|_{H_{x}^{M-1}})^{2}\|n(t)\|_{H_{x}^{M-2}}\\
&\le&
C(\|\nabla^{M}E(t)\|_{L_{x}^{2}}^{2}+1),\label{q1}\\
|I_{3}|
&\le&
\|n(t)\|_{H_{x}^{M-2}}
\|E(t)\|_{W_{x}^{M-2,\infty}}
\|\nabla^{M}E(t)\|_{L_{x}^{2}}\nonumber\\
%&\le&C\|n(t)\|_{H_{x}^{M-2}}
%\|E(t)\|_{H_{x}^{M}}\|\nabla^{M}E(t)\|_{L_{x}^{2}}\\
&\le&
C(\|\nabla^{M}E(t)\|_{L_{x}^{2}}+1)
\|\nabla^{M}E(t)\|_{L_{x}^{2}}.\label{q3}
\end{eqnarray} %%%%%%%%%%%%%%%%%%%%%%%%%%%%%%%%%%%%
An integration by parts leads
\begin{eqnarray*} %%%%%%%%%%%%%%%%%%%%%%%%%%%%%%%%%%%%
I_{2}
%&=&-2\varepsilon^{2}Im\int_{\rre^{d}}\nabla^{M-3}
%(\overline{E}\Delta^{2} E)\nabla^{M-1}ndx\\
&=&-2\varepsilon^{2}Im\int_{\rre^{d}}
\overline{E}\nabla^{M+1}E\nabla^{M-1}ndx\\
& &-2\varepsilon^{2}
{{\bf 1}}_{M\ge4}(M)
\sum_{j=0}^{M-4}
\bmat{c}M-3\\j\emat
Im\int_{\rre^{d}}\nabla^{M-3-j}
\overline{E}\nabla^{j+4}E\nabla^{M-1}ndx,
\end{eqnarray*} %%%%%%%%%%%%%%%%%%%%%%%%%%%%%%%%%%%%
where ${{\bf 1}}_{M\ge4}(M)=0$ for $M\le3$ and ${{\bf
1}}_{M\ge4}(M)=1$ for $M\ge4$. Hence
\begin{eqnarray} %%%%%%%%%%%%%%%%%%%%%%%%%%%%%%%%%%%%
\lefteqn{\left|
I_{2}+2\varepsilon^{2}Im\int_{\rre^{d}}
\overline{E}\nabla^{M+1}E\nabla^{M-1}ndx\right|}
\nonumber\\
&\le&
C\|E(t)\|_{W_{x}^{M-3,\infty}}
\|E(t)\|_{H_{x}^{M}}
\|\nabla^{M-1}n\|_{L_{x}^2}
\nonumber\\
%&\le&
%C\|E(t)\|_{H_{x}^{M-1}}
%\|E(t)\|_{H_{x}^{M}}
%\|\nabla^{M-1}n\|_{L_{x}^2}\\
%&\le&
%C\|E(t)\|_{H_{x}^{M-1}}
%(\|\nabla^{M}E(t)\|_{L_{x}^{2}}+
%\|E(t)\|_{H_{x}^{M-1}})
%\|\nabla^{M-1}n\|_{L_{x}^2}\\
&\le&C(\|\nabla^{M}E(t)\|_{L_{x}^{2}}+1)
\|\nabla^{M-1}n\|_{L_{x}^2}.
\label{q2}
\end{eqnarray} %%%%%%%%%%%%%%%%%%%%%%%%%%%%%%%%%%%%
Again an integration by parts yields
\begin{eqnarray} %%%%%%%%%%%%%%%%%%%%%%%%%%%%%%%%%%%%
I_{4}&=&
%2\varepsilon^{2}
%Im\int_{\rre^{d}}\nabla^{M-1}
%(n\overline{E})\nabla^{M+1}Edx\\
%&=&
%2\varepsilon^{2}
%Im\int_{\rre^{d}}
%\overline{E}\nabla^{M-1}n\nabla^{M+1}Edx\\
%& &+
%2\varepsilon^{2}
%\sum_{j=1}^{M-1}
%\bmat{c}M-1\\j\emat
%Im\int_{\rre^{d}}
%\nabla^{M-1-j}n\nabla^{j}\overline{E}\nabla^{M+1}Edx\\
%&=&
2\varepsilon^{2}
Im\int_{\rre^{d}}
\overline{E}\nabla^{M-1}n\nabla^{M+1}Edx
\nonumber\\
& &-
2\varepsilon^{2}
\sum_{j=1}^{M-1}\bmat{c}M-1\\j\emat
Im\int_{\rre^{d}}
\nabla^{M-j}n\nabla^{j}\overline{E}\nabla^{M}Edx
\nonumber\\
& &-
2\varepsilon^{2}
\sum_{j=1}^{M-1}\bmat{c}M-1\\j\emat
Im\int_{\rre^{d}}
\nabla^{M-1-j}n\nabla^{j+1}\overline{E}\nabla^{M}Edx.
\label{q4}
\end{eqnarray} %%%%%%%%%%%%%%%%%%%%%%%%%%%%%%%%%%%%

Now we discuss the cases $d=1,2,3$ separately.

\vskip2mm \noindent \underline{Case: $d=1$.} By the embedding
$H^{1}(\rre)\hookrightarrow L^{\infty}(\rre)$, we find
\begin{eqnarray*} %%%%%%%%%%%%%%%%%%%%%%%%%%%%%%%%%%%%
\lefteqn{
\left|I_{4}-2\varepsilon^{2}
Im\int_{\rre^{d}}
\overline{E}\nabla^{M-1}n\nabla^{M+1}Edx\right|}\\
&\le&
C\sum_{j=1}^{2}\|\nabla^{M-j}n(t)
\|_{L_{x}^{2}}
\|\nabla^{j} E(t)\|_{L_{x}^{\infty}}
\|\nabla^{M}E(t)\|_{L_{x}^{2}}\\
& &
+C\sum_{j=3}^{M}
\|\nabla^{M-j}n(t)\|_{L_{x}^{\infty}}
\|\nabla^{j} E(t)\|_{L_{x}^{2}}
\|\nabla^{M}E(t)\|_{L_{x}^{2}}\\
%&\le&
%C\|\nabla^{M-1}
%n(t)\|_{L_{x}^{2}}\|E(t)\|_{H_{x}^{2}}
%\|\nabla^{M}E(t)\|_{L_{x}^{2}}\\
%& &+
%C\|n(t)\|_{H_{x}^{M-2}}\|E(t)\|_{H_{x}^{M}}
%\|\nabla^{M}E(t)\|_{L_{x}^{2}}\\
&\le&
C(1+\|\nabla^{M}E(t)\|_{L_{x}^{2}}
+\|\nabla^{M-1}n(t)\|_{L_{x}^{2}})
\|\nabla^{M}E(t)\|_{L_{x}^{2}}.
\end{eqnarray*} %%%%%%%%%%%%%%%%%%%%%%%%%%%%%%%%%%%%
Combining the above inequality, (\ref{q1}),
(\ref{q3}) and (\ref{q2}), we have
\begin{equation*}  %%%%%%%%%%%%%%%%%%%%%%%%%%%%%%%%%%%%
\frac{d}{dt}
H_{M}(E,n)(t)
\le
CH_{M}(E,n)(t).
\end{equation*} %%%%%%%%%%%%%%%%%%%%%%%%%%%%%%%%%%%%
The Gronwall lemma yields
\begin{equation*} %%%%%%%%%%%%%%%%%%%%%%%%%%%%%%%%%%%%
H_{M}(E,n)(t)
\le CH_{M}(E,n)(0)e^{Ct}.
\end{equation*} %%%%%%%%%%%%%%%%%%%%%%%%%%%%%%%%%%%%
Since
\begin{eqnarray} %%%%%%%%%%%%%%%%%%%%%%%%%%%%%%%%%%%%
\lefteqn{\|\nabla^{M}E_{0}\|_{L_{x}^{2}}^{2}
+
\|\nabla^{M-1}n_{0}\|_{L_{x}^{2}}^{2}
\le H_{M}(E,n)(0)}\nonumber\\
&\le&
C(\|E_{0}\|_{H_{x}^{M}}
+\|n_{0}\|_{H_{x}^{M-1}}
+\|n_{1}\|_{H_{x}^{M-3}}
+\|E_{0}\|_{H_{x}^{M-1}}^{2}
\|n_{0}\|_{H_{x}^{M-2}}+1),
\nonumber\\
\label{bound}
\end{eqnarray} %%%%%%%%%%%%%%%%%%%%%%%%%%%%%%%%%%%%
we have $H^{M}\times H^{M-1}$ bound for $(E,n)$ with $d=1$.

\vskip2mm

\noindent \underline{Case: $d=2$.} By the Br\'{e}zis-Gallou\"{e}t
inequality (Lemma \ref{BG}), we have
\begin{eqnarray*} %%%%%%%%%%%%%%%%%%%%%%%%%%%%%%%%%%%%
\lefteqn{
\left|I_{4}-2\varepsilon^{2}
Im\int_{\rre^{d}}
\overline{E}\nabla^{M-1}n\nabla^{M+1}Edx\right|}\\
%&\le&
%C\|\nabla^{M-1}n(t)\|_{L_{x}^{2}}
%\|\nabla E(t)\|_{L_{x}^{\infty}}
%\|\nabla^{M}E(t)\|_{L_{x}^{2}}\\
%& &+
%C\|\nabla^{M-2}n(t)\|_{L_{x}^{4}}
%\|\nabla^{2} E(t)\|_{L_{x}^{4}}
%\|\nabla^{M}E(t)\|_{L_{x}^{2}}\\
%& &
%+C\sum_{j=3}^{M}
%\|\nabla^{M-j}n(t)\|_{L_{x}^{\infty}}
%\|\nabla^{j} E(t)\|_{L_{x}^{2}}
%\|\nabla^{M}E(t)\|_{L_{x}^{2}}\\
&\le&
C\|n(t)\|_{H_{x}^{M-1}}\|E(t)\|_{H_{x}^{M}}
+\Big\{\|E(t)\|_{H_{x}^{2}}
\sqrt{\log(e+\|\nabla^{3}E(t)\|_{L_{x}^{2}})}+1\Big\}
\\
& &+C\|E(t)\|_{H_{x}^{M}}\,
\Big\{\|\nabla^{M-2}n(t)\|_{L_{x}^{2}}^{1/2}\|\nabla^{M-1}n(t)
\|_{L_{x}^{2}}^{1/2}\\
& &\qquad\qquad\qquad\qquad\times
\sqrt{\log(e+\|\nabla^{M-1}n(t)\|_{L_{x}^{2}})}+1\Big\}\\
& &\times
\Big\{\|\nabla^{2}E(t)\|_{L_{x}^{2}}^{1/2}\|\nabla^{3}E(t)
\|_{L_{x}^{2}}^{1/2}
\sqrt{\log(e+\|\nabla^{3}E(t)\|_{L_{x}^{2}})}+1\Big\}.
\end{eqnarray*} %%%%%%%%%%%%%%%%%%%%%%%%%%%%%%%%%%%%
Combining the above inequality,(\ref{q1}),
(\ref{q3}) and (\ref{q2}), we have
\begin{equation*}  %%%%%%%%%%%%%%%%%%%%%%%%%%%%%%%%%%%%
\frac{d}{dt}
H_{M}(E,n)(t)
\le CH_{M}(E,n)(t)\log(H_{M}(E,n)).
\end{equation*} %%%%%%%%%%%%%%%%%%%%%%%%%%%%%%%%%%%%
The Gronwall lemma yields
\begin{equation*} %%%%%%%%%%%%%%%%%%%%%%%%%%%%%%%%%%%%
H_{M}(E,n)(t)
\le C(H_{M}(E,n)(0))^{e^{ct}}.
\end{equation*} %%%%%%%%%%%%%%%%%%%%%%%%%%%%%%%%%%%%
By the above inequality and (\ref{bound}),
we have $H^{M}\times H^{M-1}$ bound for $(E,n)$ with $d=2$.

\vskip2mm

\noindent \underline{Case: $d=3$.}
The conservation of mass and conservation of
Hamiltonian (\ref{Hami}) imply that
\begin{equation} %%%%%%%%%%%%%%%%%%%%%%%%%%%%%%%%
\sup_{0\le t\le T}\left(\|E(t)\|_{H^2}+\|\mathcal{N}(t)\|_{H^{1}}\right)\le C,
 \label{EN1} \end{equation} %%%%%%%%%%%%%%%%% %%%%%%%%%%%%%%%%
where
\begin{equation*} %%%%%%%%%%%%%%%%%%%%%%%%%%%%%%%%
\mathcal{N}=n+\frac{i}{\sqrt{(1-\varepsilon^2\Delta)(-\Delta)}}
\frac{\pt_{t} n}{\lambda}.
\end{equation*} %%%%%%%%%%%%%%%%% %%%%%%%%%%%%%%%%
%Denote
%\begin{equation} %%%%%%%%%%%%%%%%%%%%%%%%%%%%%%%%
%I_\varepsilon=(I-\varepsilon^2\Delta)^{-1}.
% \label{EN2} \end{equation} %%%%%%%%%%%%%%%%% %%%%%%%%%%%%%%%%
We first prove that
\begin{equation} %%%%%%%%%%%%%%%%%%%%%%%%%%%%%%%%
\|\nabla E\|_{L_t^2([0, T]; L_x^\infty)}
+\|\mathcal{N}\|_{L_t^2([0, T];  L_x^\infty)}\le C,
 \label{EN3}\end{equation}   %%%%%%%%%%%%%%%%% %%%%%%%%%%%%%%%%
then we show that
\begin{equation} %%%%%%%%%%%%%%%%%%%%%%%%%%%%%%%%
\|\nabla^2 E\|_{L_t^{\frac83}([0, T];  L_x^4)}+
\|\nabla \mathcal{N}\|_{L_t^{\frac83}([0, T];  L_x^4)}\le C.
 \label{EN4}\end{equation}   %%%%%%%%%%%%%%%%% %%%%%%%%%%%%%%%%
%Since
%\begin{equation} %%%%%%%%%%%%%%%%%%%%%%%%%%%%%%%%
%(i\partial_t+\Delta - \varepsilon^{2}\Delta^2)E = nE,
%\label{E1} \end{equation}  %%%%%%%%%%%%%%%%% %%%%%%%%%%%%%%%%
The Strichartz estimates \eqref{Str-3} and \eqref{EN1} give that
\begin{equation}\begin{array}{rll} %%%%%%%%%%%%%%%%%%%%%%%%%%%%%%%%
\displaystyle  \|\nabla^2 E\|_{L_t^2([0, T]; L_x^6)}
%&\lesssim&
%&\|\nabla^{2-\frac{2}{2}} E_0\|_{L_x^2}  + \|\nabla^{2-\frac{2}{2}-\frac{2}{2}}(n E)\|_{L_t^2([0,T]; %L_x^{6/5})}  \\
\displaystyle &\lesssim& \|E_0\|_{H^2}  + \|n E\|_{L_t^2([0,T]; L_x^{6/5})} \\
\displaystyle &\lesssim& \|E_0\|_{H^2} +  T^{1/2}\|n\|_{L_t^{\infty}([0,T]; H^1)} \|E\|_{L_t^{\infty}([0,T]; H^2)}\\
\displaystyle &\lesssim& T^{1/2} .
\end{array} \label{E1a} \end{equation}  %%%%%%%%%%%%%%%%%
From inequalities \eqref{Str-2}  and \eqref{EN1} we have 
\begin{equation}\begin{array}{rll} %%%%%%%%%%%%%%%%%%%%%%%%%%%%%%%%
\displaystyle  \| E\|_{L_t^2([0, T]; L_x^6)} &\lesssim&
 T^{1/4} \| E\|_{L_t^4([0, T]; L_x^6)} \\
\displaystyle &\lesssim& T^{1/4} \Big(\|E_0\|_{H^2}  + \|n E\|_{L_t^{4/3}([0,T]; L_x^{6/5})} \Big)
\displaystyle \lesssim   T.
\end{array} \label{E1b} \end{equation}  %%%%%%%%%%%%%%%%%
Hence we obtain
\begin{equation}\begin{array}{rll} %%%%%%%%%%%%%%%%%%%%%%%%%%%%%%%%
\displaystyle  \|\nabla E\|_{L_t^2([0, T]; L_x^\infty)} &\lesssim& \|\nabla E\|_{L_t^2([0,T]; H_6^1)}  \\
\displaystyle &\lesssim& \|\nabla E\|_{L_t^2([0,T];  L^6)} +\|\nabla^2 E\|_{L_t^2([0,T];  L^6)}  \\
\displaystyle &\sim& \|E\|_{L_t^2([0,T];  L^6)} +\|\nabla^2 E\|_{L_t^2([0,T];  L^6)}  \quad
 \lesssim T.
\end{array} \label{E2} \end{equation}  %%%%%%%%%%%%%%%%%
Therefore we obtain (\ref{EN3}) for $E$.

Next we estimate ${{\mathcal N}}$. Notice that ${{\mathcal N}}$ satisfies
%As for $n$, the equation is
% \begin{equation} %%%%%%%%%%%%%%%%%%%%%%%%%%%%%%%%
%(\lambda^{-2}\partial_{tt} -\Delta+ \varepsilon^{2}\Delta^2)n = \Delta |E|^2
%\label{N1} \end{equation}  %%%%%%%%%%%%%%%%% %%%%%%%%%%%%%%%%
%which can be  rewritten as follows.
\begin{equation} %%%%%%%%%%%%%%%%%%%%%%%%%%%%%%%%
i\partial_t\mathcal{N} - \lambda\sqrt{I_\varepsilon^{-1}(-\Delta)}\mathcal{N} =
\lambda\sqrt{I_\varepsilon(-\Delta)}|E|^2.
\label{N2} \end{equation}  %%%%%%%%%%%%%%%%% %%%%%%%%%%%%%%%%
%and
%\begin{equation} %%%%%%%%%%%%%%%%%%%%%%%%%%%%%%%%
%\big(i\partial_t - \lambda\sqrt{I_\varepsilon^{-1}(-\Delta)}\big)
%\big(\mathcal{N} + I_\varepsilon|E|^2\big) =
%i\partial_t I_\varepsilon |E|^2.
%\label{N3} \end{equation}  %%%%%%%%%%%%%%%%% %%%%%%%%%%%%%%%%
To apply the Strichartz estimates to \eqref{N2} %and \eqref{N3},
we need the following proposition.

\

\begin{prop}\label{Str-est-2} %%%%%%%%%%%%%%%%%%%%%%%%%%%%%%%% P5.2
Suppose that $(q, r)$ and $(\widetilde q, \widetilde r)$ are
Schr\"odinger admissible. If $u$ is a solution of
\begin{equation} %%%%%%%%%%%%%%%%%%%%%%%%
\left\{
\begin{array}{l}
i\partial_t\mathcal{N} + \sqrt{I_\varepsilon^{-1}(-\Delta)}\mathcal{N}
= h(t,x)
\qquad t\in\rre,\ x\in\rre^{3},\\
\mathcal{N}(0,x)=\mathcal{N}_{0}(x)
\qquad\qquad\qquad\qquad t\in\rre,\ x\in\rre^{3}
\end{array}
\right. \label{Sch-2}
%\end{eqnarray}
\end{equation}%%%%%%%%%%%%%%%%%
 for some
data $u_0$, $h$ and time $0<T<\infty$, then
\begin{equation} %%%%%%%%%%%%%%%%%%%%%%%%%%%%%%%%%%%%
\|\mathcal{N}\|_{L_t^q([0,T]; L_x^r)} \lesssim
\|\mathcal{N}_0\|_{L_x^2}+ \|h\|_{L_t^{\widetilde q'}([0,T]; L_x^{\widetilde r'})},
\label{Str-4} \end{equation} %%%%%%%%%%%%%%%%%%%%%%%%%%%%%%%%%%%%
\end{prop} %%%%%%%%%%%%%%%%%%%%%%%%%%%%%%%%%%%% P5.2

\vskip2mm

\noindent {\bf Proof of Proposition \ref{Str-lam}.} %%%%%%%%%%%%%%%%%%%%%%%%%%
Since $\sqrt{I_\varepsilon^{-1}(-\Delta)} \sim -\Delta$ in $L^{2}$, the proof for the
Proposition \ref{Str-est-2} is analogous to that of the
Strichartz estimate for the usual Schr\"{o}dinger equation,
see \cite[Theorem 3.1]{KPV2} for instance.
$\qquad\qed$

\begin{prop} \label{Str-lam} %%%%%%%%%%%%%%%%%%%%%%%%%%%%%%%%
Let $\mathcal{N}$ be a solution of
\begin{equation} %%%%%%%%%%%%%%%%%%%%%%%%
\left\{
\begin{array}{l}
i\partial_t\mathcal{N} + \lambda\sqrt{I_\varepsilon^{-1}(-\Delta)}\mathcal{N}
= \lambda h(t,x)
\qquad t\in\rre,\ x\in\rre^{3},\\
\mathcal{N}(0,x)=\mathcal{N}_{0}(x)
\qquad\qquad\qquad\qquad\quad t\in\rre,\ x\in\rre^{3}.
\end{array}
\right. \label{N4}
%\end{eqnarray}
\end{equation}%%%%%%%%%%%%%%%%%
Then we have
\begin{equation} %%%%%%%%%%%%%%%%%%%%%%%%%%%%%%%%
\|\mathcal{N}\|_{L_t^2H_6^s} \lesssim \lambda^{-1/2} \|\mathcal{N}_0\|_{H^s} +
\|h\|_{L_t^2H_{6/5}^s}.
\label{N5} \end{equation}  %%%%%%%%%%%%%%%%% %%%%%%%%%%%%%%%%
\end{prop} %%%%%%%%%%%%%%%%%%%%%%%%%%%%%%%%

\noindent {\bf Proof of Proposition \ref{Str-lam}.} %%%%%%%%%%%%%%%%%%%%%%%%%%
\noindent Set $\mathcal{N}(t) = \tilde{\mathcal{N}}(\lambda t)$ and
$\tilde{h}(t) = h(\lambda t)$. Then the
function $\tilde{\mathcal{N}}$ satisfies the equation (\ref{Sch-2}).
%\begin{equation} %%%%%%%%%%%%%%%%%%%%%%%%%%%%%%%%
%\big(i\partial_t + \Delta\big)v(t) = g(t).
%\label{N6} \end{equation}  %%%%%%%%%%%%%%%%% %%%%%%%%%%%%%%%%
Invoking the Strichartz estimate \eqref{Str-4}, we obtain
\begin{equation} %%%%%%%%%%%%%%%%%%%%%%%%%%%%%%%%
\|\tilde{\mathcal{N}}\|_{L_t^2H_6^s} \lesssim \|\mathcal{N}_0
\|_{H^s} + \|\tilde{h}\|_{L_t^2H_{6/5}^s}.
\label{N7} \end{equation}  %%%%%%%%%%%%%%%%% %%%%%%%%%%%%%%%%
Thus we have \eqref{N5} for $\mathcal{N}$. \qquad
$\qed$ %%%%%%%%%%%%%%%%%%%%%%%%%%

\vskip2mm

For \eqref{N2}, we invoke the Strichartz estimates (Proposition \ref{Str-lam})
%and Proposition \ref{Str-lam}
so that we can get 
\begin{equation}\begin{array}{rll} %%%%%%%%%%%%%%%%%%%%%%%%%%%%%%%%
\displaystyle  \|\mathcal{N}\|_{L_t^2H_6^2}
\displaystyle  &\lesssim& \lambda^{-1/2} \|\mathcal{N}_0\|_{H^2} +
\|\sqrt{I_\varepsilon(-\Delta)}|E|^2\|_{L_t^2H_{6/5}^2} \\
\displaystyle &\lesssim& \lambda^{-1/2} \|\mathcal{N}_0\|_{H^2} +
\|\Delta E \bar E\|_{L_t^2L_x^{6/5}} + \|\nabla E\cdot\nabla\bar E\|_{L_t^2L_x^{6/5}}\\
%\displaystyle &\lesssim& \lambda^{-1/2} \|\mathcal{N}_0\|_{H^2} +
%\|E\|_{L_t^4H^{2}}^2\\
\displaystyle &\lesssim& \lambda^{-1/2} \|\mathcal{N}_0\|_{H^2} +
T^{1/2} \|E\|_{L_t^\infty H^{2}}^2\\
\displaystyle &\lesssim& T^{1/2}.
\end{array} \label{N9} \end{equation} %%%%%%%%%%%%%%%%%
Notice that $\sqrt{I_\varepsilon(-\Delta)}$ is a bounded operator and
$\sqrt{I_\varepsilon^{-1}(-\Delta)} \sim -\Delta$. %For \eqref{N3}, we apply Proposition \ref{Str-lam} to derive
Hence by \eqref{N9}, we obtain
\begin{equation}
%\begin{array}{rll} %%%%%%%%%%%%%%%%%%%%%%%%%%%%%%%%
%\displaystyle
\|\mathcal{N} \|_{L_t^2[0,T]L_x^\infty} \leq \|\mathcal{N} \|_{L_t^2[0,T]H_6^2}
%\\&\lesssim& \|\mathcal{N} \|_{L_t^2H_6^2}^{1/2} \|\mathcal{N}\|_{L_t^2L_x^6}^{1/2}\\
\lesssim T^{1/2}.
%\end{array}
\label{N18} \end{equation}  %%%%%%%%%%%%%%%%%
Thus we have proved \eqref{EN3} for ${{\mathcal N}}$.

To obtain \eqref{EN4}, we interpolate between
\eqref{EN1} and \eqref{E1a} so that we get
\begin{equation}\begin{array}{rll} %%%%%%%%%%%%%%%%%%%%%%%%%%%%%%%%
\displaystyle \|\nabla^2 E \|_{L_t^{8/3}[0,T]L_x^4} \lesssim \|\nabla^2
E \|_{L_t^{2}[0,T]L_x^6}^{3/4} \|\nabla^2
E\|_{L_t^{\infty}[0,T]L_x^2}^{1/4} \lesssim T^{3/8}.
\end{array} \label{N19} \end{equation}  %%%%%%%%%%%%%%%%%
Also we interpolate between \eqref{EN1} and \eqref{N9} to derive
\begin{equation}\begin{array}{rll} %%%%%%%%%%%%%%%%%%%%%%%%%%%%%%%%
\displaystyle \|\nabla \mathcal{N} \|_{L_t^{8/3}[0,T]L_x^4} \lesssim
\|\nabla \mathcal{N} \|_{L_t^{2}[0,T]L_x^6}^{3/4} \|\nabla
\mathcal{N}\|_{L_t^{\infty}[0,T]L_x^2}^{1/4} \lesssim T^{3/8}.
\end{array} \label{N20} \end{equation}  %%%%%%%%%%%%%%%%%
Hence we obtain \eqref{EN4}.

By (\ref{q1}),
(\ref{q3}), (\ref{q2}) and (\ref{q4}),
we obtain
\begin{eqnarray*}
\lefteqn{\frac{d}{dt}H_{3}(E,n)(t)}\\
&\le&-4\varepsilon^{2}\int\nabla^{2}n\nabla\overline{E}\nabla^{3}Edx
-6\varepsilon^{2}\int\nabla n\nabla^{2}\overline{E}\nabla^{3}Edx
+CH_{3}(E,n)(t)\\
&\le&4\varepsilon^{2}\|\nabla E(t)\|_{L_{x}^{\infty}}H_{3}(E,n)(t)
+6\varepsilon^{2}\|\nabla n(t)\|_{L_{x}^{4}}\|\nabla^{2}E(t)\|_{L_{x}^{4}}
H_{3}^{1/2}(E,n)(t)\\
& &+CH_{3}(E,n)(t).
\end{eqnarray*}
Hence
\begin{eqnarray*}
\frac{d}{dt}H_{3}^{1/2}(E,n)(t)
&\le&(8\varepsilon^{2}\|\nabla E(t)\|_{L_{x}^{\infty}}+C)H_{3}^{1/2}(E,n)(t)\\
& &+12\varepsilon^{2}\|\nabla n(t)\|_{L_{x}^{4}}\|\nabla^{2}E(t)\|_{L_{x}^{4}}.
\end{eqnarray*}
The Gronwall lemma yields
\begin{eqnarray*}
\lefteqn{H_{3}^{1/2}(E,n)(t)}\\
&\le&H_{3}^{1/2}(E,n)(0)\exp(8\varepsilon^{2}
\int_{0}^{t}\|\nabla E(\tau)\|_{L_{x}^{\infty}}d\tau+Ct)\\
& &+12\varepsilon^{2}
\int_{0}^{t}\|\nabla n(\tau)\|_{L_{x}^{4}}\|\nabla^{2}E(\tau)\|_{L_{x}^{4}}\\
& &\qquad\qquad\times\exp(8\varepsilon^{2}
\int_{\tau}^{t}\|\nabla E(\sigma)\|_{L_{x}^{\infty}}d\sigma+C(t-\tau))
d\tau\\
&\le&H_{3}^{1/2}(E,n)(0)\exp(C
t^{1/2}\|\nabla E\|_{L_{t}^{2}L_{x}^{\infty}}+Ct)\\
& &+12\varepsilon^{2}t^{1/4}\|\nabla n\|_{L_{t}^{\frac83}L_{x}^{4}}
\|\nabla^{2}E\|_{L_{t}^{\frac83}L_{x}^{4}}
\exp(8\varepsilon^{2}
t^{1/2}\|\nabla E\|_{L_{t}^{2}L_{x}^{\infty}}+Ct).
\end{eqnarray*}
By (\ref{EN3}) and (\ref{EN4}), we have
\begin{eqnarray*}
H_{3}^{1/2}(E,n)(t)
\le H_{3}^{1/2}(E,n)(0)\exp(C
t^{3/2})+C\varepsilon^{2}t
\exp(Ct^{3/2}).
\end{eqnarray*}
Combining the above inequality and (\ref{bound}),
we have
\begin{eqnarray}
\sup_{0\le t\le T}\left(\|\nabla^{3}E\|_{L_{x}^{2}}
+\|\nabla^{2}n\|_{L_{x}^{2}}\right)\le C
\label{ee3}
\end{eqnarray}
for any $T\in(0,\infty)$ and $\lambda\in[1,\infty)$.
Further, combining (\ref{q1}),
(\ref{q3}), (\ref{q2}), (\ref{q4}) and (\ref{ee3}), we have
for $M\ge4$,
\begin{equation*}  %%%%%%%%%%%%%%%%%%%%%%%%%%%%%%%%%%%%
\frac{d}{dt}
H_{M}(E,n)(t)
\le CH_{M}(E,n)(t).
\end{equation*} %%%%%%%%%%%%%%%%%%%%%%%%%%%%%%%%%%%%
The Gronwall lemma and (\ref{bound}) yields
$H^{M}\times H^{M-1}$ bound (\ref{b1}) for $M\ge 4$.
%Now we need to find a way to show \eqref{EN5}.
This completes the
proof of Proposition
\ref{well}. $\qed$ %%%%%%%%%%%%%%%%%%%%%%%%%%%%%%%%%%%% Prop2.1

\begin{prop}\label{au} %%%%%%%%%%%%%%%%%%%%%%%%%%%%%%%%% Prop2.2
Let $d=1,2,3$ and $M\ge2$. Let $T$ be given in Proposition \ref{well}.
Then for any $(E_{0},n_{0},n_{1})\in {X_{M,\,d}}$ and $E_{0} \in H^{0,\,
M}(\mathbb{R}^d)$, the solution to (\ref{ZK}) constructed in Proposition
\ref{well} satisfies
\begin{eqnarray} %%%%%%%%%%%%%%%%%%%%%%%%%%%%%%%%%%%%
\sup_{0\le t\le T}
\sum_{\ell=0}^{M}
\sum_{k=0}^{3M-3\ell}
\||x|^{\ell}\nabla^{k}
E_{\lambda}(t)\|_{L_{x}^{2}}
&\le&C,\label{b2}
\end{eqnarray} %%%%%%%%%%%%%%%%%%%%%%%%%%%%%%%%%%%%
for any $\lambda\in[1,\infty)$.
\end{prop} %%%%%%%%%%%%%%%%%%%%%%%%%%%%%%%%%%%% Prop2.2

\vskip3mm \noindent
{\bf Proof of Proposition \ref{au}.} %%%%%%%%%%%%%%%%%%%%%%% Prop2.2
We abbreviate $E_{\lambda}$ to $E$. Let $L$ and $K$ be integers
satisfying $0\le L\le M$ and $0\le K\le 3M-3L$. We prove that there
exists a positive constant $C$ such that for any $\lambda\in[1,\infty)$,
\begin{eqnarray}
\sup_{0\le t\le T}
\||x|^{L}\nabla^{K}
E(t)\|_{L_{x}^{2}}
&\le&C,\label{yu}
\end{eqnarray}
by induction on $(L,K)$. Notice that the inequalities (\ref{yu}) with $L=0$
and $0\le K\le3M$ follow from (\ref{b1}).

We first prove that the inequalities (\ref{yu}) with $(L,K)=(0,1)$ and
$(L,K)=(0,3)$ imply the inequality (\ref{yu}) with $(L,K)=(1,0)$. We
assume the inequalities (\ref{yu}) hold for $(L,K)=(0,1)$ and
$(L,K)=(0,3)$. Taking the imaginary part of the inner product in $L^{2}$
between the first equation of (\ref{ZK}) and $|x|^{2}E$, we have
\begin{eqnarray}
\frac{d}{dt}
\||x|E(t)\|_{L_{x}^{2}}^{2}=
-2Im
\int_{\rre^{d}}
|x|^{2}\overline{E}\Delta Edx
+2
\varepsilon^{2}Im
\int_{\rre^{d}}
|x|^{2}\overline{E}
\Delta^{2} Edx.
\end{eqnarray}
Notice that
\begin{eqnarray}
Im\int_{\rre^{d}} |x|^2 \nabla\overline{E}\cdot\nabla^3E\,dx =
-Im\int_{\rre^{d}} \nabla|x|^2\cdot \nabla^3\overline{E} E\,dx.
\end{eqnarray}
By an integration by parts, we have
\begin{eqnarray}
\frac{d}{dt}
\||x|E(t)\|_{L_{x}^{2}}^{2}
&=&
2Im
\int_{\rre^{d}}
\nabla|x|^{2}\cdot\overline{E}
\nabla Edx
-4
\varepsilon^{2}Im
\int_{\rre^{d}}
\nabla|x|^{2}\cdot\overline{E}
\nabla^{3} Edx.
\nonumber\\
\label{p1}
\end{eqnarray}
The H\"{o}lder inequality yields
\begin{eqnarray*}
\frac{d}{dt}
\||x|E(t)\|_{L_{x}^{2}}^{2}
&\le&
C\|\nabla E(t)\|_{L_{x}^{2}}
\||x|E(t)\|_{L_{x}^{2}}
+C
\|\nabla^{3}E(t)\|_{L_{x}^{2}}
\||x|E(t)\|_{
L_{x}^{2}}\\
&\le&C\||x|E(t)\|_{
L_{x}^{2}}.
\end{eqnarray*}
Hence the Gronwall lemma implies
\begin{eqnarray*}
\||x|E(t)\|_{L_{x}^{2}}
\le\||x|E_{0}\|_{L_{x}^{2}}+Ct
\le C
\end{eqnarray*}
for any $t\in[0,T]$, where the constant $C$ depends on $T$ and
independent of $\lambda\in[1,\infty)$. This guarantees that the
inequality (\ref{yu}) holds for $(L,K)=(1,0)$.

Let $1\le k\le 3M-3$. Next we show that the inequalities (\ref{yu}) with
$(L,K)= (0,k+1),\ (0,k+3)$ and $(1,k')\ (k'=0,\cdots,k-1)$ imply the
inequality (\ref{yu}) with $(L,K)=(1,k)$. Assume that the inequalities
(\ref{yu}) hold for $(L,K)= (0,k+1),\ (0,k+3)$ and $(1,k')\
(k'=0,\cdots,k-1)$.

Applying the operator $\nabla^{k}$ to the first equation in (\ref{ZK}) and
taking the imaginary part of the inner product in $L^{2}$ between the
resulting equation and $|x|^{2}\nabla^{k}E$, together with invoking
integration by parts, we have
\begin{eqnarray}
& &\lefteqn{
\frac{d}{dt}
\||x|\nabla^{k}
E(t)\|_{L_{x}^{2}}^{2}}\nonumber
\\
=
& &2Im
\int_{\rre^{d}}
-|x|^{2}\nabla^{k}\overline{E}\cdot
\nabla^{k+2}E
+
\varepsilon^{2}
|x|^{2}\nabla^{k}\overline{E}\cdot
\nabla^{k+4}E\nonumber\\
& &\qquad\qquad+
|x|^{2}
\nabla^{k}\overline{E}\cdot
\nabla^{k}(nE)
dx.\\
=& &
2Im
\int_{\rre^{d}}
\nabla|x|^{2}\cdot
\nabla^{k}\overline{E}
\nabla^{k+1}Edx
-
2\varepsilon^{2}\nabla|x|^{2}\cdot
\nabla^{k}\overline{E}
\nabla^{k+3}Edx\nonumber\\
& &\qquad\qquad+
|x|^{2}
\nabla^{k}\overline{E}\cdot
\nabla^{k}
(nE)
dx.\nonumber
\label{p2}
\end{eqnarray}
The H\"{o}lder inequality and Sobolev inequality yield
\begin{eqnarray*}
\lefteqn{
\frac{d}{dt}
\||x|\nabla^{k}
E(t)\|_{L_{x}^{2}}^{2}}\\
%&\le&C
%\|\nabla^{k+1}E(t)\|_{L_{x}^{2}}
%\||x|\nabla^{k}E(t)\|_{L_{x}^{2}}
%+C
%\|\nabla^{k+3}E(t)\|_{L_{x}^{2}}
%\||x|\nabla^{k}E(t)\|_{L_{x}^{2}}\\
%& &
%+C\sum_{k'=0}^{k-1}
%\|\nabla^{k-k'}n(t)\|_{
%L_{x}^{\infty}}
%\||x|
%\nabla^{k'}E(t)\|_{
%L_{x}^{2}}
%\||x|
%\nabla^{k}E(t)\|_{
%L_{x}^{2}}
%\\
&\le&\!\!\!C\Big(
\|\nabla^{k+1}E(t)\|_{L_{x}^{2}}+
\|\nabla^{k+3}E(t)\|_{L_{x}^{2}}+
\sum_{k'=0}^{k-1}\|n(t)\|_{H_{x}^{3M-1}}
\||x|\nabla^{k'}E(t)\|_{L_{x}^{2}}\Big)\\
& &\times\||x|\nabla^{k}E(t)\|_{L_{x}^{2}}
\\
&\le&\!\!\!
C\||x|\nabla^{k}E(t)\|_{L_{x}^{2}}.
\end{eqnarray*}
Hence the Gronwall lemma implies
\begin{eqnarray*}
\||x|\nabla^{k}E(t)\|_{L_{x}^{2}}
\le C
\end{eqnarray*}
for any $t\in[0,T]$ with the constant $C$ depending on $T$ but
independent of $\lambda\in[1,\infty)$. This shows that the inequality
(\ref{yu}) holds for $(L,K)=(1,k)$.

Let $2\le\ell\le M$. Next we prove that the inequalities (\ref{yu}) with
$(L,K)=(\ell-1,1),\ (\ell-1,3)$ and $(\ell-2,2)$ imply the inequality
(\ref{yu}) with $(L,K)=(\ell,0)$. We assume the inequality (\ref{yu}) holds
for $(L,K)=(\ell-1,1),\ (\ell-1,3)$ and $(L,K)=(\ell-2,2)$. By an argument
similar to that in (\ref{p1}), we have
\begin{eqnarray*}
&& \frac{d}{dt}
\||x|^{\ell}E(t)\|_{L_{x}^{2}}^{2}\\
&=&
Im
\int_{\rre^{d}}
\Big(2\nabla|x|^{2\ell}\overline{E}
\nabla E
-4\varepsilon^{2}\nabla|x|^{2\ell}\overline{E}
\nabla^{3} E\\
& &\qquad\qquad\qquad\qquad-4\varepsilon^{2}
\ell(2\ell+d-2)
|x|^{2\ell-2}
\overline{E}\Delta E\Big)\, dx\\
&\lesssim&
\Big(\||x|^{\ell-1}\nabla E\|_{L_{x}^{2}}
+\||x|^{\ell-1}\nabla^{3}E\|_{L_{x}^{2}}
+\||x|^{\ell-2}\nabla^{2}E\|_{L_{x}^{2}}\Big)
\||x|^{\ell}E\|_{L_{x}^{2}}\\
&\lesssim&\||x|^{\ell}E(t)\|_{L_{x}^{2}},
\end{eqnarray*}
which implies the inequality (\ref{yu}) holds for $(L,K)=(\ell,0)$.

Let $2\le\ell\le M$ and $1\le k\le 3M-3\ell$. Finally we show that the
inequalities (\ref{yu}) with $(L,K)=(\ell-1,k+1),\ (\ell-1,k+3),\ (\ell-2,k+2)$
and $(\ell,k')\ (k'=0,\cdots,k-1)$ imply the inequality (\ref{yu}) with
$(L,K)=(\ell,k)$. Assume that the inequality (\ref{yu}) hold for
$(L,K)=(\ell-1,k+1),\ (\ell-1,k+3),\ (\ell-2,k+2)$ and $(\ell,k')\
(k'=0,\cdots,k-1)$.

The similar argument as given in (\ref{p2}), we have
\begin{eqnarray*}
\lefteqn{
\frac{d}{dt}
\||x|^{\ell}\nabla^{k}
E(t)\|_{L_{x}^{2}}^{2}}\\
&=&
2Im
\int_{\rre^{d}}
\nabla|x|^{2\ell}
\nabla^{k}\overline{E}
\nabla^{k+1}Edx
-
4\varepsilon^{2}Im
\int_{\rre^{d}}
\nabla|x|^{2\ell}
\nabla^{k}\overline{E}
\nabla^{k+3}Edx
\\
& &
-4\ell(2\ell+d-2)\varepsilon^{2}
Im
\int_{\rre^{d}}
|x|^{2\ell-2}
\nabla^{k}\overline{E}
\nabla^{k+2}Edx\\
& &+2Im\int_{\rre^{d}}
|x|^{2\ell}
\nabla^{k}\overline{E}
\nabla^{k}
(nE)
dx\\
&\lesssim&
\Big( \||x|^{\ell-1}\nabla^{k+1}E\|_{L_{x}^{2}}
+
\||x|^{\ell-1}\nabla^{k+3}E\|_{L_{x}^{2}}
+
\||x|^{\ell-2}\nabla^{k+2}E\|_{L_{x}^{2}} \Big)\\
& &\qquad\times
\||x|^{\ell}\nabla^{k}E\|_{L_{x}^{2}}\\
& &
+\sum_{k'=0}^{k-1}
\|\nabla^{k-k'}n\|_{
L_{x}^{\infty}}
\||x|^{\ell}
\nabla^{k'}E\|_{
L_{x}^{2}}
\||x|^{\ell}
\nabla^{k}E\|_{
L_{x}^{2}}
\\
&\lesssim&
\||x|^{\ell}\nabla^kE(t)\|_{
L_{x}^{2}}.
\end{eqnarray*}
This shows that the inequality (\ref{yu}) holds for $(L,K)=(\ell,k)$.
Collecting those estimates we obtain the inequality (\ref{b2}). $\qed$

\vskip3mm

Next we consider the fourth order nonlinear Schr\"{o}dinger type
equation (\ref{4NLS}). For the existence and the uniform bound of
solution for (\ref{4NLS}), we have the following.

\begin{prop}\label{dcm1}
Let $d=1,2,3$ and $M\ge2$. Then for any $E_{0}\in H^{M}(\rre^{d})$
there exists a unique solution to (\ref{4NLS}) satisfying
\begin{eqnarray*}
E_{\infty}\in C([0,\infty)
;H^{M}(\rre^{d})).
\end{eqnarray*}
Furthermore, for any $T\in(0,\infty)$ and some $C>0$, $E_{\infty}$
satisfies
\begin{eqnarray}
\sup_{0\le t\le T}
\|E_{\infty}(t)\|_{H_{x}^{M}}
\le C.\label{bb1}
\end{eqnarray}
\end{prop}

%\begin{prop}\label{dcm2}
%Let $d=2,3$ and $M\ge2$ be an integer.
%Then for any $E_{0}\in
%H^{M}(\rre^{d})$
%there exist $T^{\ast}\in(0,\infty)$
%and a unique solution
%to (\ref{4NLS}) satisfying
%\begin{eqnarray*}
%E_{\infty}\in C([0,T^{\ast})
%;H^{M}(\rre^{d})).
%\end{eqnarray*}
%Furthermore,
%for any $T\in(0,T^{\ast})$
%there exist positive constants $C$ such that
%$E_{\infty}$ satisfies
%(\ref{bb1}).
%\end{prop}

\noindent {\bf Proof of Proposition \ref{dcm1}.} The existence and
uniqueness of solution to (\ref{4NLS}) follows from the combination of
the Strichartz estimate for the unitary group $U_{\varepsilon}(t)$ and
the contraction mapping principle, see \cite{Caz} for instance. The global
existence of solution in $H^{2}$ follows from the conservations of mass
and Hamiltonian
\begin{eqnarray*}
\tilde{H}_{2}(E)(t)
=\frac12\|\nabla E(t)\|_{L_{x}^{2}}^{2}
+\frac{\varepsilon^{2}}{2}\|\Delta E(t)\|_{L_{x}^{2}}^{2}
-\frac14
\int_{\rre^{d}}
|E|^{2}(1-\varepsilon^{2}\Delta)^{-1}|E|^{2}dx
\end{eqnarray*}
and the Gagliardo-Nirenberg inequality (Lemma \ref{GN}). For $M\ge3$,
the global existence of solution in $H^{M}$ follows from the usual
energy method and the bound of $L_{t}^{\infty}H_{x}^{2}$ norm of $E$.
Since the proof is almost similar to that of Proposition \ref{well},
we omit the detail.$\qed$

%To derive the uniform estimate (\ref{bb1}), let
%$n_{\varepsilon,\infty}=-(1-\varepsilon^{2}\Delta)^{-1}
%|E_{\infty}|^{2}$. Then we see that
%$(E_{\infty},n_{\varepsilon,\infty})$ satisfies the following
%Sch\"{o}rodinger-elliptic system
%\begin{eqnarray*}
%\left\{
%\begin{array}{l}
%\medskip
%i\pt_{t}E_{\infty}
%+\Delta E_{\infty}
%-\varepsilon^{2}\Delta^{2}
%E_{\infty}
%=n_{\varepsilon,\infty}E_{\infty},
%\qquad t\in\rre,\ x\in\rre^{d},\\
%\medskip
%\displaystyle{
%-\Delta n_{\varepsilon,\infty}
%+\varepsilon^{2}\Delta^{2}n_{\varepsilon,\infty}
%=\Delta|E_{\infty}|^{2},
%\qquad t\in\rre,\ x\in\rre^{d}.}
%\end{array}
%\right.
%\end{eqnarray*}
%Since the proof of the inequality
%(\ref{b1}) is independent of
%$\lambda\in[0,\infty)$, we can prove (\ref{bb1})
%in a manner similar to
%the proof of (\ref{b1}).

%------------------------------------------------------------
%     Section 3.   Proof of Theorem
%------------------------------------------------------------
\section{\bf Proofs of Theorems \ref{main1}, \ref{main2}
and \ref{main3}}

In this section we prove Theorems \ref{main1}, \ref{main2} and
\ref{main3}.
%From (\ref{ZK}) and (\ref{4NLS}), we obtain
%\begin{eqnarray}\left\{\begin{array}{l}\medskip
%i\pt_{t}(E_{\lambda}-E_{\infty})
%+\Delta_\varepsilon(E_{\lambda}-E_{\infty})\\
%\medskip\qquad=
%-\{{I_{\varepsilon}}|E_{\lambda}|^{2}\}E_{\lambda}
%+\{{I_{\varepsilon}}|E_{\infty}|^{2}\}E_{\infty} +
%Q_{\lambda}E_{\lambda},\\
%(E_{\lambda}-E_{\infty})(0,x)=0,
%\end{array}\right.\label{I1}\end{eqnarray}
%where $Q_{\lambda}
%=n_{\lambda}+{I_{\varepsilon}}|E_{\lambda}|^{2}$.
The initial value problem \eqref{I} can be rewritten as the integral
equation \eqref{E-E}.
%\begin{equation}\begin{array}{ll}
%&E_{\lambda}(t)-E_{\infty}(t)\\
%=&i \displaystyle \int_{0}^{t} U_{\varepsilon}(t-s)\big[
%\{{I_{\varepsilon}}|E_{\lambda}|^{2}\}E_{\lambda}-
%\{{I_{\varepsilon}}|E_{\infty}|^{2}\}E_{\infty}-
%(Q_{\lambda}E_{\lambda})\big](s)ds,
%\end{array}\label{EE1}\end{equation}
%where $U_{\varepsilon}(t)= \exp(it\Delta_\varepsilon)$.
Hence we have
\begin{equation} %%%%%%%%%%%%%%%%%%%%%%%%%%%%%%%%%%%%%
\|E_{\lambda}(t)-
E_{\infty}(t)\|_{H_{x}^{m}}
\le J_1+J_2, \label{i1}
\end{equation} %%%%%%%%%%%%%%%%%%%%%%%%%%%%%%%%%%%%%%%
where
\begin{equation} %%%%%%%%%%%%%%%%%%%%%%%%%%%%%%%%%%%%%%%
\displaystyle J_1=\int_{0}^{t}
\|[\{{I_{\varepsilon}}|E_{\lambda}|^{2}\}E_{\lambda}-
\{{I_{\varepsilon}}|E_{\infty}|^{2}\}E_{\infty}](s)\|_{H_{x}^{m}}ds
\end{equation} %%%%%%%%%%%%%%%%%%%%%%%%%%%%%%%%%%%%%%%
and
\begin{equation} %%%%%%%%%%%%%%%%%%%%%%%%%%%%%%%%%%%%%%%
\displaystyle J_2=\int_{0}^{t}
\|(Q_{\lambda}E_{\lambda})(s)
\|_{H_{x}^{m}}ds.
\end{equation} %%%%%%%%%%%%%%%%%%%%%%%%%%%%%%%%%%%%%%%
We first evaluate $J_{1}$. Since the operator ${I_{\varepsilon}}$ is
bounded from $H_{x}^{m}$ to $H_{x}^{m-2}$, we obtain 
\begin{eqnarray*} &
&\|\{{I_{\varepsilon}}|E_{\lambda}|^{2}\}E_{\lambda}-
%\{{I_{\varepsilon}}|E_{\infty}|^{2}\}E_{\lambda}
%\\& &+
%\{{I_{\varepsilon}}|E_{\infty}|^{2}\}E_{\lambda}-
\{{I_{\varepsilon}}|E_{\infty}|^{2}\}E_{
\infty}\|_{H_{x}^{m}}\\
&\le&
\|{I_{\varepsilon}}
\{|E_{\lambda}|^{2}-|E_{\infty}|^{2}\}
\|_{H_{x}^{m}}
\|E_{\lambda}\|_{H_{x}^{m}}
+
\|{I_{\varepsilon}}|E_{\infty}|^{2}
\|_{H_{x}^{m}}
\|E_{\lambda}-E_{
\infty}\|_{H_{x}^{m}}\\
%&\le&C
%\||E_{\lambda}|^{2}-|E_{\infty}|^{2}
%\|_{H_{x}^{m-2}}
%\|E_{\lambda}\|_{H_{x}^{m}}
%+C
%\||E_{\infty}|^{2}
%\|_{H_{x}^{m-2}}
%\|E_{\lambda}-E_{
%\infty}\|_{H_{x}^{m}}\\
&\le&
C(\|E_{\lambda}\|_{H_{x}^{m}}^{2}
+\|E_{\infty}\|_{H_{x}^{m-2}}^{2})
\|E_{\lambda}-E_{
\infty}\|_{H_{x}^{m}}.
%\qquad m>2+d/2.
\end{eqnarray*} 
Combining the above inequalities and (\ref{b2}), we have
\begin{equation}
J_{1}\le C\int_{0}^{t}
\|E_{\lambda}(s)-E_{
\infty}(s)\|_{H_{x}^{m}}ds,
\label{i2}
\end{equation}
where the constant $C$ is independent of $\lambda\in[1,\infty)$.

To evaluate $J_{2}$, we need to estimate $Q_{\lambda}$ given in
\eqref{INT} and thus rewrite it as
\begin{eqnarray}
Q_{\lambda}(t)
&=&
\cos(\lambda t\omega_\varepsilon) f_0
+\frac{\sin(\lambda t\omega_\varepsilon)}{\lambda
\omega_\varepsilon}\nabla\cdot f_1
%\nonumber\\& &
+
\int_{0}^{t}
\frac{\sin(\lambda (t-s)\omega_\varepsilon)}{\lambda
\omega_\varepsilon}
\pt_{t}^{2}{I_{\varepsilon}}|E_{\lambda}|^{2}
(s)ds \nonumber\\
&\equiv&
Q_{\lambda}^{(0)}(t)
+Q_{\lambda}^{(1)}(t)
+Q_{\lambda}^{(2)}(t),
\label{i3}
\end{eqnarray}
where
\begin{equation}%%%%%%%%%%%%%%%%%%%%%%%%%%%%%%%%%
f_0=n_{0}+{I_{\varepsilon}}|E_{0}|^{2},\quad \nabla\cdot f_1=
n_{1}+2Im\{E_{0}\Delta_\varepsilon\overline{E}_{0}\},
\quad\text{and}\quad \omega_\varepsilon=\sqrt{-\Delta_\varepsilon}.
\label{Q0-1}
\end{equation} %%%%%%%%%%%%%%%%%%%%%%%%%%%%%%%%%
Notice that $f_1= \phi+2Im\{E_0\nabla
{I_\varepsilon^{-1}}\overline{E}_0+ \varepsilon^2\nabla
E_0\Delta\overline{E}_0\}$, where $\nabla\cdot\phi=n_1$.

 For the three terms
$Q_{\lambda}^{(j)},\ j=0,1,2$, we have the following lemmas.

\begin{lem}\label{1f} %%%%%%%%%%%%%%%%%%%%%%%%%%%%%%%%%
Let $d=1, 2, 3,$ and $m\ge2$.  For any $0\le k\le m$, the inequality
\begin{equation} \begin{array}{lll} %%%%%%%%%%%%%%%%%%%%%%%%%%%%%%%%
&&|\pt_{x}^{k}
Q_{\lambda}^{(0)}(t,x_0)|\\
&\lesssim&\left\{
\begin{array}{l}
\displaystyle{
(1+\lambda t)^{-2}(1+|x_0|)^{2}
\|f_0\|_{H_{x}^{m+ 1+[d/2]} }
\quad
if\ |x_0|\ge\lambda t/2\ or\ 0\le\lambda t\le1,}\\
\displaystyle{ (1+\lambda t)^{ -\frac{d}{2}-\sigma } \Big(
(d-1)\|f_0\|_{\dot H_{x}^{-\sigma}}+ \sum_{j=0}^{2} \||x|^{j}f_0\|_{
H_{x}^{m+j-3+[d/2]}} \Big) }\\
\qquad\qquad\qquad\qquad\qquad\qquad if\ |x_0|\le\lambda t/2,\, \lambda t>1,
\end{array}
\right.
\end{array}\label{L3.1} \end{equation} %%%%%%%%%%%%%%%%%%%%%%%%%%%%%
holds for any $t\in(0,\infty)$ and $\lambda\in[1,\infty)$, where
 $\sigma=3/2$ for $d=1$,
$0\leq\sigma<1$ for $d=2$, and $0\leq\sigma<\frac12$ for $d=3$.
\end{lem} %%%%%%%%%%%%%%%%%%%%%%%%%%%%%%%%%

\begin{lem}\label{2ff} %%%%%%%%%%%%%%%%%%%%%%%%%%%%%%%%%
Let $d=1,2,3$ and $m\ge2$. Then the inequality
\begin{equation} %%%%%%%%%%%%%%%%%%%%%%%%%%%%%%%%%
\|Q_{\lambda}^{(1)}(t)
\|_{H_{x}^{m}}
\le C \lambda^{-1}
(\|n_{1}\|_{H_{x}^{m-1}}
+\|n_{1}\|_{\dot{H}_{x}^{-1}}
+\|E_{0}\|_{H_{x}^{m+3}}^{2}) \label{Q1-1}
\end{equation} %%%%%%%%%%%%%%%%%%%%%%%%%%%%%%%%%
holds for any $t\in(0,\infty)$ and $\lambda\in[1,\infty)$.
\end{lem} %%%%%%%%%%%%%%%%%%%%%%%%%%%%%%%%%

\begin{lem}\label{2f} %%%%%%%%%%%%%%%%%%%%%%%%%%%%%%%%%
Let $d=1, 2, 3,$ and $m\ge2$.  For any $0\le k\le m$, the inequality
\begin{equation} \begin{array}{ll} %%%%%%%%%%%%%%%%%%%%%%%%%%%%%%%%
&|\pt_{x}^{k}
Q_{\lambda}^{(1)}(t,x_0)|\\
\lesssim&\!\!\left\{\!\!
\begin{array}{l}
\displaystyle
\lambda^{-1}(1+\lambda t)^{-2}(1+|x_0|)^{2} \|f_1\|_{H_{x}^{m+
2+[d/2]} \quad\displaystyle if\ |x_0|\ge\lambda t/2\ or\ 0\le\lambda t\le1,}\\
\displaystyle \lambda^{-1} (1+\lambda t)^{ -\frac{d}{2}-\sigma } \Big(
(d-1)\|f_1\|_{\dot H_{x}^{-\sigma}}+ \sum_{j=0}^{2} \||x|^{j}f_1\|_{
H_{x}^{m+j-4+[d/2]}} \Big)\\
\qquad\qquad\qquad\qquad\qquad\qquad
if\ |x_0|\le\lambda t/2,\, \lambda t>1,
\end{array}
\right.
\end{array}\label{L3.3} \end{equation} %%%%%%%%%%%%%%%%%%%%%%%%%%%%%
holds for any $t\in(0,\infty)$ and $\lambda\in[1,\infty)$, where
 $\sigma=3/2$ for $d=1$,
$0\leq\sigma<1$ for $d=2$, and $0\leq\sigma<\frac12$ for $d=3$.
\end{lem} %%%%%%%%%%%%%%%%%%%%%%%%%%%%%%%%%

\begin{lem}\label{3ff} %%%%%%%%%%%%%%%%%%%%%%%%%%%%%%%%%
Let $d=1,2,3$ and $m\ge3$. Let $T^{\ast}$ be given by Proposition
\ref{well}. Then for any $T\in(0,\infty)$ for $d=1, 2$ and for any
$T\in(0,T^{\ast})$ for $d=3$ and for any $(E_{0},n_{0},n_{1}) \in
H^{m+4}(\rre^{d})\times H^{m+3}(\rre^{d})\times H^{m+1}(\rre^{d})$,
the inequality
\begin{equation} %%%%%%%%%%%%%%%%%%%%%%%%%%%%%%%%%
\|Q_{\lambda}^{(2)}(t)
\|_{H_{x}^{m}}
\le C \lambda^{-1}
T\sup_{0\le t\le T}
(1+\|n_{\lambda}(t)\|_{H_{x}^{m}})
\|E_{\lambda}(t)\|_{H_{x}^{m+4}}^{2} \label{Q2-1}
\end{equation} %%%%%%%%%%%%%%%%%%%%%%%%%%%%%%%%%
holds for any $\lambda\in[1,\infty)$ where $(E_{\lambda},n_{\lambda})$
is the solution to (\ref{ZK}) given in Proposition \ref{well}.
\end{lem} %%%%%%%%%%%%%%%%%%%%%%%%%%%%%%%%%

\begin{lem}\label{3f} %%%%%%%%%%%%%%%%%%%%%%%%%%%%%%%%%
Let $d=1, 2, 3$ and $T^{\ast}$ be given by Proposition \ref{well}. Let
$m\ge6-[d/2]$. and $M\ge4$ be integers satisfying $3M\ge m+6$. Then
for any $T\in(0,\infty)$ for $d=1, 2$ and for some $T\in(0,T^{\ast})$ for
$d=3$, $(E_{0},n_{0},n_{1}) \in X_{3M,\,d}$, and $E_{0}\in
H^{0,M}(\rre^{d})$, the inequality
\begin{eqnarray} %%%%%%%%%%%%%%%%%%%%%%%%%%%%%%%%%
\lefteqn{|\pt_{x}^{k}
Q_{\lambda}^{(2)}(t,x_0)|}\nonumber\\
&\le&C\left\{
\begin{array}{l}
\displaystyle{
\lambda^{-2}(1+|x_0|)
\sup_{0\le t\le T}\left\{
\big(1+\|n_{\lambda}(t)\|_{H_{x}^{m+4+[\frac{d}{2}]}}\big)
\|E_{\lambda}(t)\|_{H_{x}^{m+5+[\frac{d}{2}]}}^{2}\right\}}
\\
\qquad\quad\qquad\qquad\quad\qquad\qquad\quad\qquad
if\ |x_0|\ge\lambda t/2\ or\ 0\le\lambda t\le1,\\
\displaystyle{
F_{d}(\lambda)
\sup_{0\le t\le T}\left\{
\big(1+\|n_{\lambda}(t)\|_{H_{x}^{m+2+[\frac{d}{2}]}}\big)
\left(\sum_{\ell=0}^{2}
\||x|^{\ell}E_{\lambda}(t)\|_{H_{x}^{m+3+[\frac{d}{2}]-3\ell}}\right)^{2}\right\}}
\\
\displaystyle{
+\lambda^{-2}(1+|x_0|)
\sup_{0\le t\le T}\left\{
\big(1+\|n_{\lambda}(t)\|_{H_{x}^{m+4+[\frac{d}{2}]}}\big)
\|E_{\lambda}(t)\|_{H_{x}^{m+5+[\frac{d}{2}]}}^{2}\right\}}
\\
\qquad\quad\qquad\qquad\quad\qquad\qquad\quad\qquad
if\ |x_0|\le\lambda t/2,\, \lambda t>1
\end{array}
\right.\nonumber\\
\label{y7}
\end{eqnarray} %%%%%%%%%%%%%%%%%%%%%%%%%%%%%%%%%
holds for any $t\in(0,T)$, $0\le k\le m$ and $\lambda\in[1,\infty)$,
where $F_{d}(\lambda)=\lambda^{-2}$ for $d=1$,
$\lambda^{-2}\log\lambda$ for $d=2$, and $\lambda^{-2}$ for $d=3$.
\end{lem} %%%%%%%%%%%%%%%%%%%%%%%%%%%%%%%%%

We shall prove Lemmas \ref{1f}, \ref{2ff}, \ref{2f}, \ref{3ff} and \ref{3f} in
the next section. Now we prove Theorems \ref{main1} and \ref{main2}
assuming that Lemmas \ref{1f}, \ref{2ff}, \ref{2f}, \ref{3ff} and \ref{3f}
hold.

\vskip3mm \noindent {\bf Proof of Theorems \ref{main1}, \ref{main2}
and \ref{main3}}. We first consider the case when
$n_{0}+{I_{\varepsilon}}|E_{0}|^{2}\neq0$. The inequality (\ref{T1-2})
follows from Lemmas \ref{2ff} and \ref{3ff}:
\begin{eqnarray*}%%%%%%%%%%%%%%%%%%%%%%%%%%%%
\|Q_{\lambda}(t)
-Q_{\lambda}^{(0)}(t)\|_{H_{x}^{m}}
\le
\|Q_{\lambda}^{(1)}(t)\|_{H_{x}^{m}}
+\|Q_{\lambda}^{(2)}(t)\|_{H_{x}^{m}}
\le C
\lambda^{-1}.
\end{eqnarray*}%%%%%%%%%%%%%%%%%%%%%%%%%%%%
Let us show the inequality (\ref{T1-1}). From Lemma \ref{1f},
Propositions \ref{well} and \ref{au} we have
\begin{eqnarray}%%%%%%%%%%%%%%%%%%%%%%%%%%%%
\lefteqn{
\|Q_{\lambda}^{(0)}
E_{\lambda}(t)\|_{H_{x}^{m}}}
\nonumber\\
&\le&C
(1+\lambda t)^{-2}
\|(1+|x|)^{2}
E_{\lambda}(t)\|_{H_{x}^{m}}
+C(1+\lambda t)^{-\mu}
\|E_{\lambda}(t)\|_{H_{x}^{m}}
\nonumber\\
&\le&C(1+\lambda t)^{-\mu},
\label{i4}
\end{eqnarray}%%%%%%%%%%%%%%%%%%%%%%%%%%%%
where $\mu = 2$  for $d=1$ and $\mu=d/2+\sigma$ for $d=2,3$, and the
constant $C$ depends on $E_{0}$ and $n_{0}$. We deduce from
Lemmas \ref{2ff} and \ref{3ff} that
\begin{eqnarray} %%%%%%%%%%%%%%%%%%%%%%%%%%%%%%%%%
\|(Q_{\lambda}^{(1)}
+Q_{\lambda}^{(2)})
E_{\lambda}(t)\|_{H_{x}^{m}}
\le C
\lambda^{-1}.
\label{i5}
\end{eqnarray} %%%%%%%%%%%%%%%%%%%%%%%%%%%%%%%%%
Therefore (\ref{i4}) and (\ref{i5}) together yield
\begin{eqnarray} %%%%%%%%%%%%%%%%%%%%%%%%%%%%%%%%%
J_{2}\le C
\int_{0}^{t}(1+\lambda s)^{-\mu}ds
+C\lambda^{-1}\int_{0}^{t}ds
\le C\lambda^{-1},
\label{i6}
\end{eqnarray} %%%%%%%%%%%%%%%%%%%%%%%%%%%%%%%%%
where the constant $C$ depends on $T$ but is independent of
$\lambda\in[1,\infty)$. Combining (\ref{i1}), (\ref{i2}) and (\ref{i6}), we
see
\begin{eqnarray*} %%%%%%%%%%%%%%%%%%%%%%%%%%%%%%%%%
\|E_{\lambda}(t)
-E_{\infty}(t)\|_{H_{x}^{m}}
\le C\int_{0}^{t}
\|E_{\lambda}(s)
-E_{\infty}(s)\|_{H_{x}^{m}}ds
+C\lambda^{-1}.
\end{eqnarray*} %%%%%%%%%%%%%%%%%%%%%%%%%%%%%%%%%
The Gronwall lemma implies
\begin{eqnarray*} %%%%%%%%%%%%%%%%%%%%%%%%%%%%%%%%%
\|E_{\lambda}(t)
-E_{\infty}(t)\|_{H_{x}^{m}}
\le C\lambda^{-1}\exp(CT).
\end{eqnarray*} %%%%%%%%%%%%%%%%%%%%%%%%%%%%%%%%%
Hence we have (\ref{T1-1}) and (\ref{T1-2}). For the case where
$n_{0}+{I_{\varepsilon}}|E_{0}|^{2}\equiv0$, the inequalities (\ref{T1-3}),
(\ref{T2-3}) and (\ref{T3-3}) follow from the argument similar as above.
Indeed it suffices to replace Lemmas \ref{2ff} and \ref{3ff} by Lemmas
\ref{2f} and \ref{3f}, respectively. This completes the proof of Theorems
\ref{main1} , \ref{main2} and \ref{main3}. \qquad$\qed$

%------------------------------------------------------------
%     Section 4.   Proof of Lemma
%------------------------------------------------------------
\section{\bf Estimates for $Q_{\lambda}$}

In this section we prove Lemmas \ref{1f} - \ref{3f}. We first denote some
notations which will be used throughout this section.
%Set the interval $I_0={\rre}$ for $d=1$ and $I_0=[0,\,\infty)$ for $d=2,\, 3$.
Then we set the functions
\begin{equation}  %%%%%%%%%%%%%%%%%%%%%%%%%%%%%%%%%%%%%%%%%%%%%%
\varphi_{\pm}(t, \xi) =x_0\cdot\xi\pm\lambda
t \xi_\varepsilon ,   \label{p+-}
\end{equation} %%%%%%%%%%%%%%%%%%%%%%%%%%%%%%%%%%%%%%%%%%%%%%
where
\begin{equation}  %%%%%%%%%%%%%%%%%%%%%%%%%%%%%%%%%%%%%%%%%
\xi_\varepsilon=\left\{\begin{array}{ll}
\xi(1+\varepsilon^{2}|\xi|^{2})^{1/2} &\text{for}\quad d=1,\\
|\xi|(1+\varepsilon^{2}|\xi|^{2})^{1/2} &\text{for}\quad d=2, 3.
\label{xi-ep}
\end{array}\right. \end{equation}  %%%%%%%%%%%%%%%%%%%%%%%%%%%%%%%%%

To estimate the integrals $Q_{\lambda}^{(k)}(t, x_0)$ for $k=0, 1, 2$, we
first compute the partial derivatives of $\varphi_{\pm}$ over the
variable $\xi$ and we get, for $d=1$,
\begin{eqnarray} %%%%%%%%%%%%%%%%%%%%%%%%%%%%%%%%
& &\varphi_{\pm}' =x_0\pm\lambda t
\frac{1+2\varepsilon^{2}\xi^{2}}{(1+\varepsilon^{2}\xi^{2})^{1/2}},\,\,
\qquad\varphi_{\pm}'' =\pm\lambda t \frac{\varepsilon^{2}\xi
(3+2\varepsilon^{2}\xi^{2})}{(1+\varepsilon^{2}\xi^{2})^{3/2}}
\nonumber\\
& &\varphi_{\pm}'''=\pm\lambda t
\frac{3\varepsilon^{2}}{(1+\varepsilon^{2}\xi^{2})^{5/2}}.
\label{p123=}
\end{eqnarray}  %%%%%%%%%%%%%%%%%%%%%%%%%%%
We note the oscillatory integrals $Q_{\lambda}^{(k)}(t, x_0)$ for $k=0, 1,
2$ have no stationary point in the region $|x_0|<\lambda t$.

\vskip1mm

For $d=1$, combining the identity
\begin{equation} %%%%%%%%%%%%%%%%%%%%%%%%%%%
e^{i\varphi_{\pm}(t, \xi)}=
\frac{\pt_{\xi} e^{i\varphi_{\pm}(t, \xi)}}
{i\varphi_{\pm}'(t, \xi)} \label{e^ip1}
\end{equation} %%%%%%%%%%%%%%%%%%%%%%%%%%%
and repeating the integration by parts twice, we have
\begin{equation} \begin{array}{rll}   %%%%%%%%%%%%%%%%%%%%%%%%%%%%%%%%
&\displaystyle \int_{-\infty}^\infty e^{i\varphi_{\pm}(t, \xi)}\widehat{g}(\xi)d\xi \\
%= &\displaystyle \int_{-\infty}^\infty e^{i\varphi_{\pm}(t, \xi)}
%\pt_{\xi}
%\left\{ \frac{1}{i\varphi_{\pm}'(t, \xi)} \pt_{\xi}
%\left(\frac{\widehat{g}(\xi)}{ i\varphi_{\pm}'(t, \xi)}
%\right)\right\}d\xi \\
= &\displaystyle \int_{-\infty}^\infty e^{i\varphi_{\pm}(t, \xi)}
\Big[\frac{-1}{(\varphi_{\pm}^{'})^2}
\pt_{\xi}^{2}\widehat{g} \,+\,
\frac{3\varphi_{\pm}^{''}}{(\varphi_{\pm}^{'})^3}\pt_{\xi}\widehat{g}\,+\,
\Big(\frac{\varphi_{\pm}^{'''}}{(\varphi_{\pm}^{'})^3}-\frac{3(\varphi_{\pm}^{''})^2}{(\varphi_{\pm}^{'})^4}\Big)
\widehat{g}(\xi)\Big]d\xi.
\end{array} \label{Q0jd1}  \end{equation} 
%%%%%%%%%%%%%%%%%%%%%%%%%%%%

\vskip1mm

 For $d=2, 3$, the phase function is $\varphi_\pm =
x_0\cdot\xi \pm |\xi|\sqrt{1+\varepsilon^2|\xi|^2}$ whose second
derivative $\Delta\phi_{\pm} \sim \pm\lambda t
\big((d-1)|\xi|^{-1}\langle\xi\rangle+ |\xi|\langle\xi\rangle^{-1}\big)$ is
singular at the origin. This would require that the initial data lies in $\dot
H^{-2}$ which is more than that is required in the Hamiltonian
\eqref{Hami}.  To avoid the difficulty, we modify the identity
\eqref{e^ip1} as follows:
\begin{equation} %%%%%%%%%%%%%%%%%%%%%%%%%%%
e^{i\varphi_{\pm}(t, \eta)}=
\frac{\pt_{\eta}(\eta e^{i\varphi_{\pm}(t, \eta)})}
{1+i\eta\varphi_{\pm}'(t, \eta)}. \label{e^ip}
\end{equation} %%%%%%%%%%%%%%%%%%%%%%%%%%%
Combining the identity (\ref{e^ip}) and repeating the integration by
parts twice, we have 
\begin{equation} \begin{array}{rll}   %%%%%%%%%%%%%%%%%%%%%%%%%%%%%%%%
&\displaystyle \int_0^\infty e^{i\varphi_{\pm}(t, \eta)}\widehat{g_k}(\eta)d\eta \\
%= &\displaystyle \int_0^\infty e^{i\varphi_{\pm}(t, \eta)}
%\eta\pt_{\eta}
%\left\{ \frac{\eta}{1+i\eta\varphi_{\pm}'(t, \eta)} \pt_{\eta}
%\left(\frac{\widehat{g_k}(\eta)}{ 1+i\eta\varphi_{\pm}'(t, \eta)}
%\right)\right\}d\eta \\
= &\displaystyle \int_0^\infty
e^{i\varphi_{\pm}}q_{1}\pt_{\eta}^{2}\widehat{g_k}(\eta)d\eta \,\,+\,\, \int_0^\infty
e^{i\varphi_{\pm}}q_{2}\pt_{\eta}\widehat{g_k}(\eta)d\eta  +\,\, \int_0^\infty
e^{i\varphi_{\pm}}q_{3}\pt_{\eta}\widehat{g_k}(\eta)d\eta\,\,\\
&\displaystyle +\,\, \int_0^\infty
e^{i\varphi_{\pm}}q_{4}\widehat{g_k}(\eta)d\eta
\displaystyle +\,\, \int_0^\infty
e^{i\varphi_{\pm}}q_{5}\widehat{g_k}(\eta)d\eta \\
\displaystyle \equiv &Q_{k,1}(t,
x_0)+Q_{k,2}(t,x_0)+Q_{k,3}(t,x_0)+Q_{k,4}(t,x_0)+Q_{k,5}(t,x_0),
\end{array} \label{Q0j}  \end{equation} %%%%%%%%%%%%%%%%%%%%%%%%%%%%
where the functions $g_k$ for $k=0,1, 2, $ will be specified later,
\begin{equation}\begin{array}{rll}%%%%%%%%%%%%%%%%%%%%%%%%%%%%%%%%
&\displaystyle q_{1}= \frac{\eta^{2}}{
(1+i\eta\varphi_{\pm}'(t, \eta))^{2}}, \\
&\displaystyle q_{2}= \frac{\eta}{(1+i\eta\varphi_{\pm}'(t, \eta))^{2}},\\
&\displaystyle q_{3}= \frac{3\eta^{2}(\varphi_{\pm}'(t, \eta)+\eta\varphi_{\pm}''(t, \eta))}{
-i(1+i\eta\varphi_{\pm}'(t, \eta))^{3}}, \\
&\displaystyle q_{4}= \frac{\eta(\varphi_{\pm}'(t, \eta)+3\eta\varphi_{\pm}''(t, \eta)+\eta^{2}
\varphi_{\pm}'''(t, \eta))}{-i(1+i\eta\varphi_{\pm}'(t, \eta))^{3}}, \\
&\displaystyle q_{5}= \frac{3\eta^{2}(\varphi_{\pm}'(t, \eta)+\eta\varphi_{\pm}''(t, \eta))^{2}}{
-(1+i\eta\varphi_{\pm}'(t, \eta))^{4}}.
\end{array} \label{qj=}\end{equation}   %%%%%%%%%%%%%%%%%%%%%%%%%%%%
For any $\eta\in [0, \infty)$ and in the region $|x_0|\le(\lambda t)/2$,
we have
\begin{equation}\begin{array}{rll}%%%%%%%%%%%%%%%%%%%%%%%%%%%%%%%%
&\displaystyle |\varphi_{\pm}'(t, \eta)|\sim \lambda t(1+\varepsilon|\eta|), %\\
&\displaystyle |\varphi_{\pm}^{\prime\prime}(t, \eta)|\sim \lambda t
\frac{\varepsilon^2|\eta|}{1+\varepsilon|\eta|}, \\
&\displaystyle |\varphi_{\pm}^{\prime\prime\prime}(t, \eta)|\sim \lambda t
\frac{\varepsilon^2}{(1+\varepsilon|\eta|)^5}, %\\
&\displaystyle |1+i\eta\varphi_{\pm}'(t, \eta)|\sim 1+|\eta|\lambda t(1+\varepsilon|\eta|).
\end{array} \label{p123} \end{equation}   %%%%%%%%%%%%%%%%%%%%%%%%%%
Invoking \eqref{qj=}
and \eqref{p123}, we can obtain
\begin{equation}\begin{array}{rll}%%%%%%%%%%%%%%%%%%%%%%%%%%%%%%%%
&\displaystyle |q_{1}|\sim \frac{\eta^{2}}{1+(|\eta|\lambda t)^2(1+\varepsilon|\eta|)^{2}}, %\\
&\displaystyle |q_{2}|\sim \frac{|\eta|}{1+(|\eta|\lambda t)^2(1+\varepsilon|\eta|)^{2}},\\
&\displaystyle |q_{3}|\sim \frac{\eta^{2}\lambda t(1+\varepsilon|\eta|)}
{1+(|\eta|\lambda t)^3(1+\varepsilon|\eta|)^{3}}, %\\
&\displaystyle |q_{4}|\sim \frac{|\eta|\lambda t(1+\varepsilon|\eta|)}
{1+(|\eta|\lambda t)^3(1+\varepsilon|\eta|)^{3}}, \\
&\displaystyle |q_{5}|\sim \frac{(|\eta|\lambda
t)^2(1+\varepsilon|\eta|)^2}{1+(|\eta|\lambda
t)^4(1+\varepsilon|\eta|)^{4}}.
\end{array} \label{qj}\end{equation}   %%%%%%%%%%%%%%%%%%%%%%%%%%%%
To estimate the quantities $Q_{k,j}$ for $k=0, 1, 2 $ and $j=1,
\cdot\cdot\cdot, 5$, we split each of the integrals into three parts which
are on the intervals $I_1$, $I_2$, and $I_3$ given by
\begin{equation}%%%%%%%%%%%%%%%%%%%%%%%%%%%%%%%%%%%%%%%%%%%%%%
I_1(t)=\{|\xi|\leq(\lambda t)^{-1}\},\quad %%%%%%%%%%%%%%
I_2(t)=\{(\lambda t)^{-1}<|\xi|\leq1\},\quad \text{and}\quad %%%
I_3(t)=\{|\xi|>1\}. %%%%%%%%%%%%%%%%%%%%%%%%%%%%
\label{I123}
\end{equation}  %%%%%%%%%%%%%%%%%%%%%%%%%%%%%%%%%%%%%%%%%%%%%%
 Observe that we have
\begin{equation}%%%%%%%%%%%%%%%%%%%%%%%%%%%%%%%%%%%%%%%%%%%%%%
|\eta||q_{5}|\lesssim|\eta||q_{4}|\sim|q_{3}|\lesssim|q_{2}|=|\eta|^{-1}|q_{1}|\lesssim
\left\{\begin{array}{rll}%%%%%%%%%%%%%%%%%%%%%%%%%%%%%%%%%%%%%%%%%%%
&\displaystyle |\eta| & \text{for}\,\,\eta\in I_1(t),  \\
&\displaystyle (\lambda t)^{-2}|\eta|^{-1} & \text{for}\,\,\eta\in I_2(t), \\
&\displaystyle (\lambda t)^{-2}|\eta|^{-3} & \text{for}\,\,\eta\in I_3(t).
\end{array}\right.
\label{qj<}
\end{equation}  %%%%%%%%%%%%%%%%%%%%%%%%%%%%%%%%%%%%%%%%%%%%%

\vskip3mm
\noindent
{\bf Proof of Lemma \ref{1f}} %%%%%%%%%%%%%%%%%%%%%%%%%%%%%%%%%%%%
We use the representation
\begin{equation}%%%%%%%%%%%%%%%%%%%%%%%%%%%%%%%%%%%%%%%%%%%%%%
Q_{\lambda}^{(0)}(t, x_0)
=
\left(\frac{1}{2\pi}\right)^{d/2}\frac12
Re\left[
\int_{\rre^{d}}
e^{ix_0\cdot\xi}\big(e^{i\lambda t\xi_\varepsilon}
+ e^{-i\lambda t\xi_\varepsilon}\big)
 \widehat{f_0}(\xi)d\xi\right],\label{Q(0)}
\end{equation}%%%%%%%%%%%%%%%%%%%%%%%%%%%%%%%%%%%%%%%%%%%%%%
where
\begin{equation}%%%%%%%%%%%%%%%%%%%%%%%%%%%%%%%%%%%%%%%%%%%%%%
\xi_\varepsilon=\xi\sqrt{1+\varepsilon^2\xi^2}\quad\text{for}\ d=1,\quad
\xi_\varepsilon=|\xi|\sqrt{1+\varepsilon^2|\xi|^2}\quad\text{for}\ d=2, 3. \label{xi_ep}
\end{equation}%%%%%%%%%%%%%%%%%%%%%%%%%%%%%%%%%%%%%%%%%%%%%%%%
For the case $0\le\lambda t\le1$, for any $0\le k\le m$, we have
\begin{eqnarray*}%%%%%%%%%%%%%%%%%%%%%
|\nabla^{k}Q_{\lambda}^{(0)}(t,x_0)|
&\le &C
\int_{\rre^{d}}
\langle\xi\rangle^{m}
|\widehat{f_0}(\xi)|d\xi\\
&\le &C
\|f_0\|_{H^{m+1+[d/2]}}.
%\left\{\begin{array}{rl}
%&\displaystyle \|f_0\|_{H^{m+1}}\quad\text{for}\quad d=1, \\
%&\displaystyle \|f_0\|_{H^{m+2}}\quad\text{for}\quad d=2, 3.
%\end{array}\right.
\end{eqnarray*}%%%%%%%%%%%%%%%%%%%%%%
Next we consider the case $\lambda t>1$ and $|x_0|\ge \lambda t/2$. 
Then we have
\begin{eqnarray*} %%%%%%%%%%%%%%%%%%%%%%%%%%%%%%%%%
|\nabla^{k}
Q_{\lambda}^{(0)}(t,x_0)|
&\le&C
\int_{\rre^{d}}
\langle\xi\rangle^{m}
|\widehat{f_0}(\xi)|d\xi\\
&\le&C
(\lambda t)^{-2}|x_0|^{2}\|f_0\|_{H^{m+1+[d/2]}}.
%\left\{\begin{array}{rl}
%&\displaystyle (\lambda t)^{-2}
%|x_0|^{2}\|f_0\|_{H^{m+1}}\quad\text{for}\quad d=1, \\
%&\displaystyle (\lambda t)^{-2}
%|x_0|^{2}\|f_0\|_{H^{m+2}}\quad\text{for}\quad d=2, 3.
%\end{array}\right.
\end{eqnarray*} %%%%%%%%%%%%%%%%%%%%%%%%%%%%%%%%%

Let us evaluate $Q_{\lambda}^{(0)}(t, x_0)$ in the region where
$|x_0|\le \lambda t/2$. We rewrite \eqref{Q(0)} as
%\begin{eqnarray*} %%%%%%%%%%%%%%%%%%%%%%%%%%%%%%%%%%%%
%\left(\frac{1}{2\pi}\right)^{1/2}\frac12
%Re\int_0^\infty
%\Big(e^{i\varphi_{+}(t, \eta)}+
%e^{-i\varphi_{-}(t, \eta)}\Big)
%\widehat{g_0}(\eta)
%d\eta,
%\end{eqnarray*} %%%%%%%%%%%%%%%%%%%%%%%%%%%%%%%%%%%%%%
%where $[0, \infty)={\rre}$, $\widehat{g_0}(\eta) ={{\mathcal
%F}}[f_0](\eta),$ and $\xi=\eta$  for $d=1$. Also
\begin{eqnarray*} %%%%%%%%%%%%%%%%%%%%%%%%%%%%%%%%%%%%
\left(\frac{1}{2\pi}\right)^{d/2}\frac12
Re\int_{{{\mathbb S}}^{d-1}}\int_0^\infty
\Big(e^{i\varphi_{+}(t, \eta)}+
e^{i\varphi_{-}(t, \eta)}\Big)
\widehat{g_0}(\eta)
d\eta d\sigma,
\end{eqnarray*} %%%%%%%%%%%%%%%%%%%%%%%%%%%%%%%%%%%%%%
where $\varphi_\pm$ is given in \eqref{p+-}, $\xi=\eta\omega$, and
\begin{equation} %%%%%%%%%%%%%%%%%%%%%%%%%%%%%%%%%%%%%%%%%%%%%%
\widehat{g_0}(\eta) = \widehat{f_0}(\eta\omega)\eta^{d-1}.  \label{g0}
\end{equation}  %%%%%%%%%%%%%%%%%%%%%%%%%%%%%%%%%%%%%%%%%%%%%%%%
The partial derivatives of $\widehat{g_0}$ are given by
\begin{equation} \begin{array}{rll} %%%%%%%%%%%%%%%%%%%%%%%%%%%%%%%%
&\displaystyle \pt_{\eta}\widehat{g_0}(\eta) =&
\nabla\widehat{f_0}(\eta\omega)\cdot\omega \eta^{d-1} +
(d-1)\widehat{f_0}(\eta\omega) \eta^{d-2}, \\
&\displaystyle \pt_{\eta}^2\widehat{g_0}(\eta) =& \sum_{j=1}^{d} \pt_j
\nabla\widehat{f_0}(\eta\omega)\cdot\omega \omega_j\eta^{d-1} +
2(d-1)\nabla\widehat{f_0}(\eta\omega)\cdot\omega\eta^{d-2} \\
& &+\,\,(d-1)(d-2)\widehat{f_0}(\eta\omega)\eta^{d-3}.
\end{array} \label{g0'} \end{equation} %%%%%%%%%%%%%%%%%%%%%%%%%%%%%%%%

\vskip2mm \noindent \underline{Case: $d=1$.}\,\, Invoking
\eqref{Q0jd1}, \eqref{p123},  \eqref{g0}, and \eqref{g0'},  we have
\begin{eqnarray*} %%%%%%%%%%%%%%%%%%%%%%%%%%%%%%%%%%%%%%%
|Q_{\lambda}^{(0)}(t,x_0)|
\le C
(\lambda t)^{-2} \,
\sum_{j=0}^{2}\|x^jf_0\|_{H_{x}^{j-3}}.
\end{eqnarray*}   %%%%%%%%%%%%%%%%%%%%%%%%%%%%%%%%%%%%%%%
 Combining the above argument, we obtain
\begin{eqnarray*} %%%%%%%%%%%%%%%%%%%%%%%%%%%%%%%%%%%%%%%
|\nabla^{k}Q_{\lambda}^{(0)}(t,x_0)|
\le C
(\lambda t)^{-2}\,
\sum_{j=0}^{2}\|x^jf_0\|_{H_{x}^{m+j-3}},
\end{eqnarray*}  %%%%%%%%%%%%%%%%%%%%%%%%%%%%%%%%%%%%%%%
for $0\le k\le m$.

 %\vskip2mm For the cases of $d=2,\,3$,\, we assume the vanishing condition
%$$(d-1)\eta^{d-2}|\widehat{f_0}(\eta\omega)| +
% (d-1)(d-2)\eta^{d-3}|\omega\cdot\nabla\widehat{f_0}(\eta\omega)|=0$$ at infinity \textrr{and zero}.

\vskip2mm \noindent \underline{Case: $d=2$.}\,\,  Invoking \eqref{Q0j},
\eqref{qj<}, \eqref{g0}, and \eqref{g0'},  we get
\begin{eqnarray*} %%%%%%%%%%%%%%%%%%%%%%%%%%%%%%%%%
|Q_{0,1}|&\lesssim& (\lambda t)^{-2}\big(\sqrt{\log(\lambda t)}
 \|\langle \eta\rangle^{-1} \nabla\widehat{f_0}\, \eta^{\frac12}\|_{L_{\eta}^{2}}+
\sum_{j=1}^{2} \|\pt_j\nabla\widehat{f_0}\, \eta^{\frac12}\|_{L_{\eta}^{2}}\big),\\
|Q_{0,2}|+|Q_{0,3}|&\lesssim& (\lambda t)^{-1-\sigma}
\big(\|\eta^{-\sigma} \widehat{f_0}\, \eta^{\frac12}\|_{L_{\eta}^{2}}+
 \|\langle \eta\rangle^{-2} \widehat{f_0}\, \eta^{\frac12}\|_{L_{\eta}^{2}}
  \\ & &\qquad\qquad \qquad \qquad\qquad \qquad
+\|\langle \eta\rangle^{-1}\nabla\widehat{f_0}\, \eta^{\frac12}\|_{L_{\eta}^{2}}\big),\\
|Q_{0,4}|+|Q_{0,5}| &\lesssim& (\lambda t)^{-1-\sigma}
\big(\|\eta^{-\sigma} \widehat{f_0}\, \eta^{\frac12}\|_{L_{\eta}^{2}}+
 \|\langle \eta\rangle^{-2} \widehat{f_0}\, \eta^{\frac12}\|_{L_{\eta}^{2}}\big),
\end{eqnarray*}  %%%%%%%%%%%%%%%%%%%%%%%%%%%%%%%%%%%%%%%
where $0\leq\sigma<1.$ Combining the above inequalities, we obtain
%\begin{eqnarray*}%%%%%%%%%%%%%%%%%%%%%%%%%%%%%%%%%%%
%|Q_{\lambda}^{(0)}(t,x_0)| \le \int_{{{\mathbb S}}^{1}}
%\sum_{j=1}^5|Q_{0,j}(t,x_0)| d\sigma \le C (\lambda t)^{-\mu}
%\Big(\|f_0\|_{\dot{H}_{x}^{-\sigma}}
%+\sum_{l=0}^{2}\||x|^lf_0\|_{H_{x}^{l-2}}\Big)
%\end{eqnarray*}  %%%%%%%%%%%%%%%%%%%%%%%%%%%%%%%%%%%%%%%
%and
%% where $0\leq\sigma<1,\ 1\leq\mu\leq \min\{\sigma+1,2\}$. Similarly,we have
\begin{eqnarray*} %%%%%%%%%%%%%%%%%%%%%%%%%%%%%%%%%%%%%%%
\lefteqn{|\nabla^{k}Q_{\lambda}^{(0)}(t,x_0)|}\\
&\le& C
(\lambda t)^{-1-\sigma}
\Big(\|f_0\|_{\dot{H}_{x}^{-\sigma}}
+\sum_{l=0}^{1}\||x|^lf_0\|_{H_{x}^{m+l-2}}+
\sum_{j=1}^{2}\|\langle\xi\rangle^{m}
{{\mathcal F}}[xx_j f_0]
\|_{L_{\xi}^{2}}\Big),
\end{eqnarray*}  %%%%%%%%%%%%%%%%%%%%%%%%%%%%%%%%%%%%%%%
for $0\le k\le m$, where $0\leq\sigma<1$ and $\mu=\min\{\sigma+1,
2\}$.

\vskip2mm \noindent \underline{Case: $d=3$.}\,\, Analogously, we have
\begin{eqnarray*} %%%%%%%%%%%%%%%%%%%%%%%%%%%%%%%%%
 |Q_{0,1}|&\lesssim& (\lambda t)^{-\frac32-\sigma}\big(
\| \eta^{-\sigma} \widehat{f_0}\, \eta\|_{L_{\eta}^{2}}+ \|\langle \eta\rangle^{-2}
\widehat{f_0}\, \eta\|_{L_{\eta}^{2}}
+\|\langle \eta\rangle^{-1} \nabla\widehat{f_0}\, \eta\|_{L_{\eta}^{2}}
 \\ & &
+\sum_{j=1}^{3} \|\pt_j\nabla\widehat{f_0}\, \eta\|_{L_{\eta}^{2}}\big),\\
|Q_{0,2}|+|Q_{0,3}| &\lesssim& (\lambda t)^{-\frac32-\sigma}
\big(\|\eta^{-\sigma} \widehat{f_0}\, \eta\|_{L_{\eta}^{2}}+
 \|\langle \eta\rangle^{-2} \widehat{f_0}\, \eta\|_{L_{\eta}^{2}} +
 \|\langle \eta\rangle^{-1}\nabla \widehat{f_0}\, \eta\|_{L_{\eta}^{2}}\big),\\
|Q_{0,4}|+|Q_{0,5}| &\lesssim&  (\lambda t)^{-\frac32-\sigma}
\big(\|\eta^{-\sigma} \widehat{f_0}\, \eta\|_{L_{\eta}^{2}}+
 \|\langle \eta\rangle^{-2} \widehat{f_0}\, \eta\|_{L_{\eta}^{2}}\big),
\end{eqnarray*}  %%%%%%%%%%%%%%%%%%%%%%%%%%%%%%%%%%%%%%%
where $0\leq\sigma<1/2.$ Therefore, we obtain
%\begin{eqnarray*}%%%%%%%%%%%%%%%%%%%%%%%%%%%%%%%%%%%
%|Q_{\lambda}^{(0)}(t,x_0)| \le \int_{{{\mathbb
%S}}^{2}}\sum_{j=1}^5|Q_{3,j}(t,x_0)| d\sigma \le C (\lambda t)^{-\mu}
%\Big(\|f_0\|_{\dot{H}_{x}^{-\sigma}}
%+\sum_{l=0}^{2}\||x|^lf_0\|_{H_{x}^{l-2}}\Big)
%\end{eqnarray*}  %%%%%%%%%%%%%%%%%%%%%%%%%%%%%%%%%%%%%%%
%and
%% where $0\leq\sigma<1/2,\, 1\leq\mu\leq \min\{\sigma+3/2,2\}$. Similarly, we have
\begin{eqnarray*} %%%%%%%%%%%%%%%%%%%%%%%%%%%%%%%%%%%%%%%
\lefteqn{|\nabla^{k}Q_{\lambda}^{(0)}(t,x_0)|}\\
&\le& C
(\lambda t)^{-\frac32-\sigma}
\Big(\|f_0\|_{\dot{H}_{x}^{-\sigma}}
+\sum_{l=0}^{1}\||x|^lf_0\|_{H_{x}^{m+l-2}}+
\sum_{j=1}^{3}\|\langle\xi\rangle^{m}
|{{\mathcal F}}[xx_j f_0]|
\|_{L_{\xi}^{2}}\Big),
\end{eqnarray*}  %%%%%%%%%%%%%%%%%%%%%%%%%%%%%%%%%%%%%%%
for $0\le k\le m$, where $0\leq\sigma<1/2$.
% and $1\leq\mu<\min\{\sigma+3/2, 2\}$.

This completes the proof of Lemma \ref{1f}. $\qed$  %%%%%%%%%%%%%%%% pf of Lem 3.1

\vskip3mm

%\noindent {\bf Proof of Lemma \ref{2f}}. Since %%%%%%%%%%%%%  pf of Lem 3.2
%\begin{eqnarray*} %%%%%%%%%%%%%%%%%%%%%%%%%%%%%%%%%%%%
%n_1+2Im[E_{0}
%\Delta_\varepsilon\overline{E}_{0}]=
%\nabla\cdot\Big\{\phi + 2Im\Big[(E_{0}
%\nabla I_\varepsilon \overline{E}_{0})
%+\varepsilon^{2} \sum_{k=1}^{d}
% (\pt_{k}E_{0}
%\nabla\pt_{k}\overline{E}_{0})\Big]\Big\}\equiv
%\nabla\cdot f_1,
%\end{eqnarray*} %%%%%%%%%%%%%%%%%%%%%%%%%%%%%%%%%%%%
%we have
%\begin{eqnarray*} %%%%%%%%%%%%%%%%%%%%%%%%%%%%%%%%%%%%
%\|Q_{\lambda}^{(1)}(t)\|_{H_{x}^{m}}
%&\le&
%C\lambda^{-1}
%\|f_1
%\|_{H_{x}^{m}}\\
%&\le&C\lambda^{-1}
%(\|n_{1}\|_{H^{m-1}}
%+\|n_{1}\|_{\dot{H}^{-1}}
%+\|E_{0}\|_{H^{m+3}}^{2}).
%\end{eqnarray*} %%%%%%%%%%%%%%%%%%%%%%%%%%%%%%%%%%%%

%%We can rewrite $Q_{\lambda}^{(1)}$ as follows.
%%\begin{equation}\begin{array}{rl}%%%%%%%%%%%%%%%%%%%%%%%%%%%
%%Q_{\lambda}^{(1)}(t, x_0)&\sim\displaystyle \lambda^{-1} Im \int
%%e^{ix_0\xi}
%%\Big(e^{it\lambda\xi_\varepsilon}-e^{-it\lambda\xi_\varepsilon}\Big)
%%\mathcal{F}[\omega_\varepsilon^{-1} \nabla\cdot f_1](\xi)\,d\xi\\
%%&= \displaystyle \lambda^{-1} Im \int
%%\Big(e^{i\varphi_+(\xi)}-e^{-i\varphi_-(\xi)}\Big)
%%\mathcal{F}[\omega_\varepsilon^{-1} \nabla\cdot f_1](\xi)\,d\xi.
%%\end{array}\label{Q2-3}
%%\end{equation}   %%%%%%%%%%%%%%%%%%%%%%%%%%%%
%%Analogous to the proof of Lemma \ref{1f}, we can prove that
%%\eqref{Q2-2} holds. This completes the proof of Lemma \ref{2f}.
%%$\qed$
%%%%%%%%%%%%%%%%%%%%%%%%%%%%%%%%%%%%%%%%%%%%%%%%%%%%%%

\vskip2mm \noindent {\bf Proof of Lemma \ref{2ff}}. From (\ref{i3}) and
(\ref{Q0-1}), we have
\begin{eqnarray*} %%%%%%%%%%%%%%%%%%%%%%%%%%%%%%%%%%%%
\|Q_{\lambda}^{(1)}(t)\|_{H_{x}^{m}}
&\le&
C\lambda^{-1}
\|f_1
\|_{H_{x}^{m}}\\
&\le&C\lambda^{-1}
(\|n_{1}\|_{H^{m-1}}
+\|n_{1}\|_{\dot{H}^{-1}}
+\|E_{0}\|_{H^{m+3}}^{2}).\qquad\qed
\end{eqnarray*} %%%%%%%%%%%%%%%%%%%%%%%%%%%%%%%%%%%%

\vskip2mm
\noindent {\bf Proof of Lemma \ref{2f}}.
The proof is analogous to that of Lemma \ref{1f}.  Now
\begin{equation} %%%%%%%%%%%%%%%%%%%%%%%%%%%%%%%%%
Q_{\lambda}^{(1)}(t, x_0)
=
\left(\frac{1}{2\pi}\right)^{d/2}
\frac{-1}{2\lambda}Im\left[
\int_{\rre^{d}}
e^{ix_0\cdot\xi}\big(e^{i\lambda t\xi_\varepsilon}
- e^{i\lambda t\xi_\varepsilon}\big)
\frac{\xi}{\xi_\varepsilon}\cdot\widehat{f_1}(\xi)
 d\xi\right],\label{Q(1)}
\end{equation} %%%%%%%%%%%%%%%%%%%%%%%%%%%%%%%%%
where $\xi_\varepsilon$ is given as in \eqref{xi_ep}. For the case
$0\le\lambda t\le1$, for any $0\le k\le m$, we have
\begin{eqnarray*}%%%%%%%%%%%%%%%%%%%%%
|\nabla^{k}Q_{\lambda}^{(1)}(t,x_0)|
&\le &C \lambda^{-1}
\|f_1\|_{H^{m+[d/2]}}.
\end{eqnarray*}%%%%%%%%%%%%%%%%%%%%%%
Next for the case $\lambda t>1$ and $|x_0|\ge \lambda t/2$,
we have
\begin{eqnarray*} %%%%%%%%%%%%%%%%%%%%%%%%%%%%%%%%%
|\nabla^{k}
Q_{\lambda}^{(1)}(t,x_0)|
&\le&C \lambda^{-1}
(\lambda t)^{-2}|x_0|^{2}\|f_1\|_{H^{m+[d/2]}}.
\end{eqnarray*} %%%%%%%%%%%%%%%%%%%%%%%%%%%%%%%%%

 Let us evaluate $Q_{\lambda}^{(1)}(t, x_0)$ in the
region where $|x_0|\le \lambda t/2$, and rewrite it as follows
%\begin{eqnarray*} %%%%%%%%%%%%%%%%%%%%%%%%%%%%%%%%%%%%
%\left(\frac{1}{2\pi}\right)^{1/2}\frac1{2\lambda}
%Im\int_0^\infty
%\Big(e^{i\varphi_{+}(t, \eta)}-
%e^{-i\varphi_{-}(t, \eta)}\Big)
%\widehat{g_1}(\eta)
%d\eta,
%\end{eqnarray*} %%%%%%%%%%%%%%%%%%%%%%%%%%%%%%%%%%%%%%
%where $[0, \infty)={\rre}$, $\widehat{g_1}(\eta) ={{\mathcal
%F}}[f_1](\eta\omega)/\sqrt{1+\varepsilon^2\eta^2} $, and $\xi=\eta$ for
%$d=1$. For $d=2, 3$, we have
\begin{eqnarray*} %%%%%%%%%%%%%%%%%%%%%%%%%%%%%%%%%%%%
\left(\frac{1}{2\pi}\right)^{d/2}\frac1{2\lambda}
Im\int_{{{\mathbb S}}^{d-1}}\int_0^\infty
\Big(e^{i\varphi_{+}(t, \eta)}-
e^{i\varphi_{-}(t, \eta)}\Big)
\widehat{g_1}(\eta)
d\eta d\sigma,
\end{eqnarray*} %%%%%%%%%%%%%%%%%%%%%%%%%%%%%%%%%%%%%%
where $\varphi_\pm$ is given in \eqref{p+-} $\xi=\eta\omega$, and
\begin{equation} %%%%%%%%%%%%%%%%%%%%%%%%%%%%%%%%%%%%%%%%%%%%%%
\widehat{g_1}(\eta) = \omega\cdot\widehat{f_1}(\eta\omega)\eta^{d-1}/\sqrt{1+\varepsilon^2\eta^2}. \label{g1}
\end{equation}  %%%%%%%%%%%%%%%%%%%%%%%%%%%%%%%%%%%%%%%%%%%%%%%%

\vskip2mm

% Invoking \eqref{p123=}, \eqref{e^ip}, and repeating the
%integration by parts twice, we have
%\begin{equation} \begin{array}{rll}   %%%%%%%%%%%%%%%%%%%%%%%%%%%%%%%%
%&\displaystyle Im\int_0^\infty e^{i\varphi_{\pm}(t, \eta)}\widehat{g_1}(\eta)d\eta \\
%=&\displaystyle Im\int_0^\infty
%e^{i\varphi_{\pm}(t, \eta)}q_{1}\pt_{\eta}^{2}{{\mathcal
%F}}[g_1](\eta)d\eta \,\,+\,\, Im\int_0^\infty
%e^{i\varphi_{\pm}(t, \eta)}q_{2}\pt_{\eta}\widehat{g_1}(\eta)d\eta \\
%&\displaystyle +\,\, 3Re\int_0^\infty e^{i\varphi_{\pm}(t,
%\eta)}q_{3}\pt_{\eta}\widehat{g_1}(\eta)d\eta\,\,
 %+\,\, Re\int_0^\infty
%e^{i\varphi_{\pm}(t, \eta)}q_{4}\widehat{g_1}(\eta)d\eta \\
%&\displaystyle -\,\, 3Im\int_0^\infty
%e^{i\varphi_{\pm}(t, \eta)}q_{5}\widehat{g_1}(\eta)d\eta \\
%\equiv&\displaystyle  \,\,Q_{1,1}(t, x_0)+Q_{1,2}(t,x_0)+Q_{1,3}(t,x_0)+Q_{1,4}(t,x_0)+Q_{1,5}(t,x_0), \\
%\end{array} \label{Q1j}  \end{equation}   %%%%%%%%%%%%%%%%%%%%%%%%%%%%
%where $q_{j}$ is given in \eqref{qj=} for $j=1,\,\cdot\cdot\cdot,\, 5$.

  For any $\eta\in [0, \infty)$ and in the region $|x_0|\le(\lambda
t)/2$, we estimate the quantities $Q_{1,j}$ for $j=1, \cdot\cdot\cdot, 5$.
We split the interval $[0, \infty)$ into three parts  $I_1(t)$, $I_2(t)$, and
$I_3(t)$ which are given in \eqref{I123}. Let $f_1=\big(f_{1,
1},\cdot\cdot\cdot, f_{1, d}\big).$
%we assume the vanishing condition
%$$(d-1)r^{d-2}|\widehat{f_0}(r\omega)| +
%(d-1)(d-2)r^{d-3}|\omega\cdot\nabla\widehat{f_0}(r\omega)|=0$$
%at infinity. Then
We estimate the partial derivatives of $\widehat{g_1}$,
\begin{equation} \begin{array}{rl} %%%%%%%%%%%%%%%%%%%%%%%%%%%%%%%%%%%%%%%%%%%%%%
\big|\partial_{\eta}\widehat{g_1}({\eta})\big|\lesssim&\displaystyle

\sum_{j=1}^{d} \frac{|\nabla\widehat{f_{1,j}}|}{\langle\eta\rangle}
\eta^{d-1}+ \big(|d-1|+|d-2|\eta^2\big)
\frac{|\widehat{f_{1}}|}{\langle\eta\rangle^3}\eta^{d-2}, \\
\displaystyle \big|\partial_{\eta}^2\widehat{g_1}({\eta})\big|\lesssim&
\displaystyle \sum_{j,\,\ell=1}^{d}
\frac{|\partial_{\ell}\nabla\widehat{f_{1,j}}|}{\langle\eta\rangle}\eta^{d-1}
+\,\, \big(|d-1|+|d-2|\eta^2\big)
\sum_{j=1}^{d}\frac{|\nabla\widehat{f_{1,j}}|}{\langle\eta\rangle^3}\eta^{d-2} \\
&\displaystyle +\,\, \big(|d-1||d-2|+\eta^2+|d-2||d-3|\eta^4\big)
\frac{|\widehat{f_{1}}|}{\langle\eta\rangle^5}\eta^{d-3}.
\end{array} \label{g1'-g1"}  \end{equation}    %%%%%%%%%%%%%%%%%%%%%%%%%%%%

\vskip2mm \noindent \underline{Case: $d=1$.}\,\, Invoking
\eqref{Q0jd1}, \eqref{p123}, \eqref{g1} - \eqref{g1'-g1"},  we have
\begin{eqnarray*} %%%%%%%%%%%%%%%%%%%%%%%%%%%%%%%%%%%%%%%
|Q_{\lambda}^{(1)}(t,x_0)|
\le C
\lambda^{-1}(\lambda t)^{-2} \,
\sum_{j=0}^{2}\|x^jf_1\|_{H_{x}^{j-4}}.
\end{eqnarray*}   %%%%%%%%%%%%%%%%%%%%%%%%%%%%%%%%%%%%%%%
Combining the above argument, we obtain
\begin{eqnarray*} %%%%%%%%%%%%%%%%%%%%%%%%%%%%%%%%%%%%%%%
|\nabla^{k}Q_{\lambda}^{(1)}(t,x_0)|
\le C
\lambda^{-1}(\lambda t)^{-2}\,
\sum_{j=0}^{2}\|x^jf_1\|_{H_{x}^{m+j-4}},
\end{eqnarray*}  %%%%%%%%%%%%%%%%%%%%%%%%%%%%%%%%%%%%%%%
for $0\le k\le m$.

\vskip2mm \noindent \underline{Case: $d=2$.}\,\,  Analogously we get
\begin{eqnarray*}%%%%%%%%%%%%%%%%%%%%%%%%%%%%%%%%%%%%%%%%%%%
|Q_{1,1}|&\lesssim&(\lambda t)^{-2}
\displaystyle\Big(\sum_{j,\,\ell=1}^{2}
\|\langle\eta\rangle^{-1}\partial_{ \ell}\nabla\widehat{f_{1,j}}\eta^{\frac12}\|_{L_{\eta}^{2}}+
\displaystyle\sum_{j=1}^{2}
\|\langle\eta\rangle^{-2}\nabla \widehat{f_{1,j}}\eta^{\frac12}\|_{L_{\eta}^{2}}\\
& &\qquad\qquad+\|\langle\eta\rangle^{-3}\widehat{f_{1}}\eta^{\frac12}\|_{L_{\eta}^{2}}\Big), \\
|Q_{1,2}|+|Q_{1,3}|&\lesssim&
(\lambda t)^{-1-\sigma}   \Big( \displaystyle\sum_{j=1}^{2}
\|\langle\eta\rangle^{-2}\nabla \widehat{f_{1,j}}\eta^{\frac12}\|_{L_{\eta}^{2}}+
\|\langle\eta\rangle^{-3}\widehat{f_{1}}\eta^{\frac12}\|_{L_{\eta}^{2}}\\
& &\qquad\qquad+
\|\eta^{-\sigma}\widehat{f_{1}}\eta^{\frac12}\|_{L_{\eta}^{2}} \Big),\\
|Q_{1,4}|+|Q_{1,5}|&\lesssim&
 (\lambda t)^{-1-\sigma}
\Big(\|\langle\eta\rangle^{-3}\widehat{f_{1}}\eta^{\frac12}\|_{L_{\eta}^{2}}+
\|\eta^{-\sigma}\widehat{f_{1}}\eta^{\frac12}\|_{L_{\eta}^{2}} \Big),
\end{eqnarray*}%%%%%%%%%%%%%%%%%%%%%%%%%%%%%%%%%%%%%%%%%%%%%%
where $0\leq\sigma<1.$ Combining the above inequalities, we obtain
%\begin{eqnarray*}%%%%%%%%%%%%%%%%%%%%%%%%%%%%%%%%%%%
%|Q_{\lambda}^{(1)}(t,x_0)| \le \int_{{{\mathbb S}}^{1}}
%\sum_{j=1}^5|Q_{2,j}(t,x_0)| d\sigma \le C\lambda^{-1}(\lambda t)^{-\mu}
%\Big(
%\|f_1\|_{\dot{H}_x^{-\sigma}} + \sum_{\ell=0}^{2}\||x|^\ell f_1\|_{H_x^{\ell-3}}
%\Big)
%\end{eqnarray*}  %%%%%%%%%%%%%%%%%%%%%%%%%%%%%%%%%%%%%%%
%and
\begin{eqnarray*} %%%%%%%%%%%%%%%%%%%%%%%%%%%%%%%%%%%%%%%
|\nabla^{k}Q_{\lambda}^{(1)}(t,x_0)|
\le C
\lambda^{-1}(\lambda t)^{-1-\sigma}
\Big(\|f_1\|_{\dot{H}_{x}^{-\sigma}}
+\sum_{\ell=0}^{2}\||x|^\ell f_1\|_{H_{x}^{m+\ell-3}}\Big),
\end{eqnarray*}  %%%%%%%%%%%%%%%%%%%%%%%%%%%%%%%%%%%%%%%
for $0\le k\le m$ and $0\leq\sigma<1$.
% and$1\leq\mu<\min\{\sigma+1, 2\}$.

\vskip2mm \noindent \underline{Case: $d=3$.}\,\,  Analogously we have
\begin{eqnarray*}%%%%%%%%%%%%%%%%%%%%%%%%%%%%%%%%%%%%%%%%%%%
|Q_{1,1}|&\lesssim&
 (\lambda t)^{-\frac32-\sigma}
\displaystyle\Big(\sum_{j,\,\ell=1}^{3}
\|\langle\eta\rangle^{-1}\partial_{ \ell}\nabla\widehat{f_{1,j}}\eta\|_{L_{\eta}^{2}}+
\displaystyle\sum_{j=1}^{3}
\|\langle\eta\rangle^{-2}\nabla \widehat{f_{1,j}}\eta\|_{L_{\eta}^{2}}\\
& &\qquad\qquad+
\|\langle\eta\rangle^{-3}\widehat{f_{1}}\eta\|_{L_{\eta}^{2}}\Big), \\
|Q_{1,2}|+|Q_{1,3}|&\lesssim&
 (\lambda t)^{-\frac32-\sigma}\Big( \displaystyle\sum_{j=1}^{3}
\|\langle\eta\rangle^{-2}\nabla \widehat{f_{1,j}}\eta\|_{L_{\eta}^{2}}
+
\|\langle\eta\rangle^{-3}\widehat{f_{1}}\eta\|_{L_{\eta}^{2}}
\\
& &\qquad\qquad+
\|\eta^{-\sigma}\widehat{f_{1}}\eta\|_{L_{\eta}^{2}} \Big),\\
|Q_{1,4}|+|Q_{1,5}|&\lesssim&
 (\lambda t)^{-\frac32-\sigma}
\Big(\|\langle\eta\rangle^{-3}\widehat{f_{1}}\eta\|_{L_{\eta}^{2}}+
\|\eta^{-\sigma}\widehat{f_{1}}\eta\|_{L_{\eta}^{2}} \Big),
\end{eqnarray*}%%%%%%%%%%%%%%%%%%%%%%%%%%%%%%%%%%%%%%%%%%%%%%
where $0\leq\sigma<1/2.$ 
Note that the extra decay rate $t^{1/2}$ compared 
to $d=2$ is due to the increase of the space dimension. 
Combining the above inequalities, we obtain
%\begin{eqnarray*}%%%%%%%%%%%%%%%%%%%%%%%%%%%%%%%%%%%
%|Q_{\lambda}^{(1)}(t,x_0)| \le \int_{{{\mathbb S}}^{1}}
%\sum_{j=1}^5|Q_{2,j}(t,x_0)| d\sigma \le C\lambda^{-1}(\lambda t)^{-\mu}
%\Big(
%\|f_1\|_{\dot{H}_x^{-\sigma}} + \sum_{\ell=0}^{2}\||x|^\ell f_1\|_{H_x^{\ell-3}}
%\Big)
%\end{eqnarray*}  %%%%%%%%%%%%%%%%%%%%%%%%%%%%%%%%%%%%%%%
%and
\begin{eqnarray*} %%%%%%%%%%%%%%%%%%%%%%%%%%%%%%%%%%%%%%%
|\nabla^{k}Q_{\lambda}^{(1)}(t,x_0)|
\le C
\lambda^{-1}(\lambda t)^{-\frac32-\sigma}
\Big(\|f_1\|_{\dot{H}_{x}^{-\sigma}}
+\sum_{\ell=0}^{2}\||x|^\ell f_1\|_{H_{x}^{m+\ell-3}}\Big),
\end{eqnarray*}  %%%%%%%%%%%%%%%%%%%%%%%%%%%%%%%%%%%%%%%
for $0\le k\le m$, where $0\leq\sigma<1/2$.
% and $1\leq\mu<\min\{\sigma+3/2, 2\}$.

This completes the proof of Lemma \ref{2f}. $\qed$      %%%%%%%%%%%%%% pf of Lem 3.2
%%%%%%%%%%%%%%%%%%%%%%%%%%%%%%%%%%%%%%%%%%%%%%%%%%%%%%

\vskip2mm \noindent {\bf Proof of Lemma \ref{3ff}}. For the sake of
convenience, we drop the indices of $E_{\lambda}$ and $n_{\lambda}$.

From (\ref{ZK}), we see
\begin{equation} \begin{array}{rll}   %%%%%%%%%%%%%%%%%%%%%%%%%%%%%%%%
\displaystyle \pt_{t}^{2}|E|^{2}=& -2\pt_{t}Im[\overline{E}\Delta_\varepsilon E] \\
=&\displaystyle  -2\pt_{t}Im[\nabla\cdot\{\overline{E} I_\varepsilon
\nabla E\}] -2\varepsilon^{2} \pt_{t}\sum_{k=1}^{d}Im
[\nabla\cdot(\pt_{k}\overline{E}\nabla\pt_{k}E)] \\
\displaystyle =& 2\nabla\cdot Re[\{ \Delta_\varepsilon \overline{E} -n
\overline{E}\} I_\varepsilon \nabla E + \overline{E} I_\varepsilon \nabla
\{-\Delta_\varepsilon E +nE\}] \\
\displaystyle &+2\varepsilon^{2} \nabla\cdot Re\sum_{k=1}^{d}
\Big[\pt_{k}\{ \Delta_\varepsilon \overline{E} -n \overline{E}\}
\nabla\pt_{k}E+\pt_{k}\overline{E} \nabla\pt_{k}\{ -\Delta_\varepsilon E +n E\}\Big] \\
\displaystyle \equiv& \nabla\cdot f_2.
\end{array} \label{f2}  \end{equation}   %%%%%%%%%%%%%%%%%%%%%%%%%%%%
The Sobolev embedding yields
\begin{eqnarray*} %%%%%%%%%%%%%%%%%%%%%% %%%
\|(-\Delta)^{-1/2}
I_\varepsilon^{-3/2}
\pt_{t}^{2}|E|^{2}
\|_{H_{x}^{m}}
\le C\|f_{2}\|_{H_{x}^{m-3}}
%%%%%%%%%%%%%%%%%%%%%%
%&\le&C
%\|E\|_{H_{x}^{m+1}}
%\|E\|_{H_{x}^{m}}
%+C
%\|n\|_{H_{x}^{m-3}}
%\|E\|_{H_{x}^{m-3}}
%\|E\|_{H_{x}^{m}}
%+C
%\|E\|_{H_{x}^{m-3}}
%\|E\|_{H_{x}^{m+4}}
%\\ %%%%%%%%%%%%%%%%%%%%%%
%& &
%+C \|n\|_{H_{x}^{m}} \|E\|_{H_{x}^{m-3}} \|E\|_{H_{x}^{m}} +C
%\|E\|_{H_{x}^{m+2}} \|E\|_{H_{x}^{m-1}} +C \|n\|_{H_{x}^{m-2}}
%\|E\|_{H_{x}^{m-2}} \|E\|_{H_{x}^{m-1}}
%\\ %%%%%%%%%%%%%%%%%%%%%%
%& &
%+C
%\|E\|_{H_{x}^{m-2}}
%\|E\|_{H_{x}^{m+3}}
%+C
%\|n\|_{H_{x}^{m-1}}
%\|E\|_{H_{x}^{m-2}}
%\|E\|_{H_{x}^{m-1}}
%\\ %%%%%%%%%%%%%%%%%%%%%%
\le C
(1+\|n\|_{H_{x}^{m}})
\|E\|_{H_{x}^{m+4}}^{2},
\end{eqnarray*} %%%%%%%%%%%%%%%%%%%%%% %%%%%
for $m\ge3$.
% if $d=1$ and $m\ge5$ if $d=2,3$, where we used that $H^{m-3}$
%is algebra.
Hence we have
\begin{eqnarray*} %%%%%%%%%%%%%%%%%%%%%%%%%%%%%%%%%
\|Q_{\lambda}^{(2)}\|_{H_{x}^{m}}
&\le&C
\int_{0}^{T^{\ast}}
\lambda^{-1}
\|(-\Delta)^{-1/2}I_\varepsilon^{-3/2}
\pt_{t}^{2}
|E|^{2}
\|_{H_{x}^{m}}ds\\
&\le&C\lambda^{-1}T
\sup_{0\le t\le T}
\{(1+\|n(t)\|_{H_{x}^{m}})
\|E(t)\|_{H_{x}^{m+4}}^{2}\},
%\\
%&\le&C\lambda^{-1},
\end{eqnarray*} %%%%%%%%%%%%%%%%%%%%%%%%%%%%%%%%%
where the constant $C$ depends on $T^{\ast}$ but is independent of
$\lambda\in[1,\infty)$.
%We can rewrite $Q_{\lambda}^{(2)}$ in a way analogous to
%$Q_{\lambda}^{(1)}$, Following the proof of Lemma \ref{1f},
%we can prove that \eqref{Q3-2} holds.
 This completes the proof of Lemma \ref{3ff}. $\qquad\qed$

\vskip2mm \noindent {\bf Proof of Lemma \ref{3f}}. %%%%%%%%%%%%%% pf of Lem 3.3
Now we rewrite $Q_{\lambda}^{(2)}$ as follows.
\begin{equation} %%%%%%%%%%%%%%%%%%%%%%%%%%%%%%%%%
Q_{\lambda}^{(2)}(t, x_0)
=C
%=\left(\frac{1}{2\pi}\right)^{d/2}\frac{-1}{2\lambda}
\lambda^{-1}
Im\left[
\int_0^t \int_{\rre^{d}}
e^{ix_0\cdot\xi}\big(e^{i\lambda (t-s)\xi_\varepsilon}
- e^{-i\lambda (t-s)\xi_\varepsilon}\big)
\frac{\xi\cdot\widehat{f_2}(s,
\xi)}{\xi_\varepsilon(1+\varepsilon^2\xi^2)}
 d\xi ds\right],\label{Q(2)}
\end{equation} %%%%%%%%%%%%%%%%%%%%%%%%%%%%%%%%%
where $\xi_\varepsilon$ is given in \eqref{xi_ep} and $f_2$ is
given in \eqref{f2}. For the case $0\le
t\le\lambda^{-1}\max\{1,2|x_{0}|\}$,
%$0\le\lambda (t-s)\le1$,
for any $0\le k\le m$, we have
\begin{eqnarray}%%%%%%%%%%%%%%%%%%%%%
\lefteqn{|\nabla^{k}Q_{\lambda}^{(2)}(t,x_0)|}\nonumber\\
&\le &C \lambda^{-2}(1+|x_{0}|)\sup_{0\le t\le T}
\|f_2(t)\|_{H^{m-2+[d/2]}}\nonumber\\
&\le &C \lambda^{-2}(1+|x_{0}|)\sup_{0\le t\le T}
\left\{
(1+\|n_{\lambda}(t)\|_{H_{x}^{m+4+[d/2]}})
\|E_{\lambda}(t)\|_{H_{x}^{m+5+[d/2]}}^{2}\right\}.
\nonumber\\
\label{y10}
\end{eqnarray}%%%%%%%%%%%%%%%%%%%%%%
%Next we consider the case
%$\lambda (t-s)>1$ and $|x_0|\ge \lambda(t-s)/2$.
%\begin{eqnarray*}
%|\nabla^{k}
%Q_{\lambda}^{(2)}(t,x_0)|
%&\le&C \lambda^{-1}
%(\lambda t)^{-2}|x_0|^{2}\sup_{[0, T^{\ast}]}\|f_2(s)\|_{H^{m-2+[d/2]}}.
%\end{eqnarray*}
In the region $t\ge\lambda^{-1}\max\{1,2|x_{0}|\}$, we split
$Q_{\lambda}^{(2)}$ into the following two pieces
\begin{eqnarray*}%%%%%%%%%%%%%%%%%%%%%%%%%%%%%%%%%%%%
&&Q_{\lambda}^{(2)}(t, x_0)\\
&\sim&
\lambda^{-1}
Im\biggl[
\Big(\int_0^{t-\lambda^{-1}{b}}
+\int_{t-\lambda^{-1}{b}}^{t}\Big)
\int_{\rre^{d}}
e^{ix_0\cdot\xi}\big(e^{i\lambda (t-s)\xi_\varepsilon}
- e^{-i\lambda (t-s)\xi_\varepsilon}\big)\\
& &\qquad\qquad\times
\frac{\xi\cdot\widehat{f_2}(s,
\xi)}{\xi_\varepsilon(1+\varepsilon^2\xi^2)}
 d\xi ds\biggl]\\
&\equiv&Q_{\lambda}^{(2,1)}(t, x_0)+Q_{\lambda}^{(2,2)}(t, x_0),
\end{eqnarray*}%%%%%%%%%%%%%%%%%%%%%%%%%%%%%%%%%%%%
where $b=\max\{1,2|x_{0}|\}$.  For $Q_{\lambda}^{(2,2)}$, we
can easily see that
\begin{eqnarray}%%%%%%%%%%%%%%%%%%%%%
\lefteqn{|\nabla^{k}Q_{\lambda}^{(2,2)}(t,x_0)|}\nonumber\\
&\le &C \lambda^{-2}(1+|x_{0}|)\sup_{0\le t\le T}
\|f_2(t)\|_{H^{m-2+[d/2]}}\nonumber\\
&\le &C \lambda^{-2}(1+|x_{0}|)\sup_{0\le t\le T}
\left\{
(1+\|n_{\lambda}(t)\|_{H_{x}^{m+4+[d/2]}})
\|E_{\lambda}(t)\|_{H_{x}^{m+5+[d/2]}}^{2}\right\},
\nonumber\\
\label{y11}
\end{eqnarray}%%%%%%%%%%%%%%%%%%%%%%
where $0\le k\le m$.

%$|x_0|\le \lambda (t-s)/2$,
We rewrite $Q_{\lambda}^{(2,1)}(t, x_0)$ as
%\begin{eqnarray*} %%%%%%%%%%%%%%%%%%%%%%%%%%%%%%%%%%%%
%\left(\frac{1}{2\pi}\right)^{1/2}\frac1{2\lambda}
%Im\int_0^\infty
%\Big(e^{i\varphi_{+}(t, \eta)}-
%e^{-i\varphi_{-}(t, \eta)}\Big)
%\widehat{g_1}(\eta)
%d\eta,
%\end{eqnarray*} %%%%%%%%%%%%%%%%%%%%%%%%%%%%%%%%%%%%%%
%where $[0, \infty)={\rre}$, $\widehat{g_1}(\eta) ={{\mathcal
%F}}[f_1](\eta\omega)/\sqrt{1+\varepsilon^2\eta^2} $, and $\xi=\eta$ for
%$d=1$. For $d=2, 3$, we have
\begin{eqnarray*} %%%%%%%%%%%%%%%%%%%%%%%%%%%%%%%%%%%%
C\lambda^{-1}
Im\int_0^{t-\lambda^{-1}{b}}\int_{{{\mathbb S}}^{d-1}}\int_0^\infty
\Big(e^{i\varphi_{+}(t-s, \eta)}-
e^{i\varphi_{-}(t-s, \eta)}\Big)
\widehat{g_2}(\eta)
d\eta\, d\sigma\, ds,
\end{eqnarray*} %%%%%%%%%%%%%%%%%%%%%%%%%%%%%%%%%%%%%%
where $\varphi_\pm$ is given in \eqref{p+-}, $\xi=\eta\omega$, and
\begin{equation} %%%%%%%%%%%%%%%%%%%%%%%%%%%%%%%%%%%%%%%%%%%%%%
\widehat{g_2}(\eta) = \eta^{d-1}(1+\varepsilon^2\eta^2)^{-3/2}
\omega\cdot\widehat{f_2}(s, \eta\omega). \label{g2}
\end{equation} %%%%%%%%%%%%%%%%%%%%%%%%%%%%%%%%%%%%%%%%%%%%%%%%

\vskip2mm

% Invoking \eqref{p123=}, \eqref{e^ip}, and repeating the
%integration by parts twice, we have
%\begin{equation} \begin{array}{rll}   %%%%%%%%%%%%%%%%%%%%%%%%%%%%%%%%
%&\displaystyle Im\int_0^\infty e^{i\varphi_{\pm}(t-s, \eta)}\widehat{g_2}(\eta)d\eta \\
%=&\displaystyle Im\int_0^\infty
%e^{i\varphi_{\pm}(t-s, \eta)}q_{1}\pt_{\eta}^{2}{{\mathcal
%F}}[g_2](\eta)d\eta \,\,+\,\, Im\int_0^\infty
%e^{i\varphi_{\pm}(t-s, \eta)}q_{2}\pt_{\eta}\widehat{g_2}(\eta)d\eta \\
%&\displaystyle +\,\, Re\int_0^\infty
%e^{i\varphi_{\pm}(t-s, \eta)}q_{3}\pt_{\eta}{{\mathcal
%F}}[g_2](\eta)d\eta\,\,
% +\,\, Re\int_0^\infty
%e^{i\varphi_{\pm}(t-s, \eta)}q_{4}\widehat{g_2}(\eta)d\eta \\
%&\displaystyle -\,\, Im\int_0^\infty
%e^{i\varphi_{\pm}(t-s, \eta)}q_{5}\widehat{g_2}(\eta)d\eta \\
%\equiv&\displaystyle  \,\,Q_{2,1}(t, x_0)+Q_{2,2}(t,x_0)+Q_{2,3}(t,x_0)+Q_{2,4}(t,x_0)+Q_{2,5}(t,x_0),
%\end{array} \label{Q2j}  \end{equation}   %%%%%%%%%%%%%%%%%%%%%%%%%%%%
%where $q_{j}$ are given in \eqref{qj=} for $j=1,\,\cdot\cdot\cdot,\, 5$.

  For any $\eta\in [0, \infty)$ and in the region $|x_0|\le\lambda
(t-s)/2$, we estimate the quantities $Q_{2,j}$ for $j=1, \cdot\cdot\cdot,
5$. We split the interval $[0, \infty)$ into three parts  $I_1(t-s)$,
$I_2(t-s)$, and $I_3(t-s)$ which are given in \eqref{I123}. We then
estimate the partial derivatives of $\widehat{g_2}$,
\begin{equation} \begin{array}{rl}  %%%%%%%%%%%%%%%%%%%%%%%%%%%%%%%%%%%%%%%%%%%%%%
\big|\partial_{\eta}\widehat{g_2}({\eta})\big|\lesssim&\displaystyle
\frac{\eta^{d-1} |\nabla\widehat{f_{2}}|}{\langle\eta\rangle^{3}}
 + \big(|d-1|+\eta^2\big)
\frac{\eta^{d-2}|\widehat{f_{2}}|}{\langle\eta\rangle^{3}},\\
\big|\partial_{\eta}^2\widehat{g_2}({\eta})\big|\lesssim&\displaystyle
\displaystyle \sum_{j=1}^{d} \frac{\eta^{d-1}|\partial_{j}\nabla\widehat{f_{2}}|}{\langle\eta\rangle^{3}}+\,\,
\big(|d-1|+\eta^2\big)
\frac{\eta^{d-2}|\nabla\widehat{f_{2}}|}{\langle\eta\rangle^{5}}\\
& +\,\,  \displaystyle
\big(|d-1||d-2|+\eta^2+\eta^4\big)
\frac{\eta^{d-3}|\widehat{f_{2}}|}{\langle\eta\rangle^{7}}.
\end{array} \label{g2'-g2"}  \end{equation}   %%%%%%%%%%%%%%%%%%%%%%%%%%%%

\vskip2mm \noindent \underline{Case: $d=1$.}\,\,   Invoking
\eqref{Q0jd1}, \eqref{p123},   \eqref{g2} - \eqref{g2'-g2"},  we have
\begin{eqnarray} %%%%%%%%%%%%%%%%%%%%%%%%%%%%%%%%%%%%%%%
|Q_{\lambda}^{(2,1)}(t,x_0)|
&\le& C
\lambda^{-1}\int_0^{t-\lambda^{-1}{b}}(1+\lambda (t-s))^{-2}
\sum_{\ell=0}^{2}\|x^{\ell}f_2(s)\|_{H_{x}^{\ell-6}}\,ds\nonumber\\
&\le& C
\lambda^{-2}T^{1/2}\sup_{0\le t\le T}\,
\sum_{\ell=0}^{2}\|x^{\ell}f_2(t)\|_{H_{x}^{\ell-6}}.
\end{eqnarray}   %%%%%%%%%%%%%%%%%%%%%%%%%%%%%%%%%%%%%%%
Combining the above argument, for $0\le k\le m$, we obtain
\begin{eqnarray} %%%%%%%%%%%%%%%%%%%%%%%%%%%%%%%%%%%%%%%
|\nabla^{k}Q_{\lambda}^{(2,1)}(t,x_0)|
&\le& C
\lambda^{-1}\int_0^{t-\lambda^{-1}{b}}(1+\lambda (t-s))^{-2}
\sum_{\ell=0}^{2}\|x^lf_2(s)\|_{H_{x}^{m+\ell-6}}\,ds\nonumber\\
&\le& C
\lambda^{-2}(T)^{1/2}\sup_{0\le t\le T}\,
\sum_{\ell=0}^{2}\|x^lf_2(t)\|_{H_{x}^{m+\ell-6}}.
\label{y1}
\end{eqnarray}  %%%%%%%%%%%%%%%%%%%%%%%%%%%%%%%%%%%%%%%
\vskip2mm
\noindent \underline{Case: $d=2$.}\,\,  Analogously we get
\begin{eqnarray*}%%%%%%%%%%%%%%%%%%%%%%%%%%%%%%%%%%%%%%%%%%%
|Q_{2,1}|&\lesssim&(\lambda (t-s))^{-2}\sqrt{\log (\lambda(t-s))}\\
& &\times\Big(\sum_{j=1}^{2} \|\langle\eta\rangle^{-3}\partial_{
j}\nabla\widehat{f_{2}}\eta^{\frac12}\|_{L_{\eta}^{2}}+
\|\langle\eta\rangle^{-4}\nabla \widehat{f_{2}}\eta^{\frac12}\|_{L_{\eta}^{2}}\|\langle\eta\rangle^{-5}\widehat{f_{2}}\eta^{\frac12}\|_{L_{\eta}^{2}}\Big), \\
|Q_{2,2}|+|Q_{2,3}|&\lesssim&
(\lambda(t-s))^{-1}\Big( \displaystyle
\|\langle\eta\rangle^{-4}\nabla \widehat{f_{2}}\eta^{\frac12}\|_{L_{\eta}^{2}}+
\|\langle\eta\rangle^{-5}\widehat{f_{2}}\eta^{\frac12}\|_{L_{\eta}^{2}}\Big),\\
|Q_{2,4}|+|Q_{2,5}|&\lesssim&
(\lambda(t-s))^{-1} \|\langle\eta\rangle^{-5}\widehat{f_{2}}\eta^{\frac12}\|_{L_{\eta}^{2}}.
\end{eqnarray*}%%%%%%%%%%%%%%%%%%%%%%%%%%%%%%%%%%%%%%%%%%%%%%
Combining the above inequalities, for $0\le k\le m$, we obtain
%\begin{eqnarray*}%%%%%%%%%%%%%%%%%%%%%%%%%%%%%%%%%%%
%|Q_{\lambda}^{(1)}(t,x_0)| \le \int_{{{\mathbb S}}^{1}}
%\sum_{j=1}^5|Q_{2,j}(t,x_0)| d\sigma \le C\lambda^{-1}(\lambda t)^{-\mu}
%\Big(
%\|f_2\|_{\dot{H}_x^{-\sigma}} + \sum_{\ell=0}^{2}\||x|^\ell f_2\|_{H_x^{\ell-3}}
%\Big)
%\end{eqnarray*}  %%%%%%%%%%%%%%%%%%%%%%%%%%%%%%%%%%%%%%%
%and
\begin{eqnarray} %%%%%%%%%%%%%%%%%%%%%%%%%%%%%%%%%%%%%%%
\lefteqn{|\nabla^{k}Q_{\lambda}^{(2.1)}(t,x_0)|}\nonumber\\
&\le& C
\lambda^{-1}\int_0^{t-\lambda^{-1}{b}}(1+\lambda (t-s))^{-1}
\sum_{\ell=0}^{2}\||x|^\ell f_2(s)\|_{H_{x}^{m+\ell-5}}\,ds
\nonumber\\
&\le& C
\lambda^{-2}\log(1+\lambda T)\sup_{0\le t\le T}
\left(\sum_{\ell=0}^{2}\||x|^\ell f_2(t)\|_{H_{x}^{m+\ell-5}}\right).
\label{y2}
\end{eqnarray}  %%%%%%%%%%%%%%%%%%%%%%%%%%%%%%%%%%%%%%%

\vskip2mm \noindent \underline{Case: $d=3$.}\,\,  Analogously we have
\begin{eqnarray*}%%%%%%%%%%%%%%%%%%%%%%%%%%%%%%%%%%%%%%%%%%%
|Q_{2,1}|&\lesssim&
(\lambda (t-s))^{-3/2} \\
& &\times\Big(\sum_{j=1}^{3}
\|\langle\eta\rangle^{-3}\partial_{j}\nabla\widehat{f_{2}}\eta\|_{L_{\eta}^{2}}+
\|\langle\eta\rangle^{-4}\nabla\widehat{f_{2}}\eta\|_{L_{\eta}^{2}}
\displaystyle +
\|\langle\eta\rangle^{-5}\widehat{f_{2}}\eta\|_{L_{\eta}^{2}}\Big), \\
|Q_{2,2}|+|Q_{2,3}|&\lesssim&
(\lambda (t-s))^{-3/2}\Big( \displaystyle
\|\langle\eta\rangle^{-4}\nabla \widehat{f_{2}}\eta\|_{L_{\eta}^{2}}+
\|\langle\eta\rangle^{-5}\widehat{f_{2}}\eta\|_{L_{\eta}^{2}} \Big),\\
|Q_{2,4}|+|Q_{2,5}|&\lesssim& (\lambda
(t-s))^{-3/2} \|\langle\eta\rangle^{-5}\widehat{f_{2}}\eta\|_{L_{\eta}^{2}}.
\end{eqnarray*}%%%%%%%%%%%%%%%%%%%%%%%%%%%%%%%%%%%%%%%%%%%%%%
Combining the above estimates, for $0\leq k \leq m$, we get
\begin{eqnarray} %%%%%%%%%%%%%%%%%%%%%%%%%%%%%%%%%%%%%%%
\lefteqn{|\nabla^{k}Q_{\lambda}^{(2,1)}(t,x_0)|}\nonumber\\
&\le& C
\lambda^{-1}\int_0^{t-\lambda^{-1}{b}}(1+\lambda (t-s))^{-3/2}
\sum_{\ell=0}^{2}\||x|^\ell f_2(s)\|_{H_{x}^{m+\ell-5}}\,ds
\nonumber\\
&\le& C
\lambda^{-2}\sup_{0\le t\le T}
\left(
\sum_{\ell=0}^{2}\||x|^\ell f_2(t)\|_{H_{x}^{m+\ell-5}}
\right).\label{y3}
\end{eqnarray}  %%%%%%%%%%%%%%%%%%%%%%%%%%%%%%%%%%%%%%%

On the other hand, the Sobolev embedding yields
\begin{eqnarray} %%%%%%%%%%%%%%%%%%%%%% %%%
&&\sup_{0\le t\le T}\, \sum_{\ell=0}^{2}
\||x|^{\ell}f_{2}(t)\|_{H_{x}^{m+\ell-6+[d/2]}} \nonumber\\
&\lesssim& \sup_{0\le t\le T} \big(1+\|n(t)\|_{H_{x}^{m+2
+[d/2] }}\big) \Big(\sum_{\ell=0}^{2}
\||x|^{\ell}E(t)\|_{H_{x}^{m+3+[d/2] -3\ell}}\Big)^{2}
%\nonumber\\
\label{y5}
\end{eqnarray}%%%%%%%%%%%%%%%%%%%%%% %%%%%
for $m\ge6-[d/2]$ if $d=1, 2, 3$. From (\ref{y10}), (\ref{y11}), (\ref{y1})
and (\ref{y5}), we obtain (\ref{y7}) for $d=1$. From (\ref{y10}),
(\ref{y11}), (\ref{y2}) and (\ref{y5}), we obtain (\ref{y7}) for $d=2$.
Further, from (\ref{y10}), (\ref{y11}), (\ref{y3}) and (\ref{y5}), we obtain
(\ref{y7}) for $d=3$. This completes the proof of Lemma \ref{3f}.
$\qquad\qed$

\vskip3mm \noindent {\bf Acknowledgments.} The authors express their
deep gratitude to Professor Kenji Nakanishi for valuable 
comments on the manuscript. This work was initiated while the third
author was visiting the Department of Mathematics at the UC, Santa
Barbara and NCKU at Taiwan whose hospitality was gratefully
acknowledged. The third author is partially supported by JSPS,
SYROVPABC and by MEXT, Grant-in-Aid for YS (A) 25707004. The first
two authors are partially supported by NSC and NCTS of Taiwan.

%------------------------%   references   %-----------------------------

%%%%%%%%%%%%%% %%%%%%%%%%%%%%%%


\begin{thebibliography}{30}

\bibitem{AA} Added, H. and Added, S., \textit{
Equations of Langmuir turbulence and nonlinear
Schr\"{o}dinger equation: smoothness and approximation}.
J. Funct. Anal. {\bf 79} (1988), 183--210.

\bibitem{BC}
Bourgain, J. and Colliander, J., \textit{On wellposedness of the Zakharov
system}. Internat. Math. Res. Notices {\bf 1996} (1996), 515--546.

\bibitem{BH}
Bejenaru, I. and Herr, S., \textit{ Convolutions of singular measures and
applications to the Zakharov system}. J. Funct. Anal. {\bf 261} (2011),
478--506.

\bibitem{BHHT}
Bejenaru, I., Herr, S., Holmer, J. and Tataru, D., \textit{On the 2D
Zakharov system with $L^{2}$- Schr\"{o}dinger data}. Nonlinearity {\bf
22} (2009), 1063--1089.

\bibitem{BGM} Br\'{e}zis H. and Gallou\"{e}t T.,
\textit{Nonlinear Schr\"{o}dinger evolution equations}.
Nonlinear Anal., Theory Methods Appl., {\bf 4} (1980), 677--681.

\bibitem{Caz}
Cazenave, T., \textit{``Semilinear Schr\"{o}dinger equations"}. Courant
Lecture Notes in Mathematics, {\bf 10}. American Mathematical Society
(2003).

\bibitem{CFW}
Chen, T.-J., Fang, Y.-F. and Wang, K.-H., \textit{Low regularity global
well-posedness for the quantum Zakharov system in 1D}. (2014)
preprint.

\bibitem{CHT}
Colliander, J., Holmer, J. and Tzirakis, N., \textit{Low regularity global
well-posedness for the Zakharov and Klein-Gordon-Schr\"{o}dinger
systems}, Trans. AMS. {\bf 360} (2008), 4619--4638.

\bibitem{EA}
El-Wakil, S.A. and Abdou, M.A., \textit{ New exact travelling wave
solutions of two nonlinear physical models}. Nonlinear Anal. {\bf 68}
(2008), 235--245.

\bibitem{GHPG} Garcia L.G., Haas F., de Oliveira L.P.L.
and Goedert J., \textit{Modified Zakharov equations
for plasmas with a quantum correction.}
Phys. Plasmas {\bf 12} 012302 (2005).

\bibitem{GTV}
Ginibre, J., Tsutsumi, Y. and Velo, G., \textit{On the Cauchy problem for
the Zakharov system}. J. Funct. Anal. {\bf 151} (1997), 384--436.

\bibitem{GM1}
Glangetas, L. and Merle, F., \textit{ Existence of self-similar blow-up
solutions for Zakharov equation in dimension two. I}. Comm. Math.
Phys. {\bf 160} (1994), 173--215.

\bibitem{GM2}
Glangetas, L. and Merle, F., \textit{ Concentration properties of blow-up
solutions and instability results for Zakharov equation in dimension two.
II}. Comm. Math. Phys. {\bf 160} (1994), 349--389.

\bibitem{GZG} Guo, Y. Zhang, J. and Guo, B.,
\textit{Global well-posedness and the classical limit
of the solution for the quantum Zakharov system}.
Z. Angew. Math. Phys. {\bf 64} (2013), 53--68.

\bibitem{GN}
Guo Z. and Nakanishi K., \textit{
Small energy scattering for the Zakharov system with radial symmetry},
Int. Math. Res. Not. IMRN {\bf 2014} (2014), 2327--2342.

\bibitem{GNW}
Guo, Z., Nakanishi, K. and Wang, S., \textit{ Global dynamics below the
ground state energy for the Zakharov system in the 3D radial case}. Adv.
Math. {\bf 238} (2013), 412--441.

\bibitem{Haas} Haas, F.,
\textit{``Quantum plasmas''}. Springer Series on Atomics,
Optical and Plasma Physics {\bf 65} (2011).

\bibitem{HS} Haas, F. and Shukla, P.K.,
\textit{Quantum and classical dynamics of Langmuir
wave packets}. Phys. Rev. E {\bf 79} 066402
(2009).

\bibitem{HPS}
Hani, Z., Pusateri, F., and Shatah, J., \textit{ Scattering for the Zakharov
system in 3 dimensions}. Comm. Math. Phys. {\bf 322} (2013), 731--753.

\bibitem{JLP}
Jiang, J.C., Lin, C.K., and Shao, S., \textit{On One Dimensional Quantum
Zakharov System}, (2013), To appear in Discrete and Continuous Dynamical 
System, Series A (2016). 

\bibitem{KT} Keel, M. and Tao, T. \textit{Endpoint Strichartz Estimate},
American Journal of Mathematics {\bf 120} (1998),  955-980.

\bibitem{KPV2}
Kenig C.E., Ponce G. and Vega L., \textit{
Oscillatory integrals and regularity
of dispersive equations}. Indiana Univ.math J.
{\bf 40} (1991), 33--69.

\bibitem{KPV}
Kenig, C.E., Ponce, G. and Vega, L., \textit{ On the Zakharov and
Zakharov-Schulman systems}, J. Funct. Anal. {\bf 127} (1995), 204--234.

\bibitem{K2}
Kishimoto, N., \textit{Local well-posedness for the Zakharov system on
the multidimensional torus}. J. Anal. Math. {\bf 119} (2013), 213--253.

\bibitem{KM}
Kishimoto, N. and Maeda, M., \textit{Construction of blow-up solutions
for Zakharov system on ${{\mathbb T}}^{2}$}, Ann. Inst. H. Poincar\'{e}
Anal. Non Lin\'{e}aire {\bf 30} (2013), 791--824.

\bibitem{Lions} P-L. Lions,
\textit{Mathematical Topics in Fluid Mechanics}, \emph{Oxford Lecture
Series in Mathematics and its Applications}, Vol. 1, Incompressible
Models, Oxford Science Publications, (1996).

\bibitem{MN1}
Masmoudi, N. and Nakanishi, K., \textit{ Energy convergence for singular
limits of Zakharov type systems}, Invent. Math. {\bf 172} (2008),
535--583.

\bibitem{MN2}
Masmoudi, N. and Nakanishi, K.,
\textit{From the Klein-Gordon-Zakharov
system to a singular nonlinear Schr\"odinger system},
Ann. I. H. Poincar\'e AN {\bf 27} (2010), 1073--1096.

\bibitem{OT1} Ozawa, T. and Tsutsumi, Y., \textit{
Existence and smoothing effect of solutions for the Zakharov equations},
Publ. Res. Inst. Math. Sci. {\bf 28} (1992), 329--361.

\bibitem{OT2} Ozawa, T. and Tsutsumi, Y., \textit{
The nonlinear Schr\"{o}dinger limit and the initial layer of the Zakharov
equations}, Differential Integral Equations {\bf 5} (1992), 721--745.

\bibitem{OT3} Ozawa, T. and Tsutsumi, Y., \textit{
Global existence and asymptotic behavior of solutions for the Zakharov
equations in three space dimensions}, Adv. Math. Sci. Appl. {\bf 3}
(1993/94), Special Issue, 301--334.

\bibitem{P} Pausader, B., \textit{
Global Well-Posedness for Energy Critical Fourth-Order Schr\"odinger
Equations in the Radial Case}, Dynamics of PDE {\bf 4}, No.3, (2007),
197-225.

\bibitem{S}
Shimomura, A., \textit{Scattering theory for Zakharov equations in
three-dimensional space with large data}, Commun. Contemp. Math.
{\bf 6} (2004), 881--899.

\bibitem{SW} Schochet, S. and Weinstein, M.I.,
\textit{The nonlinear Schr\"{o}dinger limit of the Zakharov equations
governing Langmuir turbulence}, Comm. Math. Phys. {\bf 106} (1986)
569--580.

\bibitem{U}
Ukai, S., \textit{The incompressible limit and the initial layer of the
compressible Euler equation}, J. Math. Kyoto Univ. {\bf 26} (1986),
323--331.

\bibitem{Z} Zakharov, V.E.,
\textit{Collapse of Langmuir waves}, Sov. Phys. JETP {\bf 35} (1972),
908--914.

\end{thebibliography}
\end{document}